\documentclass[final,hidelinks,onefignum,onetabnum]{siamart220329}

\usepackage{lipsum}
\usepackage{amsfonts}
\usepackage{graphicx}
\usepackage{epstopdf}
\usepackage{algorithmic}
\ifpdf
\DeclareGraphicsExtensions{.eps,.pdf,.png,.jpg}
\else
\DeclareGraphicsExtensions{.eps}
\fi
\usepackage{booktabs}
\usepackage{hyperref}

\usepackage{stmaryrd}
\usepackage{bm}
\usepackage{caption}
\usepackage{subcaption}
\usepackage{multirow}
\usepackage{enumitem}
\usepackage{tikz}
\usepackage{algorithm} 
\usepackage{algorithmic} 

\usepackage{geometry}
\newsiamthm{assumption}{Assumption}

\newsiamremark{remark}{Remark}
\newsiamremark{hypothesis}{Hypothesis}
\crefname{hypothesis}{Hypothesis}{Hypotheses}
\newsiamthm{claim}{Claim}
\newsiamremark{exmp}{Example}
\newtheorem{alg}{Algorithm}[section]

\headers{}{}

\title{A structure-preserving relaxation Crank-Nicolson finite element method for the Schr\"{o}dinger-Poisson equation
}

\author{Huini Liu\thanks{School of Mathematics and Computational Science, Xiangtan University, Xiangtan 411105, P.R.China (\email{liuhuini@smail.xtu.edu.cn})}
\and Nianyu Yi\thanks{Hunan Key Laboratory for Computation and Simulation in Science and Engineering; School of Mathematics and Computational Science, Xiangtan University, Xiangtan 411105, P.R.China (\email{yinianyu@xtu.edu.cn}).}
\and Peimeng Yin\thanks{Corresponding author. Department of Mathematical Sciences, The University of Texas at El Paso, El Paso, TX 79968, USA (\email{pyin@utep.edu}).}} 

\usepackage{amsopn}

\ifpdf
\hypersetup{
  pdftitle={},
  pdfauthor={}
}
\fi

\usepackage{comment}

\begin{document}

\maketitle

\begin{abstract}
In this paper, we propose a mass- and modified energy-conservative relaxation Crank-Nicolson finite element method for the Schr\"{o}dinger-Poisson equation. Utilizing only a single auxiliary variable, we simultaneously reformulate the distinct nonlinear terms present in both the Schr\"{o}dinger equation and the Poisson equation into their equivalent expressions, constructing a system equivalent to the original Schr\"{o}dinger-Poisson equation. Our proposed scheme, derived from this equivalent system, is implemented linearly, avoiding the need for iterative techniques to solve the nonlinear equation. Additionally, it is executed sequentially, eliminating the need to solve a coupled large linear system. We in turn rigorously derive the optimal error estimates for the proposed scheme, demonstrating second order accuracy in time and $(k+1)$th order accuracy in space when employing polynomials of degree up to $k$. Numerical experiments validate the accuracy and effectiveness of our method and emphasize its conservation properties over long-time simulations.  
\end{abstract}



\section{Introduction}\label{intro}

Consider the Schr\"{o}dinger-Poisson (SP) equation, also known as the Gross-Pitaevskii Poisson equation 
\cite{Cai2010Mean,Cotner2016Collisional,Shukla2006Formation}
\begin{subequations}\label{TargetEq}
\begin{align}
& \mathbf{i}u_{t} =-\Delta u+ \Phi u+ V(x)u +|u|^{2}u, \quad (x,t)\in~\Omega\times(0,T],\label{TargetEq1}\\
& -\Delta\Phi = \mu(|u|^{2}-c), \quad (x,t)\in \Omega \times [0,T], \label{TargetEq2}\\
& u(x,0) = u_{0}(x), \quad x\in \Omega. \label{TargetEq3}
\end{align}
\end{subequations}
Here, the symbol $\mathbf{i}=\sqrt{-1}$ represents the imaginary unit, $\Omega  \subset \mathbb{R}^{2}$ is a convex bounded domain, and $T > 0$ is the final time.
The complex-valued function $u(x, t)$ represents the single-particle wave function, while the real-valued function $\Phi(x, t)$ denotes the Poisson potential. Both functions satisfy the homogeneous Dirichlet boundary condition.
The nonlinear term $|u|^2u$ in the Schr\"{o}dinger equation is known as the self-repulsion, whereas the nonlinear term $|u|^2$ in the Poisson equation represents the charge density. 
The constant $\mu=\pm1$ is a rescaled physical constant, reflecting the nature of the underlying forcing: repulsive for $\mu > 0$ and attractive for $\mu < 0$. 
The parameter $c$ denotes a background charge of the particle independent of time $t$. 
$V(x)$ is a specified external potential, and $u_0(x)$ is the initial condition.

The Schr\"{o}dinger-Poisson equation was first introduced by Ruffini and Bonazzola \cite{Ruffini1696Systems} to study self-gravitating boson stars. 
Later on, it was explored in various fields of application, including quantum mechanics \cite{Cai2010Mean}, semiconductors
\cite{Markowich1990Semiconductor, Ringhofer2000Discrete}, plasma physics
\cite{Bertrand1980Classical, Shukla2006Formation, Shukla2011Nonlinear, Sakaguchi2020Gross}, optics \cite{Paredes2020From}. 
A significant body of literature is dedicated to the mathematical analysis and numerical approximation of the Schr\"{o}dinger-Poisson equation,
including the well-posedness \cite{Lange1995overview, Castella1997solutions, Arriola2001variational, Masaki2011Energy}.

In studies of Bose-Einstein condensates, boundary conditions for both $u$ and $\Phi$ in \eqref{TargetEq} typically vanish at infinity and are often scaled to bounded domains as homogeneous Dirichlet boundary conditions \cite{Cotner2016Collisional}.
For simplicity of presentation, we focus on the homogeneous Dirichlet boundary condition:
\begin{equation}\label{InitBound}
    u(x,t)=0 \quad \text{and}\quad \Phi(x, t)=0, \quad (x,t)\in \partial\Omega \times [0,T].
\end{equation}
However, various types of boundary conditions can be imposed on the SP equation \eqref{TargetEq}, including (homogeneous) Dirichlet boundary conditions \cite{Cotner2016Collisional, Arriola2001variational}, zero-flux (Neumann) boundary conditions \cite{Sakaguchi2020Gross}, and periodic boundary conditions \cite{Lange1995overview, Verma2021formation, Sakaguchi2020Gross}. More discussions about boundary conditions can be found in \cite{Lange1995overview, Lange1997Mixed} and the references therein. 
The method to be proposed later and its analysis are applicable to all these boundary conditions.
Under the homogeneous Dirichlet boundary conditions \eqref{InitBound}, the solution of the Schr\"{o}dinger-Poisson equation \eqref{TargetEq} preserves the mass conservation 
$$
M(t)=\int_{\Omega}|u|^{2}dx=M(0),
$$
and the energy conservation
\begin{equation}\label{inroeng}
E(t)=\int_{\Omega}\left(|\nabla u|^{2}+\frac{1}{2\mu}|\nabla\Phi|^{2}+ V(x)|u|^{2}+\frac{1}{2}|u|^{4}\right)dx=E(0),
\end{equation}
which are important invariant properties that are also desired at the discrete level. In literature, a modified energy is often selected as an alternative structure to the original energy, particularly in numerical methods that involve reformulating the Schr\"{o}dinger-Poisson equation \eqref{TargetEq} into an equivalent enlarged system \cite{Yi2022mass, Gong2022SAV}.
If the self-repulsion term $|u|^2u$ in (\ref{TargetEq}a) vanishes, several numerical methods have been proposed to handle the nonlinearity caused by the charge density $|u|^2$ in the Poisson equation, including the Strang splitting types of methods \cite{Auzinger2017Convergence, Lubich2008splitting}. To preserve the invariant properties at the discrete level, Ringhofer et al. introduced a Crank-Nicolson scheme \cite{Ringhofer2000Discrete} and employed the prediction-correction technique to handle the nonlinearity. 
An extension work of the Crank-Nicolson-type method was carried out by Ehrhardt et al. to develop an approximation for the spherically symmetric Schr\"{o}dinger-Poisson system \cite{Ehrhardt2006Fast}.  
A structure-preserving discontinuous Galerkin (DG) method proposed in \cite{Yi2022mass} also treated the nonlinear term implicitly, but an iterative technique was employed to handle the nonlinear term.
More recently, structure-preserving relaxation Crank-Nicolson types of methods were proposed for the nonlinear Schr\"{o}dinger equation \cite{Besse2004relaxation, Besse2021relaxation} and the Schr\"{o}dinger-Poisson equation \cite{Athanassoulis2023novel}.
The relaxation methods introduce an intermediate function to handle the nonlinearity and find solutions of Schr\"{o}dinger equation and Poisson equation at different time levels. Therefore, the corresponding schemes are linear.


For the nonlinear Schr\"{o}dinger-Poisson equation \eqref{TargetEq} that incorporates both the self-repulsion $|u|^2u$ and the charge density $|u|^2$, different techniques may be necessary to handle the two distinct nonlinear terms.
In addition, it is challenging to handle the two nonlinear terms while simultaneously conserving the invariant properties \cite{Wang2018splitting}.
A scalar auxiliary variable (SAV) Crank-Nicolson scheme was proposed in \cite{Gong2022SAV} that preserves both mass and modified energy properties.  
It is interesting to note that the SAV approach is only applied to the nonlinear term $|u|^2u$  while treating the nonlinear term $|u|^2$ simply implicitly.
Therefore, the method remains implicit and nonlinear, requiring iterative methods for convergence.
Another noteworthy DG method \cite{Yi2022mass} applies the relaxation techniques described in \cite{Besse2004relaxation} for the Schr\"{o}dinger equation but treats the nonlinear term $|u|^2$ in the Poisson equation implicitly. Therefore, iterative techniques are still needed to solve the coupled system formed by the Schr\"{o}dinger equation and the Poisson equation. 

It is natural to inquire whether it is possible to handle the nonlinear terms efficiently while conserving the invariant properties.
Motivated by effectiveness and the ability of the structure-preserving relaxation-type of schemes to preserve the invariants for the Schr\"{o}dinger equation and the general Schr\"{o}dinger-Poisson equation \cite{Athanassoulis2023novel,Besse2004relaxation, Besse2021relaxation}, we propose a linear and structure-preserving relaxation Crank-Nicolson finite element method tailored for solving the nonlinear Schr\"{o}dinger-Poisson equation \eqref{TargetEq}.
More specifically, we introduce only one auxiliary variable to reformulate two different nonlinear terms in two equations simultaneously: the self-repulsion term $|u|^2u$ in the Schr\"{o}dinger equation (\ref{TargetEq}a), and the charge density $|u|^2$ in the Poisson equation (\ref{TargetEq}b). This transforms the Schrödinger-Poisson equation \eqref{TargetEq} into an equivalent system, facilitating its discretization into a linear fully discrete finite element scheme. This approach conserves both mass and modified energy, while also allowing for a linear implementation without the need for iterative techniques.
To the best of our knowledge, the approach that introduces only one auxiliary variable to simultaneously reformulate different nonlinear terms in two distinct equations in a system, as described, has not been explored in the literature for the Schr\"{o}dinger-Poisson equation \eqref{TargetEq}.


Error analysis of the numerical methods for the Schrödinger-Poisson equations is crucial for assessing their stability and accuracy, but much attention has been given to optimal error analysis for the Schrödinger-Poisson equations without the self-repulsion term.
Lubich \cite{Lubich2008splitting} pioneered the error analysis of the Strang-type splitting method in the semi-discretization system. 
Auzinger et al. \cite{Auzinger2017Convergence} analyzed the convergence analysis for the fully discrete scheme for the Schr\"{o}dinger-Poisson equation by using the splitting finite element method. Later on, Zhang \cite{Zhang2013Optimal} studied the optimal error estimates of the finite difference method under proper regularity assumptions. 
The optimal $L^{2}$ error estimate of semi-discrete conservative DG scheme was also proved in \cite{Yi2022mass}.
However, limited research on error analysis has been established for numerical methods incorporating the nonlinear self-repulsion term. Gong et al. \cite{Gong2022SAV} established unconditional energy stability and performed convergence analysis for the SAV Crank-Nicolson spectral method. 

In this work, we rigorously derive optimal a priori error estimates for the relaxation Crank-Nicolson finite element method using the method of induction. 
Various tools have been introduced and developed to obtain the desired results, such as the uniform boundedness of the finite element approximations, and the dependence of the errors between different equations. 
Specially, the $L^2$ error of the solution in the Poisson equation is bounded by the $L^2$ error of auxiliary variable and an optimal spatial error bound.  
As a result, we obtain second order accuracy in time and $(k+1)$th order accuracy in space when employing polynomials of degrees no more than $k$.
To the best of our knowledge, there are currently no rigorous convergence results in the literature for relaxation Crank-Nicolson types of methods for the Schr\"{o}dinger-Poisson equation. 
The analysis technique developed in this work can be extended to other similar numerical methods, offering a broader applicability.
An extension of the error analysis for the structure-preserving relaxation Crank-Nicolson finite element method to the Schr\"{o}dinger-Poisson equation, without the self-repulsion term $|u|^2u$ in (\ref{TargetEq}a), was also provided.

The contributions, innovations, and significance of this work include:
\begin{itemize}[leftmargin=*]
    \item Different from the existing methods that use various techniques to handle the two distinct nonlinear terms in the Schrödinger-Poisson equation \eqref{TargetEq}, we employ only one technique, namely the relaxation method, for both nonlinear terms. Consequently, the proposed method is easy to implement.
    \item Though we use only one technique to handle the two different nonlinear terms, we prove that the proposed method preserves both mass and modified energy.
    \item The proposed method is implemented linearly without resorting to any iterative techniques and sequentially without the need to solve a coupled system. Therefore, it is computationally efficient and cheap.
    \item We derived the optimal error estimates for the proposed method, obtaining second-order accuracy in time and $(k+1)$th order accuracy in space for the $L^2$ errors when applying polynomials with a maximal degree $k$.
    \item We conduct numerical examples to verify the performance of the proposed method, including accuracy tests, conservation verification, and comparisons with existing results.
\end{itemize}

The organization of this paper is as follows. In Section \ref{secScheme}, we present the relaxation Crank-Nicolson finite element method for the Schr\"{o}dinger-Poisson equation, and we demonstrate the structure-preserving properties of both the continuous problem and the fully discrete scheme. In Section \ref{secErrorEstimate}, we establish the optimal error estimates in $L^{2}$ norm for the solutions of a fully discrete system, comprising second-order accuracy in time and $(k+1)$th order accuracy in space. We further extend the convergence results to the relaxation Crank-Nicolson scheme \cite{Athanassoulis2023novel} in Section \ref{secExtension}. In Section \ref{secNmecical}, some numerical experiments are carried out to validate the theoretical analysis and verify the performance of the proposed conservative method.

We employ $W^{m,p}(\Omega,\mathbb{R})$ and $W^{m,p}(\Omega,\mathbb{C})$ to denote real-valued and complex-valued Sobolev spaces, respectively. For brevity, we use $H^{m}(\Omega)$ for $W^{m,2}(\Omega,\mathbb{R})$ and $\mathbf{H}^{m}(\Omega)$ for $W^{m,2}(\Omega,\mathbb{C})$, with norms denoted by $\|\cdot\|_{m}$ and semi-norms by $|\cdot|_{m}$. When $m=0$, $\|\cdot\|$ represents the $L^{2}$ norm of either $L^{2}(\Omega)$ or $\mathbf{L}^{2}(\Omega)$.
Unless explicitly stated otherwise, the constants denoted by $C$, possibly accompanied by a suitable subscript, represent generic positive constants that are independent of $\tau$, $h$, $n$, and $N$, but may depend on final time $T$ and the regularity of exact solutions $u$ and $\Phi$.

\section{The relaxation Crank-Nicolson Finite Element Method}\label{secScheme}

In the following presentation, the inner product and norm of the standard complex-valued Hilbert space $\mathbf{L}^{2}(\Omega)$ are expressed as $\langle\cdot,\cdot\rangle$ and $\|\cdot\|$, respectively, 
$$\langle u, v\rangle:=\int_{\Omega}uv^{\ast}dx \quad \text{and}\quad \|u\|=\sqrt{\langle u,u\rangle},$$
where $v^{\ast}$ denotes the complex conjugate of $v$. Similarly, the inner product and norm of the
real-valued Hilbert space $L^{2}(\Omega)$ are defined by
$$(u, v):=\int_{\Omega}uvdx \quad \text{and}\quad \|u\|=\sqrt{(u, u)}.$$
Then the weak formulation of problem \eqref{TargetEq} reads as: find $u\in C^{1}([0,T],\mathbf{H}_{0}^{1}(\Omega))$ and
$\Phi\in H_{0}^{1}(\Omega)$, 
\begin{align}
&\mathbf{i}\left\langle u_{t},\omega\right\rangle=A_{0}\left(u, \omega\right)+\left\langle\Phi u,\omega\right\rangle
+\left\langle V(x)u,\omega\right\rangle +\left\langle|u|^{2}u, \omega\right\rangle, \quad \forall \omega\in \mathbf{H}_{0}^{1}(\Omega),\label{EqWeaka}\\
& A_{1}\left(\Phi,\chi\right) = \mu\left(|u|^{2}-c,\chi\right), \quad \forall \chi\in H_{0}^{1}(\Omega), \label{EqWeakb}
\end{align} 
where the bilinear forms $A_{0}(\cdot, \cdot)$ and $A_{1}(\cdot, \cdot)$ are defined as follows 
\begin{align}
A_{0}(\omega,v)=\langle\nabla \omega, \nabla v\rangle,  \quad \forall \omega, v \in \mathbf{H}_{0}^{1}(\Omega), \\ 
A_{1}(\phi,\chi)=(\nabla \phi, \nabla \chi), \quad \forall \phi, \chi\in H_{0}^{1}(\Omega), 
\end{align}
and they both satisfy the coercivity and continuity properties, namely, there exist constants $\gamma_{1}>0$ and $\gamma_{2}>0$ such that
\begin{equation}\label{coercon}
A_{j}(v, v)\geq \gamma_{1}\|v\|_{1}^{2}, \quad A_{j}(\omega, v)\leq \gamma_{2}\| \omega\|_{1}\| v\|_{1},  \quad j = 0, 1, 
\end{equation}
for any $\omega, v\in H_{0}^{1}(\Omega)$ or $\mathbf{H}_{0}^{1}(\Omega)$.

\subsection{Mass and conservation properties}
We begin with the review of the continuous mass and energy conservation for the Schr\"{o}dinger-Poisson equation \eqref{TargetEq}. Then, we propose a finite element method that conserves these properties.

The Schr\"{o}dinger-Poisson equation \eqref{TargetEq} is nonlinear, containing two nonlinear terms: the self-repulsion term $|u|^2u$ in the Schr\"{o}dinger equation (\ref{TargetEq}a), and the charge density $|u|^2$ in the Poisson equation (\ref{TargetEq}b). Observing that two nonlinearities share a common factor, we introduce a real auxiliary variable
$\Psi=|u|^{2}$. The Schr\"{o}dinger-Poisson equation \eqref{TargetEq} can then be equivalently written as
\begin{equation}\label{TargetEqRe}
\left\{
\begin{aligned}
&\Psi=|u|^{2}, \\
&\mathbf{i}u_{t}=-\Delta u+ \Phi u+ V(x)u +\Psi u,\\
&-\Delta\Phi=\mu(\Psi-c),
\end{aligned}
\right.
\end{equation}
whose weak formulation is to find $u\in C^{1}([0,T],\mathbf{H}_{0}^{1}(\Omega))$ and $\Psi,\Phi\in H_{0}^{1}(\Omega)$ such that\begin{subequations}\label{TargetEqWeak}
\begin{align}
&\left(\Psi, v\right)=\left(|u|^{2}, v\right),\quad \forall v\in H_{0}^{1}(\Omega)\label{TargetEqWeak0}\\
&\mathbf{i}\left\langle u_{t}, \omega\right\rangle
=A_{0}\left(u, \omega\right) 
+\left\langle \Phi u, \omega\right\rangle+\left\langle V(x) u, \omega\right\rangle
+\left\langle \Psi u, \omega\right\rangle, \quad \forall \omega\in \mathbf{H}_{0}^{1}(\Omega),\label{TargetEqWeak1}\\
&A_{1}\left(\Phi,\chi\right)= \mu\left(\Psi-c, \chi\right)
, \quad \forall \chi\in H_{0}^{1}(\Omega).\label{TargetEqWeak2}
\end{align}
\end{subequations}

Similar to \eqref{TargetEq}, the following invariants are preserved for the new Schr\"{o}dinger-Poisson system:
\begin{align}
\text{mass conservation} \quad M(t)&=\int_{\Omega}|u|^{2}dx=M(0), \label{MassConservationRe} \\
\text{energy conservation} \quad E(t)&=\int_{\Omega}\left(|\nabla u|^{2}+\frac{1}{2\mu}|\nabla\Phi|^{2}+ V(x)|u|^{2}+\frac{1}{2}\Psi^{2}\right)dx=E(0).\label{EnergyConservationRe}
\end{align}
Indeed, by substituting $\omega=u$ in \eqref{TargetEqWeak1}, we obtain
\begin{align}\label{MassProofEq1}
\mathbf{i}\left\langle u_{t}, u\right\rangle
&= A_{0}(u,u )
+\left\langle \Phi u, u\right\rangle+\left\langle V(x) u, u\right\rangle
+\left\langle \Psi u, u\right\rangle.
\end{align}
Taking the imaginary part of \eqref{MassProofEq1} yields
\begin{equation}\label{MassProofEq2}
\frac{1}{2}\frac{d}{dt}\int_{\Omega} |u|^{2}dx = 0,
\end{equation}
which proves the mass conservation \eqref{MassConservationRe}.

On the other hand, by taking $\omega=u_{t}$ in \eqref{TargetEqWeak1}, it holds
\begin{equation}\label{EnergyProofEq1}
\mathbf{i}\left\langle u_{t}, u_{t}\right\rangle
= A_{0}(u, u_{t})
+\left\langle \Phi u, u_{t}\right\rangle + \left\langle V(x) u, u_{t}\right\rangle +\left\langle \Psi u, u_{t}\right\rangle.
\end{equation}
The real part of \eqref{EnergyProofEq1} yields
\begin{equation}\label{EnergyProofEq2a}
\frac{d}{dt}\int_{\Omega}\left|\nabla u\right|^{2}dx
+\int_{\Omega}\Phi\frac{d}{dt}\left|u\right|^{2}dx
+\frac{d}{dt}\int_{\Omega}V(x)\left|u\right|^{2}dx
+\int_{\Omega} \Psi \frac{d}{dt}|u|^{2}dx=0.
\end{equation}
By taking $v = \Phi$ in \eqref{TargetEqWeak0}$_{t}$, which is a resulting equation from differentiation of \eqref{TargetEqWeak0} in t, it follows 
\begin{equation}\label{EnergyProofEq3a}
\int_{\Omega}\Phi\frac{d}{dt}\left|u\right|^{2}dx=\int_{\Omega}\Phi\frac{d}{dt}(\Psi-c)dx.
\end{equation}
Similarly, by taking $\chi = \Phi$ in \eqref{TargetEqWeak2}$_{t}$, the second term in \eqref{EnergyProofEq2a} can be rewritten as
\begin{equation}\label{EnergyProofEq3}
\int_{\Omega}\Phi\frac{d}{dt}\left|u\right|^{2}dx
=\frac{1}{\mu}\int_{\Omega}\nabla\Phi_{t}\cdot\nabla\Phi dx.
\end{equation}
Setting $v=\Psi$ in \eqref{TargetEqWeak0}$_{t}$, the last term in \eqref{EnergyProofEq2a} can be rewritten as 
\begin{equation}\label{EnergyProofEq4}
\frac{1}{2}\frac{d}{dt}\int_{\Omega} \Psi^{2}dx=\int_{\Omega} \Psi \frac{d}{dt}|u|^{2} dx.
\end{equation}
Therefore, \eqref{EnergyProofEq2a} reduces to
\begin{equation}
\frac{d}{dt}\left(\int_{\Omega}\left|\nabla u\right|^{2}dx
+ \frac{1}{2\mu}\int_{\Omega}\left|\nabla \Phi\right|^{2}dx
+\int_{\Omega}V(x)\left|u\right|^{2}dx
+\frac{1}{2}\int_{\Omega}\Psi^{2}dx\right)=0.
\end{equation}
Hence, the energy conservation \eqref{EnergyConservationRe} holds.

\subsection{Fully discrete scheme}
To preserve the properties mentioned above at the discrete level, we investigate a relaxation Crank-Nicolson finite element method in this subsection.

Let $\mathcal{T}_{h}$ be a shape-regular and quasi-uniform triangulation of $\Omega  \subset \mathbb{R}^{2}$, $K \in \mathcal{T}_{h}$ be an element, and $h:=\max_{K\in \mathcal{T}_{h}}h_{K}$ be the mesh size with $h_{K}$ being the diameter of $K$. 
We denote the real-valued finite element space $V_{h}$ by
\begin{equation}
V_{h}=\left\{v\in C(\Omega): v\in \mathbb{P}^{k}(K),\quad~\forall K\in\mathcal{T}_{h}\right\},
\end{equation}
where $\mathbb{P}^{k}$ is the space of real-valued polynomials of degree up to the $k~(k\geq 1)$.
Similarly, the complex-valued finite element space $V_{h}^{c}$, associated with the triangulation $\mathcal{T}_{h}$, is denoted by
\begin{equation}
V_{h}^{c}=\left\{v\in C(\Omega): v\in \mathbb{Q}^{k}(K),\quad~\forall K\in\mathcal{T}_{h}\right\},
\end{equation}
where $\mathbb{Q}^{k}$ denotes the space of complex-valued polynomials of degree up to $k~(k\geq 1)$ in space.

Denote by $\{t_{n} \ | \ t_{n}=n\tau, 0\leq n\leq N\}$ a uniform partition of time interval $[0,T]$ with time step size $\tau=T/N$, where $N$ is a positive integer. 
We also introduce $t_{n-1/2} = (t_{n}+t_{n-1})/2 = (n-\frac{1}{2})\tau$. 
For any function $\varphi(x,t)$ and $n\geq 0$, we denote $\varphi_h^{n-\theta} \in$ $V_h$ or $V_h^c$ as an approximation of $\varphi(x,t_{n-\theta})$, where $\theta = 0, \frac{1}{2}$.

For a sequence of functions $\{\varphi^{n}\}_{n=0}^{N}$, we define the operators
\begin{equation}
D_{\tau}\varphi^{n+1}:=\frac{\varphi^{n+1}-\varphi^{n}}{\tau}, \quad \overline{\varphi}^{n+1/2}:=\frac{\varphi^{n+1}+\varphi^{n}}{2}.
\end{equation}
The relaxation Crank-Nicolson method introduces an intermediate function, and solves the intermediate function and the solution of the Schr\"{o}dinger equation at different time levels. Therefore, the corresponding scheme can be implemented linearly.
For the linearity of the scheme when coupled with the Poisson equation, we further solve the Poisson equation in the same time level as the intermediate function.
More specifically, for given $(\Psi_{h}^{n-1/2}, u_{h}^{n},\Phi_{h}^{n-1/2})\in V_{h}\times V_{h}^{c}\times V_{h}$, the relaxation Crank-Nicolson finite element scheme, derived from \eqref{TargetEqRe} or its weak formualation \eqref{TargetEqWeak}, is to find $(\Psi_{h}^{n+1/2}, u_{h}^{n+1},\Phi_{h}^{n+1/2})\in V_{h}\times V_{h}^{c}\times V_{h}$ such that
\begin{subequations}\label{fullyVdiscrete}
\begin{align}
&(\Psi_{h}^{n+1/2}+\Psi_{h}^{n-1/2}, v_{h})
=(2|u_{h}^{n}|^{2},v_{h}), \quad \forall v_{h}\in V_{h},\label{fullyVdiscrete1}\\
&\mathbf{i}\langle D_{\tau}u_{h}^{n+1},\omega_{h}\rangle=A_{0}(\overline{u}_{h}^{n+1/2},\omega_{h})
+\langle(\Phi_{h}^{n+1/2}+V(x)+\Psi_{h}^{n+1/2})
\overline{u}_{h}^{n+1/2}, \omega_{h}\rangle, \quad \forall \omega_{h}\in V_{h}^{c},\label{fullyVdiscrete2}\\
&A_{1}(\Phi_{h}^{n+1/2}, \chi_{h})=\mu(\Psi_{h}^{n+1/2}-c,\chi_{h}), \quad \forall \chi_{h}\in V_{h},\label{fullyVdiscrete3}
\end{align}
\end{subequations}
where the initial data $u_{h}^{0}=\Pi_{h} u_{0}$ and $\Psi_{h}^{-1/2}=\Pi_{h}|u_{h}^{0}|^{2}$. Here, $\Pi_{h}: H^{1}(\Omega)\rightarrow V_{h}$ is the nodal interpolation operator.
To compute the initial energy, we need $\Phi_{h}^{-1/2}\in V_h$, which is obtained by
$$
A_{1}(\Phi_{h}^{-1/2},\chi_{h})=\mu(\Psi_{h}^{-1/2}-c,\chi_{h}), \quad \forall \chi_{h} \in V_h.
$$

\begin{lemma}
Given $(\Psi_{h}^{n-1/2}, u_{h}^{n},\Phi_{h}^{n-1/2})\in V_{h}\times V_{h}^{c}\times V_{h}$ and $\tau>0$, the relaxation Crank-Nicolson finite element scheme \eqref{fullyVdiscrete} admits a unique solution $(\Psi_{h}^{n+1/2}, u_{h}^{n+1},\Phi_{h}^{n+1/2})\in V_{h}\times V_{h}^{c}\times V_{h}$.
\end{lemma}
\begin{proof}
The scheme \eqref{fullyVdiscrete} is a finite dimensional system, whose existence is equivalent to its uniqueness, thus we only need to show its uniqueness. Assume that \eqref{fullyVdiscrete} has two possible solutions and their difference is denoted by $(\delta \Psi_{h}^{n+1/2}, \delta u_{h}^{n+1}, \delta\Phi_{h}^{n+1/2})$, then it satisfies
\begin{subequations}\label{fullyVdiscrete-}
\begin{align}
&(\delta \Psi_{h}^{n+1/2}, v_{h})
=0, \quad \forall v_{h}\in V_{h},\label{fullyVdiscrete1+}\\
&\mathbf{i}\langle \delta u_{h}^{n+1}/\tau,\omega_{h}\rangle=\frac{1}{2}A_{0}(\delta{u}_{h}^{n+1},\omega_{h})
+\frac{1}{2}\langle(\Phi_{h}^{n+1/2}+V(x)+\Psi_{h}^{n+1/2})
\delta{u}_{h}^{n+1}, \omega_{h}\rangle, \quad \forall \omega_{h}\in V_{h}^{c},\label{fullyVdiscrete2+}\\
&A_{1}(\delta\Phi_{h}^{n+1/2}, \chi_{h})=\mu(\delta\Psi_{h}^{n+1/2},\chi_{h}), \quad \forall \chi_{h}\in V_{h},\label{fullyVdiscrete3+}
\end{align}
\end{subequations}
Taking $v_h = \delta \Psi_{h}^{n+1/2}$ in \eqref{fullyVdiscrete1+} gives $\|\delta \Psi_{h}^{n+1/2}\| = 0$, namely $\delta \Psi_{h}^{n+1/2} =0$. Then \eqref{fullyVdiscrete3+} gives
\begin{equation}\label{Adphi}
A_{1}(\delta\Phi_{h}^{n+1/2}, \chi_{h}) =0,  \quad \forall \chi_{h}\in V_{h}.
\end{equation}
By taking $\chi_h = \delta\Phi_{h}^{n+1/2}$ in \eqref{Adphi} and applying \eqref{coercon}, it follows 
$$
\|\delta\Phi_{h}^{n+1/2}\|_1\leq 0,
$$
which implies $\delta\Phi_{h}^{n+1/2}=0$.
Finally, taking $\omega_h = \tau\delta u_h^{n+1}$ in \eqref{fullyVdiscrete2+} yields
$$
\mathbf{i}\| \delta u_{h}^{n+1}\|^2=\frac{\tau}{2}A_{0}(\delta{u}_{h}^{n+1},\delta{u}_{h}^{n+1})
+\frac{\tau}{2}\langle(\Phi_{h}^{n+1/2}+V(x)+\Psi_{h}^{n+1/2})
\delta{u}_{h}^{n+1}, \delta{u}_{h}^{n+1}\rangle,
$$
and the imaginary part gives $\|\delta u_{h}^{n+1}\|=0$, or equivalently, $\delta u_{h}^{n+1} = 0$.
Thus, the conclusion holds.
\end{proof}

By solving the intermediate function, the Poisson equation, and the solutions of the Schr\"{o}dinger equation at different time levels,
namely the intermediate function and the Poisson solution at $t_{n+1/2}$, and the Schr\"{o}dinger at $t_{n+1}$, the relaxation Crank-Nicolson finite element method \eqref{fullyVdiscrete} can be implemented in the following algorithm.

\begin{alg}\label{rCNalg}
The relaxation Crank-Nicolson finite element method \eqref{fullyVdiscrete} is solved sequentially and linearly as follows.
\begin{itemize}
\item Solve $\Psi_{h}^{n+1/2} \in V_{h}$ from \eqref{fullyVdiscrete1}.
\item Solve $\Phi_{h}^{n+1/2} \in V_{h}$ from \eqref{fullyVdiscrete3}.
\item Solve $u_{h}^{n+1} \in V_h^c$ from \eqref{fullyVdiscrete2}.
\end{itemize}
\end{alg}
\begin{remark}
The proposed relaxation Crank-Nicolson finite element method \eqref{fullyVdiscrete} is linear without resorting to any interaction techniques, \Cref{rCNalg} additionally implies that it does not require solving a couple system. 
\end{remark}

\subsection{Structure-preserving properties}
From the literature, to conserve the invariant properties is challenging to solve the Schr\"{o}dinger-Poisson equation \eqref{TargetEq}.
Next, we explore the conservation properties of the proposed relaxation Crank-Nicolson finite element scheme \eqref{fullyVdiscrete} and obtain the following statement.
\begin{lemma}\label{CoffiConservationLem}
For any $\tau>0$, the relaxation Crank-Nicolson finite element method \eqref{fullyVdiscrete} satisfies the discrete conservation for both mass and modified energy with $0\leq n\leq N-1$, respectively
\begin{align}
&M_{h}^{n+1}= M_{h}^{0}, \label{DicreteMassV}\\
&E_{h}^{n+1}= E_{h}^{0}, \label{DicreteEnergyV}
\end{align}
where the mass
$M_{h}^{n+1}:=\int_{\Omega}\left|u_{h}^{n+1}\right|^{2}dx$,
and the modified energy
\begin{equation*}
\begin{split}
E_{h}^{n+1}&:=A_{0}(u_{h}^{n+1}, u_{h}^{n+1})+\frac{1}{2\mu}A_{1}(\Phi_{h}^{n+3/2},\Phi_{h}^{n+1/2})
+\int_{\Omega}V(x)|u_{h}^{n+1}|^{2}dx
+\frac{1}{2}\int_{\Omega}\Psi_{h}^{n+3/2}\Psi_{h}^{n+1/2}dx.
\end{split}
\end{equation*}
\end{lemma}

\begin{proof}
Taking $\omega_{h}=\overline{u}_{h}^{n+1/2}$ in \eqref{fullyVdiscrete2} yields
\begin{equation}\label{MassEq}
\mathbf{i}\langle D_{\tau}u_{h}^{n+1},\overline{u}_{h}^{n+1/2}\rangle=A_{0}( \overline{u}_{h}^{n+1/2},\overline{u}_{h}^{n+1/2})
+\langle(\Phi_{h}^{n+1/2}+V(x)+\Psi_{h}^{n+1/2})
\overline{u}_{h}^{n+1/2},\overline{u}_{h}^{n+1/2}\rangle.
\end{equation}
Then the imaginary part of $\eqref{MassEq}$ gives
\begin{equation}
\|u_{h}^{n+1}\|^{2}-\|u_{h}^{n}\|^{2}=0,
\end{equation}
which implies the conservation of the mass \eqref{DicreteMassV}. 

Next, taking $\omega_{h}=D_{\tau}u_{h}^{n+1}$ in \eqref{fullyVdiscrete2} gives
\begin{equation}\label{EnergyEq}
\mathbf{i}\langle D_{\tau}u_{h}^{n+1}, D_{\tau}u_{h}^{n+1}\rangle=A_{0}( \overline{u}_{h}^{n+1/2}, D_{\tau}u_{h}^{n+1})
+\langle(\Phi_{h}^{n+1/2}+V(x)+\Psi_{h}^{n+1/2})
\overline{u}_{h}^{n+1/2}, D_{\tau}u_{h}^{n+1}\rangle.
\end{equation}
The real part of $\eqref{EnergyEq}$ implies
\begin{equation}\label{EnergyRew}
\left[A_{0}(u_{h}^{n+1},u_{h}^{n+1})-A_{0}(u_{h}^{n},u_{h}^{n})\right]
+\int_{\Omega}(\Phi_{h}^{n+1/2}
+V(x)+\Psi_{h}^{n+1/2})(|u_{h}^{n+1}|^{2}-|u_{h}^{n}|^{2})dx=0.
\end{equation}
Then we proceed to estimate the terms in \eqref{EnergyRew}. Upon calculation, 
\begin{equation}\label{EnConserPhiEq}
\begin{split}
&\int_{\Omega}\Phi_{h}^{n+1/2}\left(|u_{h}^{n+1}|^{2}-|u_{h}^{n}|^{2}\right)dx\\
&\quad=\int_{\Omega}\Phi_{h}^{n+1/2}\left(\frac{\Psi_{h}^{n+3/2}
+\Psi_{h}^{n+1/2}}{2}-\frac{\Psi_{h}^{n+1/2}
+\Psi_{h}^{n-1/2}}{2}\right)dx \qquad \text{by \eqref{fullyVdiscrete1} }\\
&\quad=\frac{1}{2}\int_{\Omega}\Phi_{h}^{n+1/2}(\Psi_{h}^{n+3/2}-c)
-\Phi_{h}^{n+1/2}(\Psi_{h}^{n-1/2}-c)dx\\
&\quad=\frac{1}{2\mu} A_{1}\left(\Phi_{h}^{n+3/2}, \Phi_{h}^{n+1/2} \right) -\frac{1}{2\mu}A_{1}\left(\Phi_{h}^{n+1/2}, \Phi_{h}^{n-1/2} \right) \qquad \text{by \eqref{fullyVdiscrete3} }.
\end{split}
\end{equation}
Similarly, by \eqref{fullyVdiscrete1}, it holds  
\begin{equation}\label{EnConservationPsiEq}
\begin{split}
\int_{\Omega}\Psi_{h}^{n+1/2}\left(| u_{h}^{n+1}|^{2}-| u_{h}^{n}|^{2}\right)dx=\frac{1}{2}\int_{\Omega}\left(\Psi_{h}^{n+3/2}\Psi_{h}^{n+1/2}-\Psi_{h}^{n+1/2}\Psi_{h}^{n-1/2}\right)dx.
\end{split}
\end{equation}
Plugging \eqref{EnConserPhiEq} and \eqref{EnConservationPsiEq} into \eqref{EnergyRew} and regrouping give the discrete energy conservation \eqref{DicreteEnergyV}.
\end{proof}

\begin{remark}
Handling the two nonlinear terms in the Schr\"{o}dinger-Poisson equation \eqref{TargetEq} while conserving the original energy at the discrete level remains a challenging task. The discrete energy has only been numerically verified for the splitting Chebyshev collocation method proposed in \cite{Wang2018splitting}, whereas the iterative methods in \cite{Yi2022mass, Gong2022SAV} conserve modified rather than original energies.
Although our proposed method also conserves a modified energy, it achieves this with much higher efficiency. 
\end{remark}

\section{Error estimates for the fully discrete system}\label{secErrorEstimate}
The main objective of this section is to establish the optimal error estimates of the relaxation Crank-Nicolson finite element method \eqref{fullyVdiscrete} for the Schr\"{o}dinger-Poisson equation \eqref{TargetEq}. To begin with, we review some useful results.

Recall that $\Pi_{h}: H^{1}(\Omega)\rightarrow V_{h}$ be the nodal interpolation operator. By the classical finite element approximation theory \cite{Brenner2008The}, it follows
\begin{equation}\label{interpol}
\left\|v-\Pi_{h}v\right\|
+h\left\|\nabla\left(v-\Pi_{h}v\right)\right\|+h\|v-\Pi_{h}v\|_{\infty}\leq Ch^{k+1}\|v\|_{k+1},\quad \forall v\in  H^{k+1}(\Omega).
\end{equation}
We also define the Ritz projection operator $\mathrm{R}_{h}: H_{0}^{1}(\Omega)\rightarrow V_{h}$, which satisfies 
\begin{equation}\label{Ritz}
\big(\nabla(v-\mathrm{R}_{h}v),\nabla\omega \big)=0, \quad  \forall \omega\in V_{h},
\end{equation}
and holds the projection error estimate
\begin{equation}\label{Ritz2}
\left\|v-\mathrm{R}_{h}v\right\|
+h\left\|\nabla\left(v-\mathrm{R}_{h}v\right)\right\|\leq Ch^{k+1}\|v\|_{k+1}, \quad \forall v\in H_{0}^{1}(\Omega)\cap H^{k+1}(\Omega).
\end{equation}
The following inverse inequality \cite{Ciarlet1978Finite} will be widely used in the analysis,
\begin{equation} \label{InverIeq}
\|v\|_{\infty}\leq Ch^{-1}\|v\|, \quad \forall v\in V_{h}.
\end{equation}
In addition, we also need the following result.
\begin{lemma} \label{InftyLem}
For the Ritz projection defined in \eqref{Ritz}, it holds for any $k\geq 1$,
\begin{equation}\label{Infbdd}
\|\mathrm{R}_{h}v\|_{\infty} \leq C, \quad \forall v\in H_{0}^{1}(\Omega)\cap H^{k+1}(\Omega),
\end{equation}
where $C$ depends on $\|v\|_{k+1}$ and $\|v\|_\infty$, independent of $h$.
\end{lemma}
\begin{proof}
By the embedding theorem, it follows $H^2(\Omega) \subset L^\infty(\Omega)$. Then $v\in H^{k+1}(\Omega)$ implies $v \in L^\infty(\Omega)$.
Let $\Pi_h v$ be the nodal interpolation of $v$. By \eqref{interpol}, \eqref{Ritz2}, and the triangle inequality,
\begin{align*}
\|v-\mathrm{R}_{h}v\|_{\infty} \leq & \|v-\Pi_h v \|_{\infty}  + \|\Pi_h v -\mathrm{R}_{h}v \|_{\infty} \leq  \|v-\Pi_h v \|_{\infty}  + Ch^{-1}\|\Pi_h v  -\mathrm{R}_{h}v \|\\
\leq &  \|v-\Pi_h v \|_{\infty}  + Ch^{-1}\left( \|v- \Pi_h v \| + \|v -\mathrm{R}_{h}v \| \right) \leq  Ch^k\|v\|_{k+1}.
\end{align*}
Therefore, applying the triangle inequality gives
\begin{align*}
\|\mathrm{R}_{h}v\|_{\infty} \leq \|v-\mathrm{R}_{h}v\|_{\infty} + \|v\|_{\infty} \leq C.
\end{align*}
\end{proof}
\begin{remark}
Specially, the result in \Cref{InftyLem} holds for any $v \in H_0^1(\Omega) \cap H^{s+1}(\Omega)$ with $s>0$. 
The projection errors \eqref{interpol}, \eqref{Ritz2}, \eqref{InverIeq}, and the bound \eqref{Infbdd} also hold for functions in complex-valued Sobolev space and the corresponding projections in complex-valued finite element space $V_{h}^{c}$. 
\end{remark}

%


\begin{lemma}[Discrete Gronwall's inequality \cite{Heywood1990Finite}]\label{GronwallLem}
Let $\tau$, $B$, and $a_{k}$, $b_{k}$, $c_{k}$, $\gamma_{k}$, for $k\geq0$, be nonnegative numbers satisfying
\begin{equation}
a_{n} + \tau\sum_{k=0}^{n} b_{k}\leq \tau\sum_{k=0}^{n}\gamma_{k}a_{k} + \tau\sum_{k=0}^{n}c_{k}+B, \quad \text{for} \quad 
n\geq 0. 
 \end{equation}
Suppose that $\tau \gamma_{k}<1$, for all $k$, and $\sigma_{k}=(1-\tau\gamma_{k})^{-1}$. Then
\begin{equation}
a_{n} + \tau\sum_{k=0}^{n}b_{k} \leq \exp\left(\tau\sum_{k=0}^{n}\sigma_{k}\gamma_{k}\right)\left(\tau\sum_{k=0}^{n}c_{k} + B\right). 
\end{equation}
\end{lemma} 
\begin{lemma}\label{LemmaGL2} \cite{Zouraris2023Error}
Let $v^a, v^b, z^a, z^b \in \mathbb{C}$ and $\mathrm{S}\left(v^a, v^b, z^a, z^b\right):=\left|v^a\right|^2-\left|v^b\right|^2-\left|z^a\right|^2+\left|z^b\right|^2$. Then, 
\begin{equation}\label{Eq1}
\begin{aligned}
\left\|\mathrm{S}\left(v^a, v^b, z^a, z^b\right)\right\| \leq & 2\left\|z^a-z^b\right\|_{\infty}\left\|v^b-z^b\right\| +\mathrm{H}\left(v^a, v^b, z^a, z^b\right)\left\|v^a-v^b-z^a+z^b\right\|,
\end{aligned}
\end{equation}
where $\mathrm{H}\left(v^a, v^b, z^a, z^b\right):=\left\|v^a\right\|_{\infty}
+\left\|v^b\right\|_{\infty}
+\left\|z^a-z^b\right\|_{\infty}$.
\end{lemma}

For the finite element approximation related to the Poisson problem, the following estimate holds. 
\begin{lemma}\label{LemmaL2}
Given~$f\in L^{2}(\Omega)$. If $a\in H_{0}^{1}(\Omega)$ satisfies
\begin{equation}\label{BiEq}
A_{1}(a,\chi_h)=(f,\chi_h),\quad \forall \chi_h\in V_{h}.
\end{equation}
Then there exists a constant $C>0$ such that
\begin{equation}\label{LemAL2}
\left\|a\right\|\leq C\left(\left\|f\right\|+h\min_{a_{h}\in V_{h}}\|a_{h}-a\|_{1}\right). 
\end{equation}
\end{lemma} 
\begin{proof}
Let $a_{h}\in V_{h}$ be an approximation of $a$. Then \eqref{BiEq} can be reformulated as
\begin{equation}\label{DualBilin}
A_{1}(a_{h},\chi_h) = A_{1}(a_{h}-a, \chi_h) + (f, \chi_h).
\end{equation}
Taking $\chi_h=a_{h}$ in \eqref{DualBilin} and applying \eqref{coercon} give
\begin{equation}
\|a_{h}\|_{1}^{2}\leq \frac{\gamma_{2}}{\gamma_{1}}\|a_{h}-a\|_{1}\|a_{h}\|_{1}+\frac{1}{\gamma_{1}}\|f\|\|a_h\|.
\end{equation}
Note that $\|a_h\|\leq C\|a_h\|_{1}$. We obtain 
\begin{equation}\label{ahbdd}
\|a_{h}\|_{1}\leq C\left(\|a_{h}-a\|_{1}+\|f\|\right).
\end{equation}
By using the triangle inequality $\|a\|_1-\|a-a_h\|_1 \leq \|a_h\|_1$, \eqref{ahbdd} yields
\begin{equation}\label{DualAH1}
\|a\|_{1}\leq C\left(\min_{a_{h}\in V_{h}}\|a_{h}-a\|_{1}+\|f\|\right),
\end{equation}

On the other hand, we introduce a function $\psi$ solving the elliptic problem
\begin{equation}\label{Poisson}
-\Delta \psi = a \quad \text{in}\ \Omega,\qquad \psi = 0 \quad \text{on}\ \partial \Omega,
\end{equation}
which holds the regularity estimate $\psi\in H^{2}(\Omega)$ and
\begin{equation}\label{DualH2}
    \|\psi\|_{2}\leq \|a\|.
\end{equation}
From \eqref{Poisson}, it follows
\begin{equation}
\|a\|^{2}= \int_{\Omega} a \cdot (-\Delta  \psi)dx  = \int_{\Omega}(\nabla a \cdot \nabla \psi)dx = A_1(a, \psi).
\end{equation}
Let $\psi_{I}\in V_{h}$ be a piecewise linear interpolant of $\psi$. Then 
\begin{equation}\label{L2lift}
\begin{split}
\|a\|^{2}&= A_1(a, \psi) = A_1(a, \psi_I) + A_1(a, \psi-\psi_I) \\
&=(f,\psi_{I})+A_{1}(a,\psi-\psi_{I})\\
&\leq \|f\|(\|\psi\|+ \|\psi-\psi_{I}\|)
+ \gamma_{2}\|a\|_{1}\|\psi-\psi_{I}\|_{1}\\
&\leq \|f\|(\|\psi\|+Ch^{2}\|\psi\|_{2})+Ch\|a\|_{1} \|\psi\|_{2}\\
&\leq C(\|f\|+h\|a\|_{1})\|a\|,
\end{split}
\end{equation}
where we have used the regularity \eqref{DualH2} and the projection errors
\begin{equation*}
\|\psi-\psi_{I}\|\leq Ch^{2}\|\psi\|_{2}, \quad \| (\psi-\psi_{I})\|_1\leq Ch\|\psi\|_{2}.
\end{equation*}
\eqref{L2lift} together with \eqref{DualAH1} yields the estimate \eqref{LemAL2}.
\end{proof}

We define the discrete Laplacian operator $\Delta_{h}:\mathbf{H}_{0}^{1}(\Omega)\rightarrow V_{h}^{c}$ as
\begin{equation}\label{DiscreteLaplace}
\langle -\Delta_{h} v, \chi_{h}\rangle = \langle\nabla v, \nabla \chi_{h}\rangle, \quad \forall\chi_{h}\in V_{h}^{c}.
\end{equation}
We also introduce linear operators $\mathbf{S}_{h}, \mathbf{T}_{h}: V_{h}^{c}\rightarrow V_{h}^{c}$,
\begin{align}
\left\langle\mathbf{S}_{h}v_{h},\omega_{h}\right\rangle=\left\langle\left(\mathbf I_{h}-\mathbf{i}\frac{\tau}{2}\Delta_{h}\right)v_{h}, \omega_{h}\right\rangle, \quad \forall \omega_{h}\in V_{h}^{c}, \label{Shoper1}\\
\left\langle\mathbf{T}_{h}v_{h},\omega_{h}\right\rangle =\left\langle\left(\mathbf I_{h}+\mathbf{i}\frac{\tau}{2}\Delta_{h}\right)v_{h}, \omega_{h}\right\rangle, \quad \forall \omega_{h}\in V_{h}^{c},\label{Shoper2}
\end{align}    
where $\mathbf I_{h}$ is an identity operator on $V_h^c$. Denoting by $\mathbf{O}_{h}=\mathbf{S}_{h}, \mathbf{T}_{h}$ and setting $\omega_h = v_h$ in \eqref{Shoper1} and \eqref{Shoper2} give
\begin{equation}
\mathrm{Re} (\mathbf{O}_{h}v_{h},v_{h})=\|v_{h}\|^{2}, \quad \forall v_{h} \in V_{h}^{c},
\end{equation}
which implies ker$(\mathbf{O}_{h})=\{0\}$. Therefore, the operators $\mathbf{S}_{h}$ and $\mathbf{T}_{h}$ are invertible.

Similar to \cite[Lemma 2.4]{Zouraris2023Error}, the following statement holds. 
\begin{lemma}\label{Operbdd}
The operators $\mathbf{S}_{h}$ defined in \eqref{Shoper1} and $\mathbf{T}_{h}$ in \eqref{Shoper2} are invertible and fulfill
\begin{align}
\|\mathbf{S}_{h}^{-1}(v_{h})\|\leq \|v_{h}\|, \quad \forall v_{h}\in V_{h}^{c},\\
\|\mathbf{B}_{h}(v_{h})\|\leq \|v_{h}\|, \quad  \forall v_{h}\in V_{h}^{c}, 
\end{align} 
where the linear operator $\mathbf{B}_{h}:V_{h}^{c}\rightarrow V_{h}^{c}$ is given by  
\begin{equation}\label{OperBh}
\mathbf{B}_{h}:=\mathbf{S}_{h}^{-1}\mathbf{T}_{h}.
\end{equation}
\end{lemma}

\begin{lemma}\label{ZourariResult}
Let $\mathbf{I}_h$, $\mathbf{S}_h$, and $\mathbf{B}_h$ be the operators in \eqref{Shoper1}, \eqref{Shoper2}, and \eqref{OperBh}, and let $\{y^n\}_{n=1}^N$ be a sequence in $V_h^c$ satisfying: 
\begin{equation}\label{lemcond1}	
y^{n+1}=(\mathbf{B}_{h}-\mathbf{I}_{h})y^n+\mathbf{B}_{h}y^{n-1}
+\mathbf{S}_{h}^{-1}\Gamma^{n+1},
\end{equation} 
where $\{\Gamma^{n+1}\}_{n=1}^N$ are given functions in $V_h^c$.
Then, for $n\geq 2$ it follows
\begin{equation}\label{Dtaubdd-}
\|y^{n+1}\|+\|y^{n}\|\leq 2\left\|\mathbf{S}_{h}(y^{2})\right\| + 2\left\|\mathbf{S}_{h}(y^{1})\right\| + 2 \sum_{l=2}^{n}\|\Gamma^{l+1}\|.
\end{equation}
\end{lemma}
\begin{proof}
The proof is summarized from Part 9 in the proof of Theorem 3.1 in \cite{Zouraris2023Error}, we present it here for completeness. If $n = 1$ in \eqref{Dtaubdd-}, the estimate is obvious by using \Cref{Operbdd} and $y^i = \mathbf{S}_{h}^{-1} \mathbf{S}_{h} y^i$ for $i=1,2$. 
Next, we will focus on $n \geq 2$. 
Note that \eqref{lemcond1} can be written in a vector form 
\begin{equation}
\left[
\begin{array}{cc}
y^{n+1} \\
y^{n}
\end{array}
\right] = M\left[
\begin{array}{cc}
y^{n} \\
y^{n-1} 
\end{array}
\right] + 
\left[
\begin{array}{cc}
F^{n+1}\\
0
\end{array}
\right], 
\end{equation}
where 
\begin{equation}
M = 
\left[
\begin{array}{cc}
\mathbf{B}_{h}-I_{h} & \mathbf{B}_{h} \\
\mathbf{I}_{h} & 0 
\end{array}
\right] \quad \text{and} \quad F^{n+1}:= \mathbf{S}_{h}^{-1}\Gamma^{n+1}.
\end{equation}
A simple induction argument yields
\begin{equation}\label{TwolEVELdU}
\left[
\begin{array}{cc}
y^{n+1} \\
y^{n}
\end{array}
\right] = M^{n-1}\left[
\begin{array}{cc}
y^{2} \\
y^{1} 
\end{array}
\right]+\sum_{l=2}^{n}M^{n-l}
\left[
\begin{array}{cc}
F^{l+1}\\
0
\end{array}
\right],
\end{equation}
where 
\begin{equation}\label{Gkappa}
M^{\kappa}
= \frac{1}{2}\left[\begin{array}{cc}
((-1)^{\kappa}\mathbf{I}_{h}+\mathbf{B}_{h}^{\kappa+1})\mathbf{S}_{h} & ((-1)^{\kappa+1}\mathbf{B}_{h}+\mathbf{B}_{h}^{\kappa+1})\mathbf{S}_{h}\\
((-1)^{\kappa+1}\mathbf{I}_{h}+\mathbf{B}_{h}^{\kappa})\mathbf{S}_{h} & ((-1)^{\kappa}\mathbf{B}_{h}+\mathbf{B}_{h}^{\kappa})\mathbf{S}_{h}
\end{array}
\right]. 
\end{equation}
Plugging \eqref{Gkappa} into \eqref{TwolEVELdU} yields
\begin{equation}
\begin{split}
\left[
\begin{array}{cc}
y^{n+1} \\
y^{n}
\end{array}
\right] &= 
\frac{1}{2}\left[
\begin{array}{cc}
((-1)^{n-1}\mathbf{I}_{h}+\mathbf{B}_{h}^{n})\mathbf{S}_{h} &
((-1)^{n}\mathbf{B}_{h} +\mathbf{B}_{h}^{n})\mathbf{S}_{h} \\
((-1)^{n} \mathbf{I}_{h}+\mathbf{B}_{h}^{n-1})\mathbf{S}_{h} & ((-1)^{n-1}\mathbf{B}_{h}+\mathbf{B}_{h}^{n-1})\mathbf{S}_{h}
\end{array}
\right]
\left[
\begin{array}{cc}
y^{2} \\
y^{1} 
\end{array}
\right]\\
&\quad +\frac{1}{2} \sum_{l=2}^{n}\left[\begin{array}{cc}
((-1)^{n-l}\mathbf{I}_{h}+\mathbf{B}_{h}^{n-l+1})\mathbf{S}_{h} & ((-1)^{n-l+1}\mathbf{B}_{h}+\mathbf{B}_{h}^{n-l+1})\mathbf{S}_{h}\\
((-1)^{n-l+1}\mathbf{I}_{h}+\mathbf{B}_{h}^{n-l})\mathbf{S}_{h} & ((-1)^{n-l}\mathbf{B}_{h}+\mathbf{B}_{h}^{n-l})\mathbf{S}_{h}
\end{array}
\right]
\left[
\begin{array}{cc}
F^{l+1} \\ 0
\end{array}
\right], 
\end{split}
\end{equation}
which gives for $i=n, n+1$,
\begin{equation}\label{lemcond2}
\begin{split}
y^{i}&=\frac{1}{2}\left[(-1)^{i} \mathbf{I}_{h}+\mathbf{B}_{h}^{i-1}\right]\mathbf{S}_{h}y^{2}+\frac{1}{2}\left[(-1)^{i-1}\mathbf{B}_{h}+\mathbf{B}_{h}^{i-1}\right]\mathbf{S}_{h}y^{1}+\frac{1}{2}\sum_{l=2}^{n}\left[(-1)^{i-l+1}\mathbf{I}_{h}+\mathbf{B}_{h}^{i-l}\right]\Gamma^{l+1},
\end{split}
\end{equation} 
which together with \Cref{Operbdd} yields \eqref{Dtaubdd-}.
\end{proof}

For error analysis purposes, we assume that the exact solutions $u$, $\Phi$ and $\Psi$ in \eqref{TargetEqRe} hold the following regularity 
\begin{equation}\label{RegularCondition}
\begin{split}
&u, u_{t} \in L^{\infty}(0,T; \mathbf{H}^{k+1}(\Omega)),\qquad \Psi,\Psi_{t}, \Phi, \Phi_{t} \in L^{\infty}(0,T; H^{k+1}(\Omega)),\\
&u_{tt}, \Psi_{tt}\in L^{\infty}(0,T; \mathbf{H}^{2}(\Omega)), \qquad \Psi_{ttt} \in L^{\infty}(0,T; \mathbf{L}^{2}(\Omega)), \qquad u_{ttt}, u_{tttt}\in L^{\infty}(0,T; \mathbf{L}^{2}(\Omega)).
\end{split}
\end{equation}
In addition, we also assume that the external potential $V(x) \in L^{\infty}(\Omega)$.

In view of the regularity assumptions in \eqref{RegularCondition} for the exact solution $u$, $\Psi$ and $\Phi$, and \Cref{InftyLem}, we have for any $n \geq 0$,
\begin{equation}\label{exactsInfty}
\begin{aligned}
& \|u^{n}\|_{\infty}\leq C_{u}, \qquad \|\Psi^{n-1/2}\|_{\infty}\leq C_{\Psi}, \qquad \|\Phi^{n-1/2}\|_{\infty}\leq C_{\Phi},\\
& \|\mathrm{R}_{h}u^{n}\|_{\infty}\leq D_{u}, \qquad \|\mathrm{R}_{h}\Psi^{n-1/2}\|_{\infty}\leq D_{\Psi}, \qquad \|\mathrm{R}_{h}\Phi^{n-1/2}\|_{\infty}\leq D_{\Phi},\\
\end{aligned}
\end{equation}
where the constants
\begin{align*}
& C_{u}=\sup\limits_{0\leq n \leq N}\|u^{n}\|_{\infty}, \qquad  C_{\Psi}=\sup\limits_{0\leq n \leq N}\|\Psi^{n-1/2}\|_{\infty}, \qquad  C_{\Phi}=\sup\limits_{0\leq n \leq N}\|\Phi^{n-1/2}\|_{\infty}, \\
& D_{u}=\sup\limits_{0\leq n \leq N}\|\mathrm{R}_{h}u^{n}\|_{\infty}, \qquad D_{\Psi}=\sup\limits_{0\leq n \leq N}\|\mathrm{R}_{h}\Psi^{n-1/2}\|_{\infty}, \qquad  D_{\Phi}=\sup\limits_{0\leq n \leq N}\|\mathrm{R}_{h}\Phi^{n-1/2}\|_{\infty}.
\end{align*}

Recall that the exact solution of \eqref{TargetEqRe} satisfies
\begin{subequations}\label{EaxctEq}
\begin{align}
&\left(\Psi^{n+1/2}+\Psi^{n-1/2}, v \right)=\left(S_{1}^{n},v\right)+\left(2|u^{n}|^{2}, v\right),  \label{EaxctEq1}   \\
&\mathbf{i}\left\langle D_{\tau}u^{n+1},\omega\right\rangle = A_{0}\left(\overline{u}^{n+1/2}, \omega\right)
+\left\langle\left(\Phi^{n+1/2}+V(x)
+\Psi^{n+1/2}\right)\overline{u}^{n+1/2},\omega\right\rangle+\langle R_{1}^{n+1}, \omega\rangle, \label{EaxctEq2}\\
&A_{1}\left(\Phi^{n+1/2}, \chi\right)=\mu\left(\Psi^{n+1/2}-c, \chi\right),\label{EaxctEq3}
\end{align}
\end{subequations}
for any $v,\chi \in H_{0}^{1}(\Omega)$ and $\omega \in \mathbf{H}_{0}^{1}(\Omega)$, where the consistency errors 
\begin{equation*}
S_{1}^{n}=\Psi^{n+1/2}+\Psi^{n-1/2}-2\Psi^{n},
\end{equation*}
and
\begin{equation}\label{R1form}
\begin{aligned}
R_{1}^{n+1}= &-\mathbf{i}(u_{t}^{n+1/2}-D_{\tau}u^{n+1}) +\Delta(\overline{u}^{n+1/2}-u^{n+1/2}) \\
&+(\Phi^{n+1/2}+V(x) +\Psi^{n+1/2})(u^{n+1/2}-\overline{u}^{n+1/2}).
\end{aligned}
\end{equation}
We define the errors $e_{u}^{n+1}$, $e_{\Psi}^{n+1/2}$ and $e_{\Phi}^{n+1/2}$ with $0\leq n\leq N-1$ as
\begin{align*}
e_{u}^{n+1}=u^{n+1}-u_{h}^{n+1}, \quad  e_{\Psi}^{n+1/2}=\Psi^{n+1/2}-\Psi_{h}^{n+1/2}, \quad e_{\Phi}^{n+1/2}=\Phi^{n+1/2}-\Phi_{h}^{n+1/2}.
\end{align*}
Taking $v=v_{h}$, $\omega=\omega_{h}$ and $\chi=\chi_{h}$ in \eqref{EaxctEq}, and subtracting \eqref{fullyVdiscrete} from \eqref{EaxctEq} yield
\begin{subequations}\label{Error}
\begin{align}
&(e_{\Psi}^{n+1/2}+e_{\Psi}^{n-1/2}, v_{h})=(S_{1}^{n},v_{h})
+(T_{1}^{n}, v_{h}),\label{Error1}\\
&\mathbf{i}\langle D_{\tau}e_{u}^{n+1},\omega_{h}\rangle
=A_{0}(\overline{e}_{u}^{n+1/2},\omega_{h})
+\langle G_{1}^{n+1},\omega_{h}\rangle
+\langle R_{1}^{n+1}, \omega_{h}\rangle, \label{Error2}\\
&A_{1}(e_{\Phi}^{n+1/2}, \chi_{h})=\mu(e_{\Psi}^{n+1/2}, \chi_{h}),\label{Error3}
\end{align}
\end{subequations}
where
\begin{equation*}
T_{1}^{n}=2|u^{n}|^{2}-2|u_{h}^{n}|^{2},
\end{equation*}
and
\begin{equation}\label{G1n1}
\begin{aligned}
G_{1}^{n+1}=&(\Phi^{n+1/2}\overline{u}^{n+1/2}
-\Phi_{h}^{n+1/2}\overline{u}_{h}^{n+1/2})
+V(x)(\overline{u}^{n+1/2}-\overline{u}_{h}^{n+1/2})\\
&+(\Psi^{n+1/2}\overline{u}^{n+1/2}-\Psi_{h}^{n+1/2}\overline{u}_{h}^{n+1/2}).
\end{aligned}
\end{equation}

By using the projection operator $\mathrm{R}_{h}$, the errors $e_{u}^{n+1}$, $e_{\Psi}^{n+1/2}$ and $e_{\Phi}^{n+1/2}$ can be split as
\begin{align}
&e_{u}^{n+1}
=(u^{n+1}-\mathrm{R}_{h}u^{n+1})
+(\mathrm{R}_{h}u^{n+1}-u_{h}^{n+1})
=\xi_{u}^{n+1}+\eta_{u}^{n+1},\label{eudecomp}\\
&e_{\Psi}^{n+1/2}=(\Psi^{n+1/2}-\mathrm{R}_{h}\Psi^{n+1/2})
+(\mathrm{R}_{h}\Psi^{n+1/2}-\Psi_{h}^{n+1/2})
=\xi_{\Psi}^{n+1/2}+\eta_{\Psi}^{n+1/2},\label{eudecompPsi}\\
&e_{\Phi}^{n+1/2}=(\Phi^{n+1/2}-\mathrm{R}_{h}\Phi^{n+1/2})
+(\mathrm{R}_{h}\Phi^{n+1/2}-\Phi_{h}^{n+1/2})
=\xi_{\Phi}^{n+1/2}+\eta_{\Phi}^{n+1/2}.\label{eudecompPhi}
\end{align}
Thus, the equivalent form of the error equations \eqref{Error} are presented as
\begin{subequations}
\begin{align}
&(\eta_{\Psi}^{n+1/2}+\eta_{\Psi}^{n-1/2}, v_{h})=(S_{2}^{n},v_{h})
+(T_{1}^{n}, v_{h}), \label{EquError1}   \\
&\mathbf{i}\langle D_{\tau}\eta_{u}^{n+1},\omega_{h}\rangle
=A_{0}(\overline{\eta}_{u}^{n+1/2},\omega_{h})
+\langle G_{1}^{n+1},\omega_{h}\rangle
+\langle R_{2}^{n+1},\omega_{h}\rangle, \label{EquError2}\\
& A_{1}(\eta_{\Phi}^{n+1/2},\chi_{h})=\mu(\eta_{\Psi}^{n+1/2}, \chi_{h})+ \mu(R_{3}^{n+1/2},\chi_{h}),\label{EquError3}
\end{align}
\end{subequations}
where
\begin{align*}
S_{2}^{n}:=S_{1}^{n}-(\xi_{\Psi}^{n+1/2}+\xi_{\Psi}^{n-1/2}),\quad R_{2}^{n+1}:=R_{1}^{n+1}-\mathbf i D_{\tau}\xi_{u}^{n+1},\quad R_{3}^{n+1/2}:=\xi_{\Psi}^{n+1/2},
\end{align*}
and we have used \eqref{Ritz} to get rid of the terms 
$A_{0}(\overline{\xi}_{u}^{n+1/2},\omega_{h})$ and $A_{1}(\xi_{\Phi}^{n+1/2},\chi_{h})$.
By using the projection error \eqref{Ritz2} and the mean value theorem, it holds
\begin{equation}\label{SREDtauqb1}
\begin{aligned}
\|D_{\tau}\xi_{u}^{n+1}\| = & \left\| D_\tau u^{n+1}-R_h D_\tau u^{n+1} \right\|\leq Ch^{k+1}\left\| D_\tau u^{n+1} \right\|_{k+1} \leq  Ch^{k+1}\left\| u_t(x,t^*) \right\|_{k+1},
\end{aligned}
\end{equation}
where $t^* \in (t_{n},t_{n+1})$.
Then applying the Taylor expansion and the properties of the interpolation operator, for any $n \geq 0$, gives the estimates
\begin{align}
&\|S_{2}^{n}\|\leq C(\tau^{2}+h^{k+1}),\label{SREqb1}\\
&\|R_{2}^{n+1}\|\leq C(\tau^{2}+h^{k+1}),\label{SREqb2}\\
&\|R_{3}^{n+1/2}\|\leq Ch^{k+1}.\label{SREqb3}
\end{align}

Then we obtain the following error estimates.
\begin{theorem}\label{Thm+}
Suppose that $u$, $\Psi$ and $\Phi$ satisfy the regularity conditions \eqref{RegularCondition}. If $\tau\leq Ch$, then there exists constant $\tau_{0}>0$ and $h_0>0$ such that when time step $\tau<\tau_{0}$ and mesh size $h<h_0$, the solution of the relaxation Crank-Nicolson finite element scheme \eqref{fullyVdiscrete} satisfies
\begin{align}
&\max_{0\leq n \leq N}\left\|e_u^{n}\right\|\leq C\left(\tau^{2}+h^{k+1}\right),\label{EstimateErrora}\\
&\max_{0\leq n \leq N-1}\left\|e_\Psi^{n+1/2}\right\|\leq C\left(\tau^{2}+h^{k+1}\right),\label{EstimateErrorb}\\
&\max_{0\leq n\leq N-1}\left\|e_\Phi^{n+1/2}\right\|\leq C(\tau^{2}+h^{k+1}).\label{EstimateErrorc}
\end{align}
\end{theorem}

\begin{proof}
We prove the results using the method of mathematical induction.\\

\textbf{Step 1}. 
In this step, we prove the following estimates.
\begin{align}
&\|e_{\Psi}^{1/2}\|\leq  C(\tau^{2}+h^{k+1}),\\
&\|e_{\Phi}^{1/2}\|\leq  C(\tau^{2}+h^{k+1}),\\
&\|e_u^{1}\| \leq  C\left(\tau^{2}+h^{k+1}\right),\label{uinti1}\\
&\|D_\tau \eta_u^{1}\| \leq  C\left(\tau^{2}+h^{k+1}\right)\label{Dtauetau1}.
\end{align}

For $n=0$, taking 
$v_{h}=\eta_{\Psi}^{1/2}-\eta_{\Psi}^{-1/2}$ in \eqref{EquError1} gives
\begin{equation}\label{ThmEq1}
\begin{split}
\|\eta_{\Psi}^{1/2}\|^{2}-\|\eta_{\Psi}^{-1/2}\|^{2}
&=\left(S_{2}^{0},\eta_{\Psi}^{1/2}-\eta_{\Psi}^{-1/2}\right)
+\left(T_{1}^{0}, \eta_{\Psi}^{1/2}-\eta_{\Psi}^{-1/2}\right)\\
&\leq 2\|S_{2}^{0}\|^{2} + \frac{1}{2} \|\eta_{\Psi}^{1/2}\|^{2} + \frac{1}{2}\|\eta_{\Psi}^{-1/2}\|^{2}+2\|T_{1}^{0}\|^{2}.
\end{split}
\end{equation}
Note that the following inequalities hold
\begin{equation}\label{initerr}
\begin{aligned}
&\|T_{1}^{0}\|\leq2\||u_{0}|^{2}-|u_{h}^{0}|^{2}\|\leq 2\|u_{0}+u_{h}^{0}\|_\infty \|u_{0}-u_{h}^{0}\|\leq Ch^{k+1},\\
&\|\eta_{\Psi}^{-1/2}\|\leq \|e_{\Psi}^{-1/2}\|+\|\xi_{\Psi}^{-1/2}\|
\leq Ch^{k+1},
\end{aligned}
\end{equation}
which together with \eqref{SREqb1} when plugging into \eqref{ThmEq1} yields
\begin{equation}\label{etapsi0}
\begin{split}
\|\eta_{\Psi}^{1/2}\|^{2}
\leq 3\|\eta_{\Psi}^{-1/2}\|^{2}
+4\|S_{2}^{0}\|^{2}
+4\|T_{1}^{0}\|^{2}\leq C(\tau^{2}+h^{k+1})^{2}.
\end{split}
\end{equation}
By \eqref{etapsi0} and the projection error \eqref{Ritz2} for $\|\xi_{\Psi}^{1/2}\|$,
\begin{equation}\label{errorpsi0}
\|e_{\Psi}^{1/2}\|\leq\|\eta_{\Psi}^{1/2}\|+\|\xi_{\Psi}^{1/2}\|\leq C(\tau^{2}+h^{k+1}).
\end{equation}
By applying the Lemma \ref{LemmaL2} to the \eqref{Error3} with $n=0$, we conclude the following error estimate
\begin{equation}\label{ThmEq3}
\begin{split}
\|e_{\Phi}^{1/2}\|\leq C\|e_{\Psi}^{1/2}\|+Ch^{k+1}\leq C(\tau^{2}+h^{k+1}).
\end{split}
\end{equation}
In view of $\tau \leq Ch$, \eqref{exactsInfty}, \eqref{etapsi0}, \eqref{ThmEq3}, and the inverse inequality \eqref{InverIeq}, there exist $h_1>0$ such that when $h< h_1$,
\begin{align}
\|\Psi_{h}^{1/2}\|_{\infty}
\leq & \|\mathrm{R}_{h}\Psi^{1/2}\|_{\infty}+\|\eta_{\Psi}^{1/2}\|_{\infty}
\leq \|\mathrm{R}_{h}\Psi^{1/2}\|_{\infty}+Ch^{-1}\|\eta_{\Psi}^{1/2}\|\leq D_\Psi+C_\Psi h \leq D_\Psi+1, \label{Psiinfprof}\\
\|\Phi_{h}^{1/2}\|_{\infty}\leq & \|\mathrm{R}_{h}\Phi^{1/2}\|_{\infty} +Ch^{-1} \|\eta_{\Phi}^{1/2}\|\leq D_\Phi+C_\Psi h\leq D_\Phi+1.\label{ThmEq4}
\end{align}

Taking $\omega_{h}=\overline{\eta}_{u}^{1/2}$ in \eqref{EquError2} with $n=0$ gives
\begin{equation*}
\mathbf{i}\left\langle D_{\tau}\eta_{u}^{1},\overline{\eta}_{u}^{1/2}\right\rangle
=A_{0}\left(\overline{\eta}_{u}^{1/2},\overline{\eta}_{u}^{1/2}\right)
+\left\langle G_{1}^{1},\overline{\eta}_{u}^{1/2}\right\rangle
+\left\langle R_{2}^{1}, \overline{\eta}_{u}^{1/2}\right\rangle,
\end{equation*}
where the imaginary part yields
\begin{equation}\label{ThmEq2a}
 \begin{split}
\frac{1}{2\tau}\left(\|\eta_{u}^{1}\|^{2}
-\|\eta_{u}^{0}\|^{2}\right)
&=\mathrm{Im}\left\langle G_{1}^{1}+R_{2}^{1},\overline{\eta}_{u}^{1/2}\right\rangle \\
&\leq \left\|G_{1}^{1}+R_{2}^{1}\right\|\left\|\overline{\eta}_{u}^{1/2}\right\|\\
&\leq\frac{1}{2}\left\|G_{1}^{1}+R_{2}^{1}\right\|^2 +\frac{1}{2}\left\|\overline{\eta}_{u}^{1/2}\right\|^2\\
&=\frac{1}{2}\left(\left\|G_{1}^{1}\right\|^2 +\left\|R_{2}^{1}\right\|^2 + 2\left\|G_{1}^{1}\right\|\left\|R_{2}^{1}\right\| \right)+\frac{1}{8}\left\|\eta_{u}^{1}+\eta_{u}^{0}\right\|^2\\
&\leq \left\|G_{1}^{1}\right\|^{2}
+\left\|R_{2}^{1}\right\|^{2} 
+\frac{1}{4}\left(\|\eta_{u}^{1}\|^{2}+\|\eta_{u}^{0}\|^{2}\right).
\end{split}
\end{equation}
By employing \eqref{exactsInfty}, \eqref{errorpsi0}-\eqref{ThmEq4}, 
\begin{equation}\label{ThmEq2b}
\begin{split}
\left\|G_{1}^{1}\right\|
&\leq \|\overline{u}^{1/2}\|_{\infty}\left(\|e_{\Phi}^{1/2}\|+\|e_{\Psi}^{1/2}\|\right) 
+\left(\|\Phi_{h}^{1/2}\|_{\infty}+\|V(x)\|_{\infty}
+\|\Psi_{h}^{1/2}\|_{\infty}\right)\|\overline{e}_{u}^{1/2}\|\\
&\leq C\left(\|\eta_{u}^{1} \|+\|\eta_{u}^{0}\|\right)
+C\left(\tau^{2}+h^{k+1}\right).
\end{split}
\end{equation}
Plugging \eqref{SREqb2} and \eqref{ThmEq2b} into \eqref{ThmEq2a} gives
\begin{equation}\label{ThmEq5}
\begin{split}
\frac{1}{2\tau}\left(\|\eta_{u}^{1}\|^{2}
-\|\eta_{u}^{0}\|^{2}\right)
&\leq C_{1}\left(\|\eta_{u}^{1}\|^{2}+\|\eta_{u}^{0}\|^{2}\right)
+C(\tau^{2}+h^{k+1})^{2}.
\end{split}
\end{equation}
Since the initial value $\eta_{u}^{0}=0$, \eqref{ThmEq5} leads to
\begin{equation}\label{etau1level}
\left\|\eta_{u}^{1}\right\|\leq C\tau(\tau^{2}+h^{k+1}),
\end{equation} 
as long as $\tau<\tau_{1}:=1/(2C_{1})$. Since $0<\tau<1$,  we then conclude that
\begin{equation*}
\left\|e_{u}^{1}\right\|\leq \|\xi_{u}^{1}\|+\|\eta_{u}^{1}\|\leq C(\tau^{2}+h^{k+1}).
\end{equation*}
Again, using $\eta_{u}^{0}=0$ and \eqref{etau1level} gives
\begin{align}\label{u1infbddA}
\|D_\tau \eta_u^1\| = \frac{1}{\tau} \left\|\eta_{u}^{1}\right\| \leq C(\tau^{2}+h^{k+1}).
\end{align}
Based on \eqref{InverIeq}, \eqref{exactsInfty} and \eqref{etau1level}, there exists $h_2$ such that when $h<h_2$,
\begin{equation}\label{u1infbdd}
\|u_{h}^{1}\|_{\infty}\leq \|\mathrm{R}_{h}u^{1}\|_{\infty}
+\|\mathrm{R}_{h}u^{1}-u_{h}^{1}\|_{\infty}\leq \|\mathrm{R}_{h}u^{1}\|_{\infty}
+Ch^{-1}\|\eta_{u}^{1}\|\leq D_u+C_u h \leq D_u + 1.
\end{equation}

\textbf{Step 2}. 
In this step, we prove the following estimates
\begin{align}
&\left\|e_u^{2}\right\|\leq C\left(\tau^{2}+h^{k+1}\right), \\
&\left\|D_\tau \eta_u^{2}\right\|\leq C\left(\tau^{2}+h^{k+1}\right), \\
&\max_{1\leq n \leq 2}\left\|e_\Psi^{n+1/2}\right\|\leq C\left(\tau^{2}+h^{k+1}\right), \\
&\max_{1\leq n\leq 2}\left\|e_\Phi^{n+1/2}\right\|\leq C(\tau^{2}+h^{k+1}). 
\end{align}
Taking the difference between $t_{1}$ and $t_{0}$ of \eqref{EquError1} with $n=1$ leads to 
\begin{equation}\label{mainproof5++}
\begin{split}		
\left(\eta_{\Psi}^{3/2}-\eta_{\Psi}^{-1/2},v_{h}\right)
=\left(S_{2}^{1}-S_{2}^{0},v_{h}\right)
+\left(T_{1}^{1}-T_{1}^{0}, v_{h}\right). 
\end{split}
\end{equation}	
By \Cref{LemmaGL2} and \eqref{u1infbdd}, it follows that 
\begin{equation}\label{mainproof6++}
\begin{split}	
\|T_{1}^{1}-T_{1}^{0}\|
&=2\left\||u^{1}|^{2}-|u_h^{1}|^{2}-|u^{0}|^{2}
+|u_h^{0}|^{2}\right\|\\
&\leq 2 \|u^{1}-u^{0}\|_{\infty}\|u^{1}-u_{h}^{1}\| \\
&\quad +\left(\|u_{h}^{1}\|_{\infty} +\|u_{h}^{0}\|_{\infty} +\|u^{1}-u^{0}\|_{\infty}\right)\left\|u_h^{0}-u_h^{1}-u^{0}
+u^{1}\right\|\\
&\leq C\|u^{1}-u^{0}\|_{\infty}\|e_{u}^{1}\| + C \|e_{u}^{1}-e_{u}^{0}\| + C\|u^{1}-u^{0}\|_{\infty}\|e_{u}^{1}-e_{u}^{0}\|\\
&\leq C \|e_{u}^{1}-e_{u}^{0}\| + C\|u^{1}-u^{0}\|_{\infty}\left(\|e_{u}^{1}-e_{u}^{0}\|+ \|e_{u}^{1}\|\right)\\
&\leq C\|e_{u}^{1}-e_{u}^{0}\|+C\tau\|e_{u}^{1}-e_{u}^{0}\|+C\tau\|e_{u}^{1}\|\leq C\|e_{u}^{1}-e_{u}^{0}\|+C\tau\|e_{u}^{1}\|\\
&\leq C\tau\|D_{\tau}\eta_{u}^{1}\|+C\tau\|\eta_{u}^{1}\| + C\tau h^{k+1}.
\end{split}
\end{equation}
Note that
\begin{equation}\label{S2ndiff++}
\begin{split}
\|S_{2}^{1}-S_{2}^{0}\| 
&=\left\|\left(S_{1}^{1}-(\xi_{\Psi}^{3/2}+\xi_{\Psi}^{1/2})\right)
-\left(S_{1}^{0}-(\xi_{\Psi}^{1/2}+\xi_{\Psi}^{-1/2})\right)\right\|\\
&\leq\left\|S_{1}^{1}-S_{1}^{0}\right\|
+\left\|\xi_{\Psi}^{3/2}-\xi_{\Psi}^{-1/2}\right\|.
\end{split}
\end{equation}
By using the Taylor expression at $t_1$ and the regularity assumption \eqref{RegularCondition},
\begin{equation}\label{S3psi+}
\begin{split}
\left\|S_{1}^{1}-S_{1}^{0}\right\|=&\left\|\Psi^{3/2}-2\Psi^{1}+2\Psi^{0}-\Psi^{-1/2}\right\|\leq \left\|\frac{1}{2}\int_{t_1}^{t_{3/2}}(t_{3/2}-t)^{2}\Psi_{ttt}(x,t)dt \right.\\ 
& \left.+ \int_{t_1}^{t_{0}}(t_{0}-t)^{2}\Psi_{ttt}(x,t)dt-\frac{1}{2}\int_{t_1}^{t_{-1/2}}(t_{-1/2}-t)^{2}\Psi_{ttt}(x,t)dt\right\|\leq C\tau^{3}.
\end{split}
\end{equation}
By using the mean value theorem,
\begin{equation}\label{S3psiXi+}
\begin{split}
\|\xi_{\Psi}^{3/2}-\xi_{\Psi}^{-1/2}\|
&=2\tau
\left\|\left(\frac{\Psi^{3/2}-\Psi^{-1/2}}{2\tau}\right)-\mathrm{R}_{h}\left(\frac{\Psi^{3/2}-\Psi^{-1/2}}{2\tau}\right)\right\|\\
&\leq C\tau h^{k+1}\left\|\frac{\Psi^{3/2}-\Psi^{-1/2}}{2\tau}\right\|_{k+1}\leq C\tau h^{k+1}\left\|\Psi_{t}(x,t^{*}) \right\|_{k+1},
\end{split}
\end{equation}
where $t^{*}\in (t_{-1/2}, t_{3/2})$. 
Plugging \eqref{S3psi+} and \eqref{S3psiXi+} into \eqref{S2ndiff++} leads to 
\begin{equation}\label{DtS3n+}
\left\|S_{2}^{1}-S_{2}^{0}\right\|\leq C\tau(\tau^{2}+h^{k+1}).
\end{equation}
Then, taking $v_{h}=\eta_{\Psi}^{3/2}+\eta_{\Psi}^{-1/2}$ in \eqref{mainproof5++} yields
\begin{equation*} 
\begin{split}	
\|\eta_{\Psi}^{3/2}\|^{2}
-\|\eta_{\Psi}^{-1/2}\|^{2}
&\leq \left(\|S_{2}^{1}-S_{2}^{0}\| +\|T_{1}^{1}-T_{1}^{0}\|\right) \|\eta_{\Psi}^{3/2}+\eta_{\Psi}^{-1/2}\|.
\end{split} 
\end{equation*}
By \eqref{initerr}, \eqref{etau1level},  \eqref{u1infbddA}, \eqref{mainproof6++} and \eqref{DtS3n+}, the following inequality holds 
\begin{equation}\label{Etapsi2Leve+}
\begin{split}	
\|\eta_{\Psi}^{3/2}\|
&\leq \|\eta_{\Psi}^{-1/2}\|+\|S_{2}^{1}-S_{2}^{0}\|
+\|T_{1}^{1}-T_{1}^{0}\|\\
&\leq C\tau\|D_{\tau}\eta_{u}^{1}\|
+C\tau\|\eta_{u}^{1}\| + C\tau(\tau^{2}+h^{k+1}) \leq C(\tau^{2}+h^{k+1}),
\end{split}
\end{equation}
which together with the projection error \eqref{Ritz2} for $\|\xi_{\Psi}^{3/2}\|$ yields
\begin{equation}\label{epsiLevel2+}
\|e_{\Psi}^{3/2}\|\leq\|\eta_{\Psi}^{3/2}\|+\|\xi_{\Psi}^{3/2}\|\leq C(\tau^{2}+h^{k+1}).
\end{equation}
By applying the Lemma \ref{LemmaL2} to the \eqref{Error3} with $n=1$, we also obtain the following error estimate
\begin{equation}\label{ThmEq3++}
\begin{split}
\|e_{\Phi}^{3/2}\|\leq C\|e_{\Psi}^{3/2}\|+Ch^{k+1}\leq C(\tau^{2}+h^{k+1}).
\end{split}
\end{equation}
By using the inverse inequality \eqref{InverIeq}, 
\eqref{exactsInfty}, and \eqref{Etapsi2Leve+},  there exist $h_3>0$ such that when $h<h_3$,
\begin{align}
\|\Psi_{h}^{3/2}\|_{\infty}
\leq & \|\mathrm{R}_{h}\Psi^{3/2}\|_{\infty}+\|\eta_{\Psi}^{3/2}\|_{\infty}
\leq \|\mathrm{R}_{h}\Psi^{3/2}\|_{\infty}+Ch^{-1}\|\eta_{\Psi}^{3/2}\|\leq D_\Psi+C_\Psi h \leq D_\Psi+1, \label{Psiinfprof2+}\\
\|\Phi_{h}^{3/2}\|_{\infty}\leq & \|\mathrm{R}_{h}\Phi^{3/2}\|_{\infty} +Ch^{-1} \|\eta_{\Phi}^{3/2}\|\leq D_\Phi+C_\Psi h\leq D_\Phi+1.\label{ThmEq42+}
\end{align}
Taking $\omega_{h}=\overline{\eta}_{u}^{3/2}$ in \eqref{EquError2} with $n=1$ gives
\begin{equation*}
\mathbf{i}\left\langle D_{\tau}\eta_{u}^{2},\overline{\eta}_{u}^{3/2}\right\rangle
=A_{0}\left(\overline{\eta}_{u}^{3/2},\overline{\eta}_{u}^{3/2}\right)
+\left\langle G_{1}^{2},\overline{\eta}_{u}^{3/2}\right\rangle
+\left\langle R_{2}^{2}, \overline{\eta}_{u}^{3/2}\right\rangle,
\end{equation*}
where the imaginary part yields
\begin{equation}\label{ThmEq2a1++}
 \begin{split}
\frac{1}{2\tau}\left(\|\eta_{u}^{2}\|^{2}
-\|\eta_{u}^{1}\|^{2}\right)
&=\mathrm{Im}\left\langle G_{1}^{2}+R_{2}^{2},\overline{\eta}_{u}^{3/2}\right\rangle \\
&\leq \left\|G_{1}^{2}+R_{2}^{2}\right\|\left\|\overline{\eta}_{u}^{3/2}\right\|\\
&\leq\frac{1}{2}\left\|G_{1}^{2}+R_{2}^{2}\right\|^2 +\frac{1}{2}\left\|\overline{\eta}_{u}^{3/2}\right\|^2\\
&=\frac{1}{2}\left(\left\|G_{1}^{2}\right\|^2 +\left\|R_{2}^{2}\right\|^2 + 2\left\|G_{1}^{2}\right\|\left\|R_{2}^{2}\right\|\right)+\frac{1}{8}\left\|\eta_{u}^{2}+\eta_{u}^{1}\right\|^2\\
&\leq \left\|G_{1}^{2}\right\|^{2}
+\left\|R_{2}^{2}\right\|^{2} 
+\frac{1}{4}\left(\|\eta_{u}^{2}\|^{2}+\|\eta_{u}^{1}\|^{2}\right).
\end{split}
\end{equation}
By applying \eqref{exactsInfty}, \eqref{epsiLevel2+}-\eqref{ThmEq42+}, we have 
\begin{equation}\label{G2bound+}
\begin{split}
\left\|G_{1}^{2}\right\|
&\leq \|\overline{u}^{3/2}\|_{\infty}\left(\|e_{\Phi}^{3/2}\|+\|e_{\Psi}^{3/2}\|\right) 
+\left(\|\Phi_{h}^{3/2}\|_{\infty}+\|V(x)\|_{\infty}
+\|\Psi_{h}^{3/2}\|_{\infty}\right)\|\overline{e}_{u}^{3/2}\|\\
&\leq C\left(\|\eta_{u}^{2} \|+\|\eta_{u}^{1}\|\right)
+C\left(\tau^{2}+h^{k+1}\right).
\end{split}
\end{equation}
Plugging \eqref{SREqb2} and \eqref{G2bound+} into \eqref{ThmEq2a1++} gives
\begin{equation} \label{Eta2level}
\begin{split}
\frac{1}{2\tau}\left(\|\eta_{u}^{2}\|^{2}
-\|\eta_{u}^{1}\|^{2}\right)
&\leq C_{2}\left(\|\eta_{u}^{2}\|^{2}+\|\eta_{u}^{1}\|^{2}\right)
+C(\tau^{2}+h^{k+1})^{2}.
\end{split}
\end{equation}
In view of \eqref{etau1level}, \eqref{Eta2level} leads to
\begin{equation} \label{eta2levle++}
\left\|\eta_{u}^{2}\right\|\leq C\tau(\tau^{2}+h^{k+1}),
\end{equation} 
as long as $\tau<\tau_{2}:=1/(2C_{2})$. Since $0<\tau<1$,  we then conclude that
\begin{equation*}
\left\|e_{u}^{2}\right\|\leq \|\xi_{u}^{2}\|+\|\eta_{u}^{2}\|\leq C(\tau^{2}+h^{k+1}).
\end{equation*}
By using the triangle inequality and \eqref{eta2levle++} gives
\begin{align} \label{DtauU2}
\|D_\tau \eta_u^2\| = \frac{1}{\tau} \left\|\eta_{u}^{2}-\eta_{u}^{1}\right\|\leq \frac{1}{\tau} \left(\left\|\eta_{u}^{2}\right\| +\left\|\eta_{u}^{1}\right\| \right)\leq C(\tau^{2}+h^{k+1}).
\end{align}
With \eqref{InverIeq}, \eqref{exactsInfty} and \eqref{eta2levle++}, there exists $h_4>0$ such that when $h<h_4$,
\begin{equation} \label{UH2}
\|u_{h}^{2}\|_{\infty}\leq \|\mathrm{R}_{h}u^{2}\|_{\infty}
+\|\mathrm{R}_{h}u^{2}-u_{h}^{2}\|_{\infty}\leq \|\mathrm{R}_{h}u^{2}\|_{\infty}
+Ch^{-1}\|\eta_{u}^{2}\|\leq D_u + 1.
\end{equation}

Next, we take the difference between $t_{2}$ and $t_{1}$ of \eqref{EquError1} with $n=2$ and $v_{h}=\eta_{\Psi}^{5/2}+\eta_{\Psi}^{1/2}$, which yields 
\begin{equation*} 
\begin{split}	
\|\eta_{\Psi}^{5/2}\|^{2}
-\|\eta_{\Psi}^{1/2}\|^{2}
&\leq \left(\|S_{2}^{2}-S_{2}^{1}\| +\|T_{1}^{2}-T_{1}^{1}\|\right) \|\eta_{\Psi}^{5/2}+\eta_{\Psi}^{1/2}\|.
\end{split} 
\end{equation*} 
Similar to \eqref{mainproof6++}-\eqref{DtS3n+}, by applying \Cref{ZourariResult}, \eqref{u1infbddA}, \eqref{u1infbdd}, \eqref{DtauU2} and \eqref{UH2}, we have
\begin{equation}\label{S21level}	
\begin{split}
\|S_{2}^{2}-S_{2}^{1}\| 
&=\left\|\left(S_{1}^{2}-(\xi_{\Psi}^{5/2}+\xi_{\Psi}^{3/2})\right)
-\left(S_{1}^{1}-(\xi_{\Psi}^{3/2}+\xi_{\Psi}^{1/2})\right)\right\|\\
&\leq\left\|S_{1}^{2}-S_{1}^{1}\right\|
+\left\|\xi_{\Psi}^{5/2}-\xi_{\Psi}^{1/2}\right\|\\
&\leq C\tau(\tau^{2}+h^{k+1}), 
\end{split}
\end{equation}
and
\begin{equation}\label{T12level}	
\|T_{1}^{2}-T_{1}^{1}\|\leq C\tau\|D_{\tau}\eta_{u}^{2}\|+C\tau\|D_{\tau}\eta_{u}^{1}\|+C\tau\|\eta_{u}^{2}\|+C\tau h^{k+1}\leq C(\tau^{2}+h^{k+1}).
\end{equation}
By combining with \eqref{etapsi0}, \eqref{S21level} and \eqref{T12level}, we have 
\begin{equation}\label{Psilevel3}
\begin{split} 
\|\eta_{\Psi}^{5/2}\|
\leq\|\eta_{\Psi}^{1/2}\|+\|S_{2}^{2}-S_{2}^{1}\|
+\|T_{1}^{2}-T_{1}^{1}\|\leq C(\tau^{2}+h^{k+1}).
\end{split}
\end{equation}
With the projection estimate \eqref{Ritz2}, we get 
\begin{equation} 
\|e_{\Psi}^{5/2}\|\leq\|\eta_{\Psi}^{5/2}\|+\|\xi_{\Psi}^{5/2}\|\leq C(\tau^{2}+h^{k+1}).
\end{equation}
Applying the Lemma \ref{LemmaL2} to the \eqref{Error3} with $n=2$, it holds 
\begin{equation} 
\begin{split}
\|e_{\Phi}^{5/2}\|\leq C\|e_{\Psi}^{5/2}\|+Ch^{k+1}\leq C(\tau^{2}+h^{k+1}).
\end{split}
\end{equation}

\textbf{Step 3}.
We assume that the estimates in \eqref{EstimateErrora}-\eqref{EstimateErrorc} hold for $0\leq n \leq m$ with $m\geq 2$ as follows
\begin{align}
&\max_{0\leq n \leq m}\left\|e_u^{n}\right\|\leq C\left(\tau^{2}+h^{k+1}\right), \label{MaxUleveln}\\
&\max_{1\leq n \leq m}\left\|D_\tau \eta_u^{n}\right\|\leq C\left(\tau^{2}+h^{k+1}\right), \label{MaxUlevelnD} \\
&\max_{0\leq n \leq m}\left\|e_\Psi^{n+1/2}\right\|\leq C\left(\tau^{2}+h^{k+1}\right), \label{MaxPsileveln}\\
&\max_{0\leq n\leq m}\left\|e_\Phi^{n+1/2}\right\|\leq C(\tau^{2}+h^{k+1}). \label{MaxPhileveln}
\end{align}
By using \eqref{exactsInfty} and the inverse inequality \eqref{InverIeq},  there exists $h_5>0$ such that when $h<h_5$, it holds for $0\leq n \leq m$,
\begin{align}
&\|u_{h}^{n}\|_{\infty}\leq \|\mathrm{R}_{h}u^{n}\|_{\infty}
+1 \leq D_u+1, \label{uninfbdd}\\
&\|\Psi_{h}^{n+1/2}\|_{\infty}\leq \|\mathrm{R}_{h} \Psi^{n+1/2}\|_{\infty}
+1\leq D_\Psi+1,\label{Psiinfbdd}\\
&\|\Phi_{h}^{n+1/2}\|_{\infty}\leq \|\mathrm{R}_{h}\Phi^{n+1/2}\|_{\infty}
+1 \leq D_\Phi+1.\label{Phiinfbdd}
\end{align}

Next, we establish that the estimates \eqref{MaxUleveln}-\eqref{MaxPhileveln} also hold for $n=m+1$. 
Taking the difference between \eqref{EquError2} at $t_{m+1}$ and $t_{m-1}$ that gives 
\begin{equation*}
\begin{split}	
\mathbf{i}\langle D_{\tau}\eta_{u}^{m+1}-D_{\tau}\eta_{u}^{m-1}, \omega_{h}\rangle
=&A_{0}(\overline{\eta}_{u}^{m+1/2}-\overline{\eta}_{u}^{m-3/2}, \omega_{h}) \\
&+\langle G_{1}^{m+1}-G_{1}^{m-1},\omega_{h}\rangle
+\langle R_{2}^{m+1}-R_{2}^{m-1},\omega_{h}\rangle\\
=&\frac{\tau}{2}A_{0}(D_{\tau}\eta_{u}^{m+1}+2D_{\tau}\eta_{u}^{m}+D_{\tau}\eta_{u}^{m-1}, \omega_{h})\\
&+\langle G_{1}^{m+1}-G_{1}^{m-1},\omega_{h}\rangle
+\langle R_{2}^{m+1}-R_{2}^{m-1},\omega_{h}\rangle,
\end{split}
\end{equation*}  
which can be written pointwisely as
\begin{equation}\label{OperEqU}
\begin{split}	
D_{\tau}\eta_{u}^{m+1}-D_{\tau}\eta_{u}^{m-1} 
&=\mathbf{i}\frac{\tau }{2} \Delta_{h}\left(D_{\tau}\eta_{u}^{m+1}
+2D_{\tau}\eta_{u}^{m} + D_{\tau}\eta_{u}^{m-1}\right)
+\Gamma_{1}^{m+1}+\Gamma_{2}^{m+1},
\end{split}
\end{equation}  
where $\Gamma_{1}^{m+1}:=-\mathbf{i}P_{h}(R_{2}^{m+1}-R_{2}^{m-1})$, $\Gamma_{2}^{m+1} :=-\mathbf{i}P_{h}(G_{1}^{m+1}-G_{1}^{m-1})$, and $P_{h}:\mathbf{L}^{2}(\Omega)\rightarrow V_{h}^{c}$ denotes the $L^{2}$ projection. By applying $\mathbf{S}_h^{-1}$ to \eqref{OperEqU} and using the operators introduced in \eqref{Shoper1}, \eqref{Shoper2} and \eqref{OperBh}, it follows
\begin{equation}\label{opereqn}	
D_{\tau}\eta_{u}^{m+1}=(\mathbf{B}_{h}-\mathbf{I}_{h})D_{\tau}\eta_{u}^{m}+\mathbf{B}_{h}D_{\tau}\eta_{u}^{m-1}
+\mathbf{S}_{h}^{-1}\sum_{j=1}^{2}\Gamma_{j}^{m+1}.  
\end{equation}
Applying \cref{ZourariResult} to \eqref{opereqn} gives
\begin{equation}\label{Dtaubdd}
\|D_{\tau}\eta_{u}^{m+1}\|+\|D_{\tau}\eta_{u}^{m}\|\leq 2\left\|\mathbf{S}_{h}(D_{\tau}\eta_{u}^{2})\right\| + 2\left\|\mathbf{S}_{h}(D_{\tau}\eta_{u}^{1})\right\| + 2\sum_{n=2}^{m}\left(\|\Gamma_{1}^{n+1}\| +\|\Gamma_{2}^{n+1}\|\right).
\end{equation}

\textbf{Step 4}. In this step, we use the standard integral remainder of Taylor expansion to estimate $\|\Gamma_{1}^{n+1}\|$ and $\|\Gamma_{2}^{n+1}\|$ in \eqref{Dtaubdd} based on the regularity assumption in \eqref{RegularCondition}. By definition, 
\begin{align}
\|\Gamma_1^{n+1}\| \leq & \|R_{2}^{n+1}-R_{2}^{n-1}\| \leq \|R_{1}^{n+1}-R_{1}^{n-1}\|+\|D_{\tau}\xi_{u}^{n+1} -D_{\tau}\xi_{u}^{n-1}\|,\label{Gm1norm} \\
\|\Gamma_2^{n+1}\| \leq & \|G_{1}^{n+1}-G_{1}^{n-1}\|.\label{Gm2norm}
\end{align}
We first estimate $\|\Gamma_1^{n+1}\|$.
From \eqref{R1form}, we obtain
\begin{equation}\label{appendEq2}
\begin{split}
    \|R_{1}^{n+1}-R_{1}^{n}\|&\leq\left\|(u_{t}^{n+1/2}-D_{\tau}u^{n+1})-(u_{t}^{n-1/2}-D_{\tau}u^{n})\right\|\\
    &\quad +\left\|\Delta(\overline{u}^{n+1/2}-u^{n+1/2}-\overline{u}^{n-1/2}+u^{n-1/2})\right\|\\
    &\quad +\left\|(\Phi^{n+1/2}+V(x)+\Psi^{n+1/2})(u^{n+1/2}-\overline{u}^{n+1/2})\right. \\
     &\quad\quad-\left.
    (\Phi^{n-1/2}+V(x)+\Psi^{n-1/2})(u^{n-1/2}-\overline{u}^{n-1/2})\right\|.
   \end{split}
\end{equation}
Next, we apply the Taylor expression to each term in \eqref{appendEq2} at $t_{n}$. For the first term, it follows 
\begin{equation}\label{appendEq3}
    \begin{split}
    &\left\|(u_{t}^{n+1/2}-D_{\tau}u^{n+1})-(u_{t}^{n-1/2}-D_{\tau}u^{n})\right\|\\
    &\quad= \left\|u_{t}^{n+1/2}-u_{t}^{n-1/2}-\frac{1}{\tau}\left(u^{n+1}-2u^{n}+u^{n-1}\right)\right\|\\
    & \quad\leq\left\|\frac{1}{2!}\int_{t_{n}}^{t_{n+1/2}}(t_{n+1/2}-t)^{2}u_{tttt}(t)dt-\frac{1}{2!}\int_{t_{n}}^{t_{n-1/2}}(t_{n-1/2}-t)^{2}u_{tttt}(t)dt \right.\\
    & \qquad-\left.\frac{1}{3!}\times\frac{1}{\tau}\int_{t_{n}}^{t_{n+1}}(t_{n+1}-t)^{3}u_{tttt}(t)dt -\frac{1}{3!}\times\frac{1}{\tau}\int_{t_{n}}^{t_{n-1}}(t_{n-1}-t)^{3}u_{tttt}(t)dt\right\|\\
    &\quad\leq\left\| \frac{\tau^{3}}{16} \int_{0}^{1}\left(1-s\right)^{2}u_{tttt}\left(t_{n}+\frac{\tau}{2}s\right)ds+\frac{\tau^{3}}{16}\int_{0}^{1}\left(1-s\right)^{2}u_{tttt}\left(t_{n}-\frac{\tau}{2}s\right)ds\right.\\
    &\qquad\left.-\frac{\tau^{3}}{6}\int_{0}^{1}(1-s)^{3}u_{tttt}(t_{n}+\tau s)ds - \frac{\tau^{3}}{6}\int_{0}^{1}(1-s)^{3}u_{tttt}(t_{n}-\tau s)ds\right\|\leq C\tau^{3}.
  \end{split}
\end{equation} 
For the second term, it holds
\begin{equation}\label{appendEq4}
    \begin{split}
    &\left\|\Delta(\overline{u}^{n+1/2}-u^{n+1/2}-\overline{u}^{n-1/2}+u^{n-1/2})\right\|\\
    & \quad\leq \left\|\frac{1}{2}\times\frac{1}{2!}\int_{t_{n}}^{t_{n+1}}(t_{n+1}-t)^{2}u_{ttt}(t)dt -\frac{1}{2!}\int_{t_{n}}^{t_{n+1/2}}(t_{n+1/2}-t)^{2}u_{ttt}(t)dt\right.\\
    &\qquad \left. +\frac{1}{2}\times\frac{1}{2!}\int_{t_{n}}^{t_{n-1}}(t_{n-1}-t)^{2}u_{ttt}(t)dt +\frac{1}{2!}\int_{t_{n}}^{t_{n-1/2}}(t_{n-1/2}-t)^{2}u_{ttt}(t)dt\right\|_{H^{2}}
    \leq C\tau^{3}.
    \end{split}
\end{equation} 
For the third item, it follows
\begin{equation}\label{appendEq4a}
\begin{split}
& \left\|(\Phi^{n+1/2}+V+\Psi^{n+1/2})(u^{n+1/2}-\overline{u}^{n+1/2})-(\Phi^{n-1/2}+V+\Psi^{n-1/2})(u^{n-1/2}-\overline{u}^{n-1/2})\right\|\\
&\quad\leq \left\|\Phi^{n+1/2}+V(x)+\Psi^{n+1/2}\right\|_{\infty}\left\|u^{n+1/2}-\overline{u}^{n+1/2}-u^{n-1/2}+\overline{u}^{n-1/2}\right\|\\
&\quad + \left\|(\Phi^{n+1/2}-\Phi^{n-1/2})+(\Psi^{n+1/2}-\Psi^{n-1/2})\right\|_{\infty}\left\|u^{n-1/2}-\overline{u}^{n-1/2}\right\|.
\end{split}
\end{equation}
Similar to \eqref{appendEq4}, it holds
\begin{align}
&\left\|u^{n-1/2}-\overline{u}^{n-1/2}\right\|\leq C\tau^{2}, \label{appendEq5+}\\
&\left\|u^{n+1/2}-\overline{u}^{n+1/2}-u^{n-1/2}+\overline{u}^{n-1/2}\right\|\leq C\tau^{3}.\label{appendEq5}
\end{align}
In addition, Taylor's theorem and the regularity assumption \eqref{RegularCondition} imply
\begin{equation}\label{appendEq5++}
\begin{aligned}
    & \left\|(\Phi^{n+1/2}-\Phi^{n-1/2})+(\Psi^{n+1/2}-\Psi^{n-1/2})\right\|_{\infty} =\left\| \int_{t_{n-1/2}}^{t_{n+1/2}}\Phi_{t}(s)ds + \int_{t_{n-1/2}}^{t_{n+1/2}}\Psi_{t}(s)ds \right\|_{\infty} \leq C\tau.
\end{aligned}
\end{equation}
Therefore, using \eqref{appendEq3}-\eqref{appendEq5++} and the regularity assumption \eqref{exactsInfty}, we conclude
\begin{equation}\label{appendEq5a}
\left\|R_{1}^{n+1}-R_{1}^{n}\right\|\leq C\tau^{3}. 
\end{equation}  

Moreover, by using the projection error estimate \eqref{Ritz2}, it follows
\begin{equation}\label{appendEq6}
\begin{split}
\|D_{\tau}\xi_{u}^{n+1}-D_{\tau}\xi_{u}^{n-1}\|&=\left\|\mathrm{R}_{h}\left(\frac{u^{n+1}-u^{n}-u^{n-1} +u^{n-2} }{\tau}\right)- \frac{u^{n+1}-u^{n}-u^{n-1} +u^{n-2}}{\tau}\right\|\\
&\leq C\frac{1}{\tau}\int_{0}^{\tau} \left(\int_{t_{n-2}+s}^{t_{n}+s}\left\|\mathrm{R}_{h}u_{tt}(t)-u_{tt}(t) \right\|dt\right)ds \leq C\tau h^{k+1},
\end{split}
\end{equation}
where we have used
\begin{equation}
u^{n+1}-u^{n}-u^{n-1}+u^{n-2} = \int_{0}^{\tau}\left(\int_{t_{n-2}+s}^{t_{n}+s}u_{tt}(t)dt\right)ds.  
\end{equation} 
\eqref{appendEq5a} and \eqref{appendEq6} together with \eqref{Gm1norm} imply
\begin{equation}\label{OperEqU6}
\|\Gamma_1^{n+1}\|\leq \|R_{1}^{n+1}-R_{1}^{n}\|+ \|R_{1}^{n}-R_{1}^{n-1}\|+ \|D_{\tau}\xi_{u}^{n+1}-D_{\tau}\xi_{u}^{n-1}\|\leq C\tau(\tau^2+h^{k+1}).
\end{equation}

Next, we estimate $\|\Gamma_2^{n+1}\|$. From \eqref{G1n1}, it follows
\begin{equation}\label{appendEq9}
\begin{split}
\left\|G_{1}^{n+1}-G_{1}^{n-1}\right\|
\leq&\left\|(\Psi^{n+1/2}\overline{u}^{n+1/2}-\Psi_{h}^{n+1/2}\overline{u}_{h}^{n+1/2})-(\Psi^{n-3/2}\overline{u}^{n-3/2}-\Psi_{h}^{n-3/2}\overline{u}_{h}^{n-3/2})\right\|\\
&+\left\|V(x)\right\|_{\infty}\left\|(\overline{u}^{n+1/2}-\overline{u}_{h}^{n+1/2})-(\overline{u}^{n-3/2}-\overline{u}_{h}^{n-3/2})\right\|\\
&+\left\|(\Phi^{n+1/2}\overline{u}^{n+1/2}-\Phi_{h}^{n+1/2}\overline{u}_{h}^{n+1/2})-(\Phi^{n-3/2}\overline{u}^{n-3/2}-\Phi_{h}^{n-3/2}\overline{u}_{h}^{n-3/2})\right\|.
\end{split}
\end{equation}
For the first term in \eqref{appendEq9},
\begin{equation}\label{appendEq11}
\begin{split}
&\left\|(\Psi^{n+1/2}\overline{u}^{n+1/2}-\Psi_{h}^{n+1/2}\overline{u}_{h}^{n+1/2})-(\Psi^{n-3/2}\overline{u}^{n-3/2}-\Psi_{h}^{n-3/2}\overline{u}_{h}^{n-3/2})\right\|\\
&\leq\left\|(\Psi^{n+1/2}-\Psi^{n-3/2})(\overline{u}^{n-3/2}-\overline{u}_{h}^{n-3/2})\right\| +\left\|(\Psi^{n+1/2}-\Psi_{h}^{n+1/2})(\overline{u}^{n+1/2}-\overline{u}^{n-3/2})\right\|\\
&\quad+\left\|\overline{u}_{h}^{n-3/2}( \Psi^{n+1/2}-\Psi^{n-3/2}-\Psi_{h}^{n+1/2}+\Psi_{h}^{n-3/2})\right\|\\
&\quad+\left\|\Psi_h^{n+1/2}\left(\overline{u}^{n+1/2}-\overline{u}^{n-3/2}-\overline{u}_{h}^{n+1/2}+\overline{u}_{h}^{n-3/2}\right)\right\| := K_1 + K_2 +K_3+K_4. 
\end{split}
\end{equation}
By the Taylor expansion, the split \eqref{eudecomp}, \eqref{eudecompPsi}, and the projection errors, it is easy to obtain
\begin{align}
K_1 &\leq\frac{C\tau}{2}\left\|u^{n-1}+u^{n-2}-u_{h}^{n-1}-u_{h}^{n-2}\right\|\leq C\tau\left(\|\eta_{u}^{n-1}\|+ \|\eta_{u}^{n-2}\|\right)+C\tau h^{k+1}, \label{appendEq12}\\ 
K_2 &\leq\left\|\Psi^{n+1/2}-\Psi_{h}^{n+1/2}\right\|\left\|\frac{(u^{n+1}+u^{n})-(u^{n-1}+u^{n-2})}{2}\right\|\leq C\tau\left(\tau^2+h^{k+1}\right), \label{appendEq13} 
\end{align}
where we have used \eqref{MaxPsileveln}.
By \eqref{uninfbdd}, and using the mean value theorem, it holds
\begin{equation}\label{appdexEUthr}
\begin{aligned}
K_3 \leq & C\|e_\Psi^{n+1/2}-e_\Psi^{n-3/2}\| \leq C\|\eta_\Psi^{n+1/2}-\eta_\Psi^{n-3/2}\| + {C}{\tau} \left\|\frac{\xi_\Psi^{n+1/2}-\xi_\Psi^{n-3/2}}{\tau}\right\| \\
\leq & C\|\eta_\Psi^{n+1/2}-\eta_\Psi^{n-3/2}\| + C\tau h^{k+1}.
\end{aligned}
\end{equation}
Then, taking the difference of \eqref{EquError1} between two time levels and using \eqref{eudecompPsi} yields  
\begin{equation}\label{appdexEpsi}
\begin{split}
(\eta_{\Psi}^{n+1/2}-\eta_{\Psi}^{n-3/2}, v_{h})=(S_{2}^{n}-S_{2}^{n-1},v_{h}) + (T_{1}^{n}-T_{1}^{n-1}, v_{h}).
\end{split}
\end{equation}
Similar to \eqref{mainproof6++} and \eqref{DtS3n+}, it follows
\begin{equation}\label{mainproof6}
\begin{split}	
    \|T_{1}^{n}-T_{1}^{n-1}\|\leq C\tau\|D_{\tau}\eta_{u}^{n}\|+C\tau\|\eta_{u}^{n}\| +C\tau h^{k+1},
\end{split}
\end{equation}	
\begin{equation}\label{DtS3n}
    \left\|S_{2}^{n}-S_{2}^{n-1}\right\|\leq C\tau(\tau^{2}+h^{k+1}).
\end{equation}
By setting $v_{h}=\eta_{\Psi}^{n+1/2}-\eta_{\Psi}^{n-3/2}$ in \eqref{appdexEpsi}, applying Cauchy-Schwartz inequality, and using \eqref{mainproof6} and \eqref{DtS3n},
\begin{equation}\label{appdexEpsib}
\begin{split}
\|\eta_{\Psi}^{n+1/2}-\eta_{\Psi}^{n-3/2}\| 
&\leq\|S_{2}^{n}-S_{2}^{n-1}\|
+\|T_{1}^{n}-T_{1}^{n-1}\| \leq C\tau\left(\|D_{\tau}\eta_{u}^{n}\|
+\|\eta_{u}^{n}\|\right) 
+C\tau(\tau^{2}+h^{k+1}).  
\end{split}
\end{equation}
Plugging \eqref{appdexEpsib} into \eqref{appdexEUthr} gives
\begin{equation}\label{appdexEUthr2}
\begin{split}
K_3 &\leq C\tau\left(\|D_{\tau}\eta_{u}^{n}\|
+\|\eta_{u}^{n}\|\right) 
+C\tau(\tau^{2}+h^{k+1}).
\end{split}
\end{equation}
By \eqref{SREDtauqb1} and \eqref{Psiinfbdd}, $K_4$ in \eqref{appendEq11} gives
\begin{equation}\label{appdexEUfour}
\begin{split}
K_4 \leq& C \left\|(\overline{u}^{n+1/2}-\overline{u}_{h}^{n+1/2})-(\overline{u}^{n-3/2}-\overline{u}_{h}^{n-3/2})\right\| \\
\leq & C\tau\left\| D_\tau e_u^{n+1} +2D_\tau e_u^n + D_\tau e_u^{n-1}\right\| \\
\leq & C\tau\left(\left\| D_\tau \eta_u^{n+1}\right\| +\left\|D_\tau \eta_u^n\right\| + \left\|D_\tau \eta_u^{n-1}\right\| \right)+ C\tau h^{k+1}. 
\end{split}
\end{equation}
Plugging \eqref{appendEq12}, \eqref{appendEq13}, \eqref{appdexEUthr2} and \eqref{appdexEUfour} into \eqref{appendEq11} implies
\begin{equation}\label{appendEq14a}
\begin{split}
&\left\|(\Psi^{n+1/2}\overline{u}^{n+1/2}-\Psi_{h}^{n+1/2}\overline{u}_{h}^{n+1/2})-(\Psi^{n-3/2}\overline{u}^{n-3/2}-\Psi_{h}^{n-3/2}\overline{u}_{h}^{n-3/2})\right\| \\
& \quad \leq C\tau \left( \|D_{\tau}\eta_{u}^{n+1}\|+ \|D_{\tau}\eta_{u}^{n}\|+\|D_{\tau}\eta_{u}^{n-1}\| +\|\eta_{u}^{n}\| +\|\eta_{u}^{n-1}\|+\|\eta_{u}^{n-2}\| \right)+ C\tau(\tau^{2}+h^{k+1}).  
\end{split}
\end{equation}

Similar to $K_4$ in \eqref{appdexEUfour}, the second term in \eqref{appendEq9} yields
\begin{equation}\label{appendEq10}
\begin{split}
& \left\|\overline{e}_u^{n+1/2} - \overline{e}_u^{n-3/2} \right\|
\leq  C\tau\left(\left\| D_\tau \eta_u^{n+1}\right\| +\left\|D_\tau \eta_u^n\right\| + \left\|D_\tau \eta_u^{n-1}\right\| \right)+ C\tau h^{k+1}. 
\end{split}
\end{equation} 

Similar to \eqref{appendEq11}, the estimate of the third term in \eqref{appendEq9}  is given by
\begin{equation}\label{appendEq14}
\begin{split}
&\left\|(\Phi^{n+1/2}\overline{u}^{n+1/2}-\Phi_{h}^{n+1/2}\overline{u}_{h}^{n+1/2})-(\Phi^{n-3/2}\overline{u}^{n-3/2}-\Phi_{h}^{n-3/2}\overline{u}_{h}^{n-3/2})\right\| \\
&\quad\leq C\tau(\tau^{2}+h^{k+1}) + C\tau \left(\|D_{\tau}\eta_{u}^{n+1}\|+ \|D_{\tau}\eta_{u}^{n}\|+\|D_{\tau}\eta_{u}^{n-1}\| +\|\eta_{u}^{n}\| +\|\eta_{u}^{n-1}\|+\|\eta_{u}^{n-2}\| \right).
\end{split}
\end{equation}
Thereby, by \eqref{appendEq14a}, \eqref{appendEq10} and \eqref{appendEq14}, it holds
\begin{equation}\label{OperEqU4}
\begin{split}
\|\Gamma_2^{n+1}\| \leq & \left\|G_{1}^{n+1}-G_{1}^{n-1}\right\|\leq C\tau\left(\|\eta_{u}^{n}\|+\|\eta_{u}^{n-1}\|+\|\eta_{u}^{n-2}\|\right)\\
&+C\tau\left(\|D_{\tau}\eta_{u}^{n+1}\|+\|D_{\tau}\eta_{u}^{n}\|+\|D_{\tau}\eta_{u}^{n-1}\|\right)+C\tau(\tau^{2}+h^{k+1}).
\end{split}
\end{equation}

\textbf{Step 5}. 
In this step, we show that the estimates \eqref{MaxUleveln} and \eqref{MaxUlevelnD} hold for $n=m+1$, that is
\begin{align}
&\left\|e_u^{m+1}\right\|\leq C\left(\tau^{2}+h^{k+1}\right), \label{MaxUleveln-}\\
&\left\|D_\tau \eta_u^{m+1}\right\|\leq C\left(\tau^{2}+h^{k+1}\right). \label{MaxUlevelnD-}
\end{align}

Taking $n=0$ in \eqref{EquError2} and using \eqref{DiscreteLaplace} and $\eta_{u}^{0}=0$ yield
\begin{equation}\label{OperEqU3a-}
\begin{split}
\frac{1}{\tau}\langle \mathbf{S}_h \eta_{u}^{1},\omega_{h}\rangle
= -\mathbf{i} \langle G_{1}^{1},\omega_{h}\rangle
-\mathbf{i}\langle R_{2}^{1},\omega_{h}\rangle,
\end{split}
\end{equation}
which by taking $\omega_{h} = \mathbf{S}_h \eta_{u}^{1}$ in \eqref{OperEqU3a-} and using the estimates \eqref{SREqb2}, \eqref{ThmEq2b}, \eqref{etau1level} yields
\begin{equation}\label{OperEqU3a+}
\begin{split}
\|\mathbf{S}_{h}\eta_{u}^{1}\| \leq \tau\left(\|G_{1}^{1}\|+\|R_{2}^{1}\|\right)\leq C\tau(\tau^{2}+h^{k+1}).
\end{split}
\end{equation}
Using $\eta_{u}^{0}=0$ again gives
\begin{equation}\label{OperEqU3a}
\begin{split}
\|\mathbf{S}_{h}(D_{\tau}\eta_{u}^{1})\|=\frac{1}{\tau}\|\mathbf{S}_{h}\eta_{u}^{1}\|\leq C(\tau^{2}+h^{k+1}).
\end{split}
\end{equation}
Moreover, it also holds 
\begin{equation}\label{OperEqU3b+}
\begin{split}
\|\mathbf{S}_{h}\eta_{u}^{2}\|&\leq\|\mathbf{T}_{h}\eta_{u}^{1}\|
+C\tau\left(\|G_{1}^{2}\|+\|R_{2}^{2}\|\right)\leq\|(2\mathbf{I}_{h}-\mathbf{S}_{h})\eta_{u}^{1}\| +C\tau\left(\|G_{1}^{2}\|+\|R_{2}^{2}\|\right)\\
&\leq 2\|\eta_{u}^{1}\|+\|\mathbf{S}_{h}\eta_{u}^{1}\| + C\tau(\tau^{2}+h^{k+1}), 
\end{split}
\end{equation}
where we have used \eqref{SREqb2}, \eqref{etau1level}, \eqref{ThmEq3++}, \eqref{OperEqU3a+} and
\begin{equation}\label{leveG12}
\begin{split}
\left\|G_{1}^{2}\right\|
&\leq \|\overline{u}^{3/2}\|_{\infty}\left(\|e_{\Phi}^{3/2}\|+\|e_{\Psi}^{3/2}\|\right) 
+\left(\|\Phi_{h}^{3/2}\|_{\infty}+\|V(x)\|_{\infty}
+\|\Psi_{h}^{3/2}\|_{\infty}\right)\|\overline{e}_{u}^{3/2}\|\\
&\leq C\left(\|\eta_{u}^{2} \|+\|\eta_{u}^{1}\|\right)
+C\left(\tau^{2}+h^{k+1}\right). 
\end{split}
\end{equation}
Similar to \eqref{ThmEq5}, by using \eqref{SREqb2} and \eqref{leveG12}, it holds 
\begin{equation}\label{2TauU2}
\begin{split}
\frac{1}{2\tau}\left(\|\eta_{u}^{2}\|^{2}
-\|\eta_{u}^{1}\|^{2}\right) 
& =\mathrm{Im}\left\langle G_{1}^{2}+R_{2}^{2},\overline{\eta}_{u}^{3/2}\right\rangle \leq \left\|G_{1}^{2}+R_{2}^{2}\right\|\left\|\overline{\eta}_{u}^{3/2}\right\|\\
&\leq\frac{1}{2}\left\|G_{1}^{2}+R_{2}^{2}\right\|^2 +\frac{1}{2}\left\|\overline{\eta}_{u}^{3/2}\right\|^2\\
&\leq \frac{1}{2}\left(\left\|G_{1}^{2}\right\|^2 +\left\|R_{2}^{2}\right\|^2 + 2\left\|G_{1}^{2}\right\|\left\|R_{2}^{2}\right\| \right)+\frac{1}{8}\left\|\eta_{u}^{2}+\eta_{u}^{1}\right\|^2\\
&\leq\left\|G_{1}^{2}\right\|^2+ \left\|R_{2}^{2}\right\|^2 + \frac{1}{4}\left(\|\eta_{u}^{2}\|^{2}+\|\eta_{u}^{1}\|^{2}\right)\\
&\leq C\left(\|\eta_{u}^{2}\|^{2}+\|\eta_{u}^{1}\|^{2}\right)+C(\tau^{2}+h^{k+1})^{2}.
\end{split}
\end{equation}
As long as $\tau<\tau_{3}:=\min\{\tau_1, 1/(2C)\}$, plugging \eqref{etau1level} into \eqref{2TauU2} implies
\begin{equation}\label{eta2level2}	 
	\|\eta_{u}^{2}\|\leq C\tau(\tau^{2}+h^{k+1}).
\end{equation}
 Then, it holds
\begin{equation}\label{OperEqU3b}
\begin{split}
\|\mathbf{S}_{h}\eta_{u}^{2}\| \leq C\tau(\tau^{2}+h^{k+1}).
\end{split}
\end{equation}
By using \eqref{OperEqU3a+} and \eqref{OperEqU3b}, 
\begin{equation}\label{OperEqU3bcc}
\begin{split}
\|\mathbf{S}_{h}(D_{\tau}\eta_{u}^{2})\|
&\leq\frac{1}{\tau}\left(\|\mathbf{S}_{h}\eta_{u}^{2}\|+\|\mathbf{S}_{h}\eta_{u}^{1}\|\right)\leq C(\tau^{2}+h^{k+1}).
\end{split}
\end{equation}
Plugging \eqref{OperEqU6}, \eqref{OperEqU4}, \eqref{OperEqU3a} and \eqref{OperEqU3bcc} into \eqref{Dtaubdd} and using initial estimates in Step 1 yield
\begin{equation}\label{OperEqU3}
\begin{split}
&\|D_{\tau}\eta_{u}^{m+1}\|+\|D_{\tau}\eta_{u}^{m}\|\leq C\tau \sum_{n=1}^{m} \left(\|\eta_{u}^{n}\|+\|D_{\tau}\eta_{u}^{n+1}\|+ \|D_{\tau}\eta_{u}^{n}\|\right)+C(\tau^{2}+h^{k+1}).   
\end{split}
\end{equation}

Setting $\omega_{h}=\overline{\eta}_{u}^{n+1/2}$ in \eqref{EquError2}, and taking its imaginary part give
\begin{equation}\label{etauNlevel}
\begin{split}
\frac{1}{2\tau}\left(\|\eta_{u}^{n+1}\|^{2}-\|\eta_{u}^{n}\|^{2}\right)
&=\mathrm{Im}\left\langle G_{1}^{n+1},\overline{\eta}_{u}^{n+1/2}\right\rangle+\mathrm{Im} \left\langle R_{2}^{n+1},\overline{\eta}_{u}^{n+1/2}\right\rangle\\
&\leq\frac{1}{2}\left\|G_{1}^{n+1}\right\| \|\eta_{u}^{n+1}+\eta_{u}^{n}\|
+\frac{1}{2}\left\|R_{2}^{n+1}\right\|
\|\eta_{u}^{n+1}+\eta_{u}^{n}\|.
\end{split}
\end{equation}
Similar to \eqref{ThmEq2b}, we have
\begin{equation}\label{G1m}
\begin{split}
\|G_{1}^{n+1}\|
\leq & \left(\|\Phi_{h}^{n+1/2}\|_{\infty}+\|V(x)\|_{\infty}+\|\Psi_{h}^{n+1/2}\|_{\infty}\right)\|\overline{e}_{u}^{n+1/2}\|\\
&+\|\overline{u}^{n+1/2}\|_{\infty}\left(\|e_{\Phi}^{n+1/2}\|+\|e_{\Psi}^{n+1/2}\|\right)
\leq C\left(\|\eta_{u}^{n+1}\|+\|\eta_{u}^{n}\|\right)+C(\tau^{2}+h^{k+1}),
\end{split}
\end{equation}
where we have used the boundedness \eqref{Psiinfbdd} and \eqref{Phiinfbdd}, and the estimates \eqref{MaxPsileveln} and \eqref{MaxPhileveln}.

Applying \eqref{SREqb2} and \eqref{G1m} upon simplification,  \eqref{etauNlevel} gives
\begin{equation}\label{etauNlevelau}
\begin{split}
\|\eta_{u}^{n+1}\|-\|\eta_{u}^{n}\|&\leq \tau \left\|G_{1}^{n+1}\right\|+ \tau\left\|R_{2}^{n+1}\right\|\leq C\tau\left(\|\eta_{u}^{n+1}\|+\|\eta_{u}^{n}\|\right)
+C\tau\left(\tau^{2}+h^{k+1}\right),
\end{split}
\end{equation}
which upon summing up \eqref{etauNlevelau} from $n=1$ to $m$ leads to  
\begin{equation}\label{upsicombine2}
\begin{split}
\|\eta_{u}^{m+1}\| 
\leq & \|\eta_{u}^{1}\|
+C\left(\tau^{2}+h^{k+1}\right)+ C\tau\sum_{n=1}^{m}\left(\|\eta_{u}^{n+1}\|+\|\eta_{u}^{n}\|\right).
\end{split}
\end{equation} 
The summation of \eqref{OperEqU3} and \eqref{upsicombine2} yields
\begin{equation}\label{upsicombine3}
\begin{split}
&\|\eta_{u}^{m+1}\|+  \|D_{\tau}\eta_{u}^{m+1}\|+\|D_{\tau}\eta_{u}^{m}\| \leq C\tau \sum_{n=1}^{m} \left( \|\eta_{u}^{n+1}\| + \|D_{\tau}\eta_{u}^{n+1}\|+ \|D_{\tau}\eta_{u}^{n}\| \right)+C(\tau^{2}+h^{k+1}).    
\end{split}
\end{equation}
By Gronwall's inequality in \Cref{GronwallLem}, there exists $\tau_{4}>0$independent of $m$ such that when $\tau<\tau_{4}$,
\begin{equation}\label{Gronwalluerr}
\begin{split} 
&\|\eta_{u}^{m+1}\|+ \|D_{\tau}\eta_{u}^{m+1}\|+\|D_{\tau}\eta_{u}^{m}\| \leq C\left(\tau^{2}+h^{k+1}\right),
\end{split}
\end{equation}
where $C$ depends on $T$ and is independent of $m$.
The estimate \eqref{Gronwalluerr}, together with the projection error, implies the estimates \eqref{MaxUleveln-} and \eqref{MaxUlevelnD-}.

\textbf{Step 6}.
Last, we show that \eqref{MaxPsileveln} and \eqref{MaxPhileveln} also hold for $n=m+1$.
By \eqref{InverIeq}, \eqref{exactsInfty}, \eqref{upsicombine2}, and $\tau \leq Ch$, there exist $h_6>0$ depending on $T$ but independent of $m$ such that when $h< h_6$, 
\begin{equation*}
\begin{split}
&\|u_{h}^{m+1}\|_{\infty}\leq \|\mathrm{R}_{h}u^{m+1}\|_{\infty}
+Ch^{-1}\|\eta_{u}^{m+1}\|\leq D_u+C_uh\leq D_u+1.
\end{split}
\end{equation*}
Setting $v_{h}=\eta_{\Psi}^{n+1/2}+\eta_{\Psi}^{n-3/2}$ in \eqref{appdexEpsi} gives
\begin{equation}\label{mainproof8a}
\begin{split}	
\|\eta_{\Psi}^{n+1/2}\|
-\|\eta_{\Psi}^{n-3/2}\|
&\leq \|S_{2}^{n}-S_{2}^{n-1}\|
+\|T_{1}^{n}-T_{1}^{n-1}\|\leq C\tau\left(\|D_{\tau}\eta_{u}^{n}\|
+\|\eta_{u}^{n}\|\right)+C\tau (\tau^{2}+h^{k+1}).
\end{split}
\end{equation}
Summing up \eqref{mainproof8a} from $n=1$ to $n=m+1$ gives
\begin{equation}\label{etaerrm1}
\|\eta_{\Psi}^{m+3/2}\|\leq C\tau \sum_{n=1}^{m+1}\left(\|D_{\tau}\eta_{u}^{n}\| 
+\|\eta_{u}^{n}\|\right) + C(\tau^2+h^{k+1}) \leq C(\tau^2+h^{k+1}),
\end{equation}
where we have used \eqref{initerr}, \eqref{etapsi0}, \eqref{MaxUleveln}, \eqref{MaxUlevelnD}, and \eqref{MaxPsileveln} with $n=m$.
The estimate \eqref{etaerrm1} together with the projection error implies
\begin{equation}\label{erropsiL2Phi}
\begin{split}
&\|e_{\Psi}^{m+3/2}\|\leq \|\xi_{\Psi}^{m+3/2}\|+ \|\eta_{\Psi}^{m+3/2}\|\leq C\left(\tau^{2}+h^{k+1}\right).
\end{split}
\end{equation}

Lemma \ref{LemmaL2}, \eqref{erropsiL2Phi}, and \eqref{Error3} further give
\begin{equation*}
\|e_{\Phi}^{m+3/2}\|
\leq C\|e_{\Psi}^{m+3/2}\|+Ch^{k+1}\leq C(\tau^{2}+h^{k+1}).
\end{equation*}

Therefore, the estimates \eqref{MaxUleveln}-\eqref{MaxPhileveln} hold for $n=m+1$, if $\tau_{0}=\max\{\tau_{i}\}_{i=1}^{4}$ and $h_0=\min\{h_j\}_{j=1}^{6}$, which depend on $T$ but are independent of $N$. 
This completes the proof.
\end{proof}


\section{Extension}\label{secExtension}

The model equation \eqref{TargetEq} without the self-repulsion term $|u|^{2}u$ and the external potential will degenerate to the Schr\"{o}dinger-Poisson equation with constant coefficients \cite{Athanassoulis2023novel}
\begin{subequations}\label{CoeffSystem}
	\begin{align}
		&\mathbf{i}u_{t}=-\alpha \Delta u +\beta \Phi u, \quad (x,t) \in \Omega \times (0,T],\\
		&\Delta \Phi = |u|^{2}-c,\quad x\in \Omega,\\
  & u(x,0) = u_{0}(x), \quad x\in \Omega, \\
&u(x,t)=0 \quad \text{and}\quad \Phi(x)=0, \quad x\in \partial\Omega,
	\end{align}
\end{subequations}
where the parameter $\alpha>0$, $\beta \in \mathbb{R}$. 

Introducing an auxiliary variable $\Psi$, the system \eqref{CoeffSystem} can be equivalently expressed as    
\begin{equation}\label{CoeffSystemRe}
	\left\{
	\begin{aligned}
		&\Psi=|u|^{2}, \\
		&\mathbf{i}u_{t}=-\alpha\Delta u+ \beta\Phi u,\\
		&\Delta\Phi=\Psi-c.
	\end{aligned}
	\right.
\end{equation}  
Then, the proposed relaxation finite element method \eqref{fullyVdiscrete} for the nonlinear Schr\"{o}dinger-Poisson equation \eqref{TargetEq} reduces to the Besse-style relaxation Crank-Nicolson finite element method \cite{Athanassoulis2023novel},
\begin{subequations}\label{CoeffSystemFull}
	\begin{align}
		&\left(\Psi_{h}^{n+1/2}+\Psi_{h}^{n-1/2}, v_{h}\right) = \left(2|u_{h}^{n}|^{2}, v_{h}\right), \quad \forall v_{h} \in V_{h}, \\
		&\mathbf{i}\left\langle D_{\tau} u_{h}^{n+1}, \omega_{h}\right\rangle = \alpha A_{0}\left(\overline{u}_{h}^{n+1/2},\omega_{h}\right) +\beta\left\langle \Phi_{h}^{n+1/2} \overline{u}_{h}^{n+1/2}, \omega_{h}\right\rangle , \quad \forall \omega_{h} \in V_{h}^{c},\\
		& A_{1}\left(\Phi_{h}^{n+1/2}, \chi_{h}\right) = -\left((\Psi_{h}^{n+1/2}-c), \chi_{h}\right), \quad \forall \chi_{h} \in V_{h},  
	\end{align}
\end{subequations}
The following results hold for the scheme above.
\begin{lemma}\cite{Athanassoulis2023novel}
For any $\tau>0$, the relaxation Crank-Nicolson finite element method \eqref{CoeffSystemFull} satisfies the discrete conservation for both mass and modified energy with $0\leq n\leq N-1$, respectively
\begin{align}
&M_{h}^{n+1}= M_{h}^{0}, \\
&E_{h}^{n+1}= E_{h}^{0}, 
\end{align}
where the mass
$M_{h}^{n+1}=\int_{\Omega}\left|u_{h}^{n+1}\right|^{2}dx$,
and the modified energy 
$$
E_{h}^{n+1}= \alpha A(u_{h}^{n+1}, u_{h}^{n+1})+\frac{\beta}{2}A(\Phi_{h}^{n+3/2},\Phi_{h}^{n+1/2}).
$$
\end{lemma}

Following the convergence analysis of the proposed scheme \eqref{fullyVdiscrete} for the nonlinear Schr\"{o}dinger-Poisson equation \eqref{TargetEq}, we can extend the current error estimates to 
the scheme \eqref{CoeffSystemFull} for the Schr\"{o}dinger-Poisson equation \eqref{CoeffSystem}. More specifically, we derive the following results. 
\begin{theorem}\label{ThmExtension}
Suppose that $u$, $\Psi$ and $\Phi$ satisfy the regularity conditions \eqref{RegularCondition}. If $\tau\leq Ch$, then there exists a constant $\tau_{0}>0$ and $h_{0}>0$ such that when time step $\tau<\tau_{0}$ and mesh size $h<h_{0}$, the solutions of the relaxation Crank-Nicolson finite element scheme \eqref{CoeffSystemFull} satisfy the following estimates
\begin{align}
&\max_{0\leq n \leq N-1}\left\|e_{u}^{n+1}\right\|\leq C\left(\tau^{2}+h^{k+1}\right),\\
&\max_{0\leq n \leq N-1}\left\|e_{\Psi}^{n+1/2}\right\|\leq C\left(\tau^{2}+h^{k+1}\right),\\
&\max_{0\leq n\leq N-1}\left\|e_{\Phi}^{n+1/2}\right\|\leq C(\tau^{2}+h^{k+1}).
\end{align}
\end{theorem}
The proof is similar to that of Theorem \ref{Thm+}, thus we omit it here.

\begin{remark}
    The proposed method and the error analysis also have the potential to be applied to other types of equations, such as the Gross-Pitaevskii-Poisson equation \cite{Verma2021formation} and the Gross-Pitaevskii-Poisson system \cite{Sakaguchi2020Gross}. The Gross-Pitaevskii-Poisson equation incorporates a nonlocal mean density and additionally conserves the momentum, adding complexity beyond \eqref{TargetEq}, while the Gross-Pitaevskii-Poisson system involves the interaction between positive and negative bosonic ions. We leave these explorations for future work. 
\end{remark}

\section{Numerical experiments}\label{secNmecical}
In this section, we present numerical experiments to validate our theoretical analysis. This includes an examination of the convergence rates and the conservation properties of the relaxation Crank-Nicolson finite element method. All numerical examples are implemented using the FEALPy package \cite{Wei2017FEALPy}.

We consider the two-dimensional Schr\"{o}dinger-Poisson equation on $\Omega=[-8,8]^{2}$,
\begin{equation}\label{SPex}
\begin{aligned}
&\mathbf{i}u_{t}(x_{1},x_{2},t)=-\frac{1}{2}\Delta u + \Phi(x_{1},x_{2},t) u+V(x_{1},x_{2})u+|u|^{2} u, \quad (x_{1},x_{2})\in \Omega, \\
&-\Delta\Phi(x_{1},x_{2},t)= |u|^{2}-1, \quad (x_{1},x_{2})\in \Omega, \\
&u(x_{1},x_{2},t)=0, \quad (x_{1},x_{2})\in \partial\Omega, \\
&\Phi(x_{1},x_{2},t)= 0, \quad (x_{1},x_{2})\in \partial\Omega,\\
&u(x_{1},x_{2},0)= u_{0}(x_{1},x_{2}) = \frac{1}{\sqrt{2\pi}}
e^{\frac{-x_{1}^{2}+x_{2}^{2}}{4}}(x_{1}+\mathbf{i}x_{2}), \quad (x_{1},x_{2}) \in ~\Omega.
\end{aligned}
\end{equation}
Here, we consider three different external potentials $V(x_{1},x_{2})=V_i(x_{1},x_{2})$, $i=0,1,2$ with 
$V_{0}(x_{1},x_{2})=0$, $V_{1}(x_{1},x_{2})=\frac{x_{1}^{2}+x_{2}^{2}}{2}$, and $V_{2}(x_{1},x_{2})=\frac{x_{1}^{2}-x_{2}^{2}}{2}$.

{\bf Test case 1.}
To validate the accuracy and convergence rate of the relaxation Crank-Nicolson finite element method, we take the $\mathbb{Q}^{k}$ polynomials with $k=1,2$. As the exact solution is unavailable, we compute the time discretization errors as 
$\|u_h^{T/\tau}-u_h^{T/(2\tau)}\|$, where $u_h^{T/\tau}$ is finite element solution at $t=T$ with time step $\tau$.
Table \ref{TableTimeOrder} reports the time discretization error in $L^{2}$ norm and the order of accuracy, utilizing a sufficiently small fixed spatial mesh size. Based on the obtained results, it is evident that the proposed method exhibits second-order accuracy in time.

\begin{table}[!ht]
\caption{Time discretization errors with $T=0.1$ and $V(x_{1},x_{2}) = V_{2}(x_{1},x_{2})$.}
\centering
\begin{tabular}{cccccc}
\toprule
$\tau$ & 1.0e-02 & 5.0e-03 & 2.5e-03 \\
\midrule
$\|u_h^{T/\tau}-u_h^{T/(2\tau)}\|$ & 6.3247e-03 & 1.5870e-03 & 3.9710e-04\\
Order & -- & 1.99 & 2.00   \\
\bottomrule
\end{tabular}
\label{TableTimeOrder}
\end{table}

{\bf  Test case 2.}
In Table \ref{TableSpaceOrder}, we compute the spatial discretization errors $\|u_{NC}-u_{2NC}\|$ between the two-level approximations at final time $T=0.1$ with a sufficiently small fixed time step, where $u_{NC}$ denotes the numerical solution on $NC\times NC$ meshes. It is observed that the proposed method demonstrates $(k+1)$th order accuracy in space.
\begin{table}[!ht]
\caption{Spatial discretization errors with $T=0.1$ and $V(x_{1},x_{2})=V_{2}(x_{1},x_{2})$.}
\centering
\begin{tabular}[c]{ccccccc}
\toprule
$\mathbb{Q}^k$ & $\|u_{50}-u_{100}\|$ & Order & $\|u_{100}-u_{200}\|$& Order  & $\|u_{200}-u_{400}\|$ & Order \\
\midrule
$k=1$ & 1.8851e-02 & --  & 4.7363e-03 & 1.99 & 1.1856e-03 & 2.00\\
$k=2$ & 6.4115e-04 & --  & 8.0828e-05 & 2.99 & 1.0144e-05 & 2.99\\
\bottomrule
\end{tabular}
\label{TableSpaceOrder}
\end{table}

{\bf  Test case 3.}
Subsequently, we apply the proposed method using a mesh with $NC=80$ for spatial discretization and a time step of $\tau=2\times 10^{-3}$, based on $\mathbb{Q}^2$ polynomials, to verify the performance of our numerical scheme in preserving mass and energy conservation properties. For $0\leq n\leq N-1$, we define the mass change and energy change as follows:
\begin{equation}\label{mechange}
\text{Mass Change}=\left|\frac{M_{h}(t_{n})-M_{h}(0)}{M_{h}(0)}\right|, \quad \text{Energy Change}=\left|\frac{E_{h}(t_{n})-E_{h}(0)}{E_{h}(0)}\right|.
\end{equation}

\begin{figure}[htbp]
	\centering
	\begin{subfigure}{0.4\linewidth}
		\centering
		\includegraphics[width=1\linewidth]{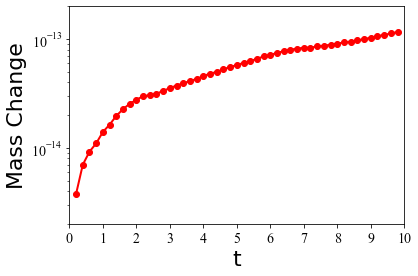}
	\end{subfigure}
	\centering
	\begin{subfigure}{0.4\linewidth}
		\centering
		\includegraphics[width=1\linewidth]{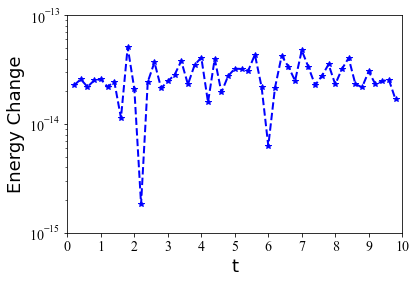}
	\end{subfigure}
	\caption{Evolution of the mass and modified energy with $V(x_{1},x_{2})=V_{0}(x_{1},x_{2})$.}\label{RelativeErrorCase4}
\end{figure}

\begin{figure}[htbp]
	\centering
	\begin{subfigure}{0.4\linewidth}
		\centering
		\includegraphics[width=1\linewidth]{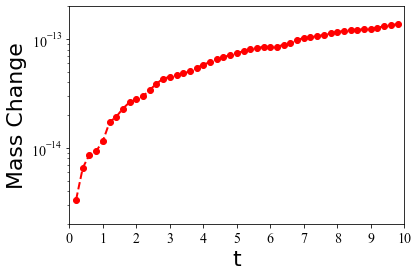}
	\end{subfigure}
	\centering
	\begin{subfigure}{0.4\linewidth}
		\centering
		\includegraphics[width=1\linewidth]{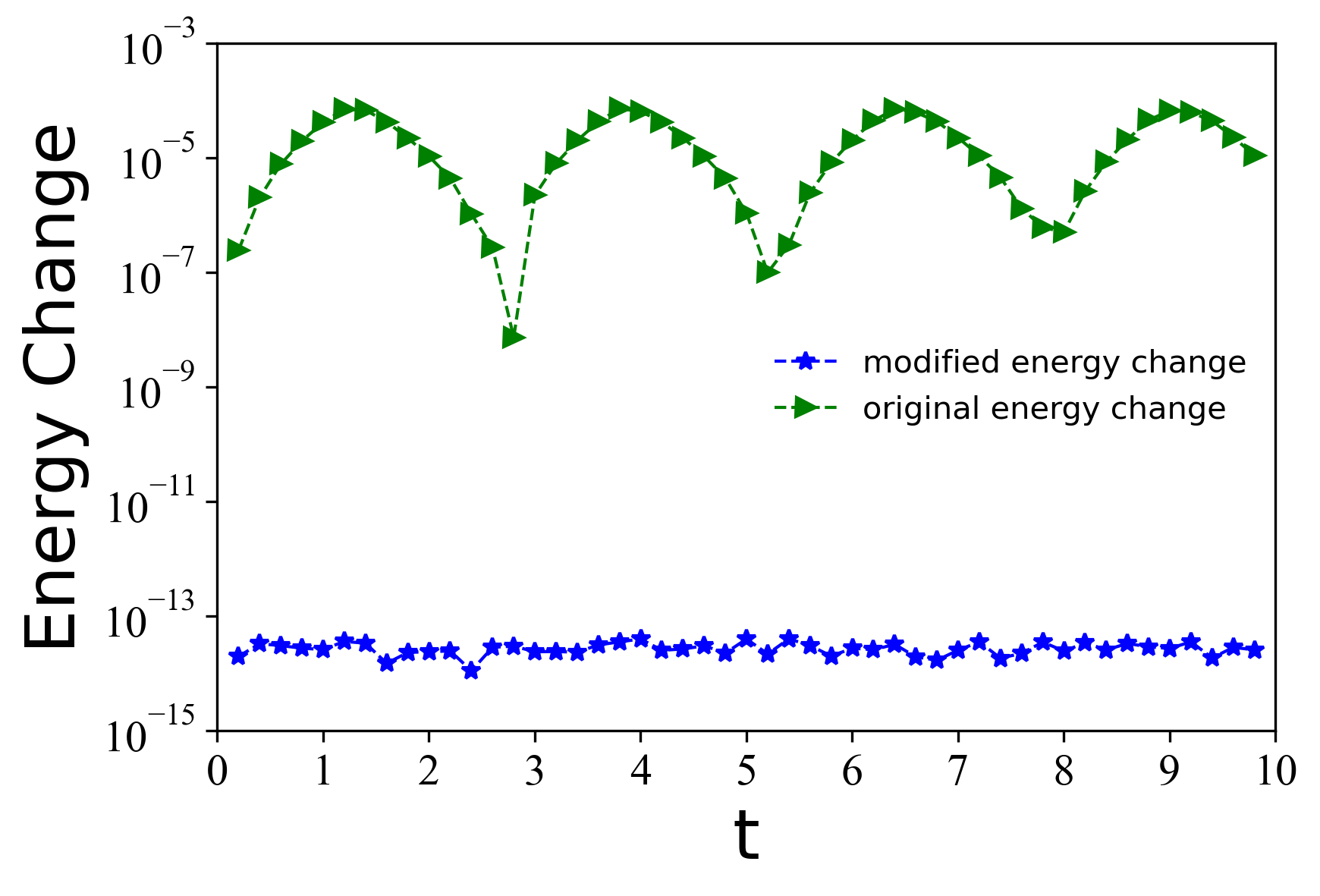}
	\end{subfigure}
	\caption{Evolution of the mass and energy with $V(x_{1},x_{2})=V_{1}(x_{1},x_{2})$ .}
	\label{RelativeErrorCase2a}
\end{figure}

\begin{figure}[htbp]
	\centering
	\begin{subfigure}{0.4\linewidth}
		\centering
		\includegraphics[width=1\linewidth]{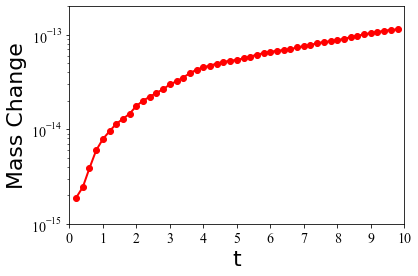}
	\end{subfigure}
	\centering
	\begin{subfigure}{0.4\linewidth}
		\centering
		\includegraphics[width=1\linewidth]{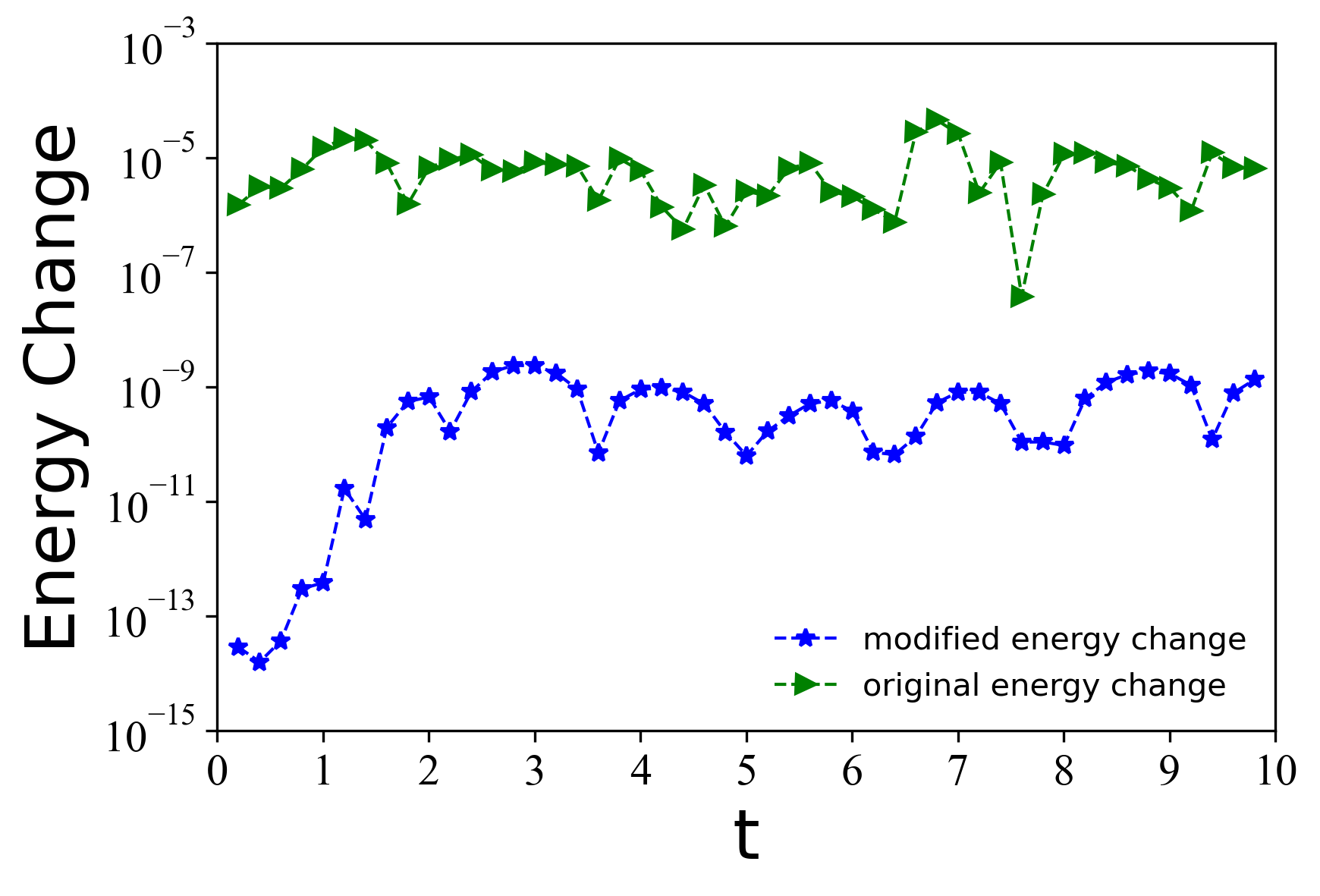}
	\end{subfigure}
	\caption{Evolution of the mass and energy with $V(x_{1},x_{2})=V_{2}(x_{1},x_{2})$ .}
	\label{RelativeErrorCase2b}
\end{figure}

The discrete mass and energy, as defined in \Cref{CoffiConservationLem}, are computed for 
$V(x_{1},x_{2})=V_{i}(x_{1},x_{2})$, with $i=0,1,2$, and the changes in mass and energy are illustrated in \Cref{RelativeErrorCase4} - \Cref{RelativeErrorCase2b}, respectively. 
Although the case with $V(x_{1},x_{2})=V_{2}(x_{1},x_{2})$ shows a relatively larger energy error compared with other cases, as seen in \Cref{RelativeErrorCase2b}, the results suggest that both mass and modified energy are well preserved at the discrete level for all cases.

For cases of $V(x_{1},x_{2})=V_1(x_{1},x_{2})$ and $V(x_{1},x_{2})=V_2(x_{1},x_{2})$, we also compute a direct approximation of the original energy in \eqref{inroeng} at $t_n$, defined as
\begin{equation}\label{Energyti2a}
 \tilde{E}_{h}^n:=\int_{\Omega}\left(\frac{1}{2}|\nabla u^{n}_{h}|^{2}+\frac{1}{2\mu}|\nabla \bar{\Phi}^{n}_{h}|^{2}+ V(x)|u^{n}_{h}|^{2}+\frac{1}{2}|u^{n}_{h}|^{4}\right)dx, \quad 0\leq n \leq N,
\end{equation}
where
\[
\bar{\Phi}^{n}_{h} = \frac{{\Phi}^{n+1/2}_{h}+{\Phi}^{n-1/2}_{h}}{2}.
\]
For both cases, the changes in the approximated original energy $\tilde{E}_{h}^n$, defined similarly to the energy change in \eqref{mechange}, are also shown in \Cref{RelativeErrorCase2a} and \Cref{RelativeErrorCase2b}. Although the changes in the directly approximated original energy are relatively larger than those of the modified energy, the original energy remains well preserved in both cases.

\begin{figure}[htbp]
	\centering
	\begin{subfigure}{0.325\linewidth}
		\centering
		\includegraphics[width=1.16\linewidth]{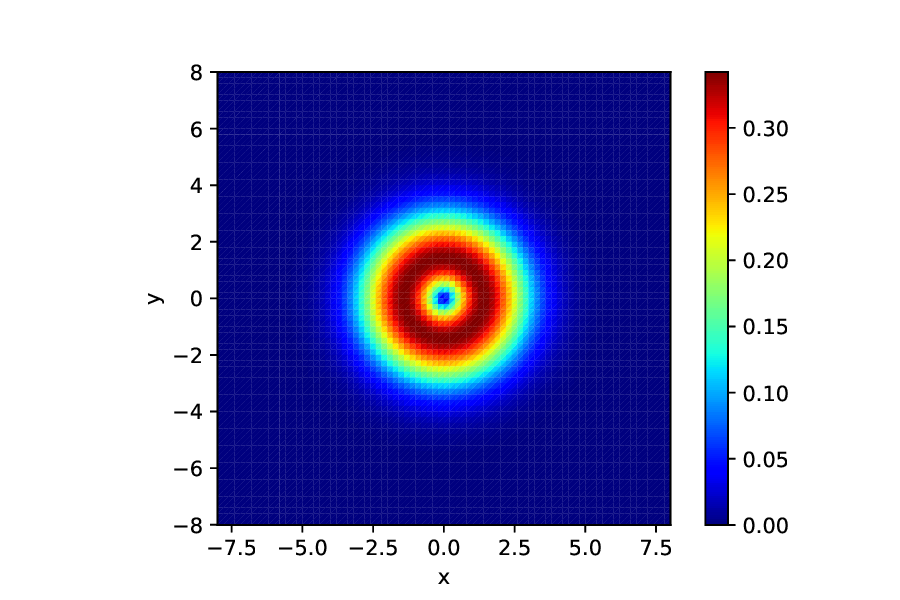}
		\caption{$t=0$}
	\end{subfigure}
	\centering
	\begin{subfigure}{0.325\linewidth}
		\centering
		\includegraphics[width=1.16\linewidth]{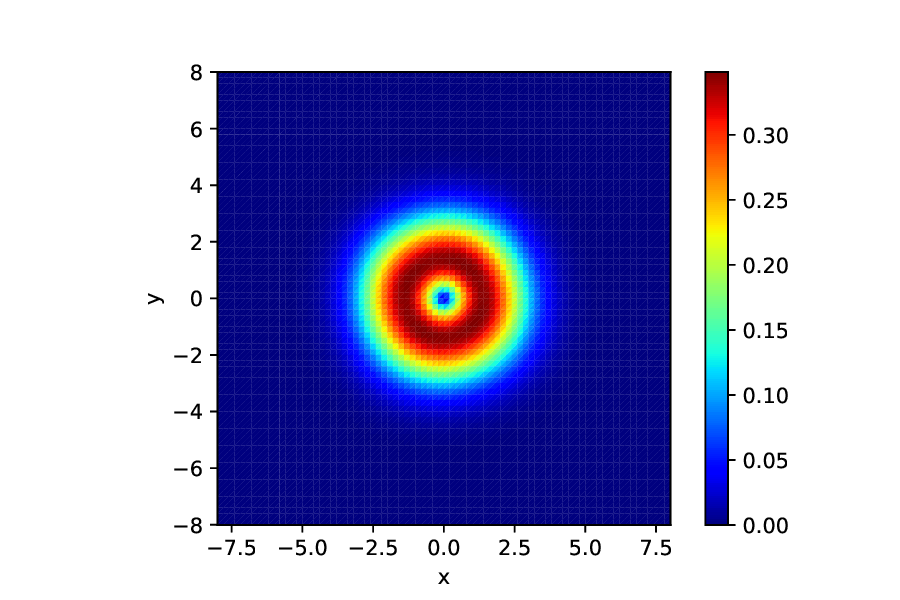}
		\caption{$t=5$}
	\end{subfigure}
	\centering
	\begin{subfigure}{0.325\linewidth}
		\centering
		\includegraphics[width=1.16\linewidth]{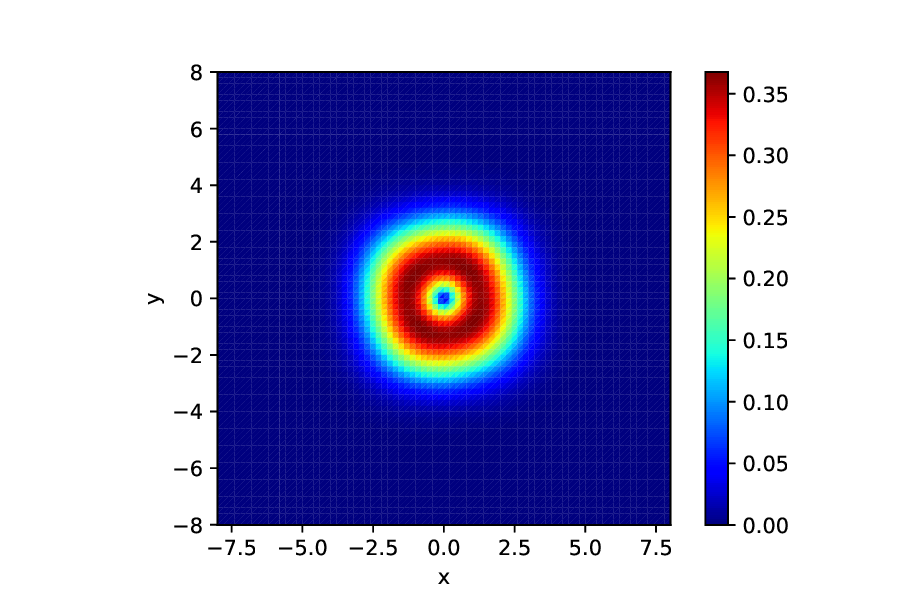}
		\caption{$t=10$}
	\end{subfigure}
	\caption{The patterns evolution of the wave function $|u(x_{1},x_{2},t)|$ with $V(x_{1},x_{2})=V_{0}(x_{1},x_{2})$.}
	\label{Case4Solution}
\end{figure}

\begin{figure}[htbp]
	\centering
	\begin{subfigure}{0.325\linewidth}
		\centering
		\includegraphics[width=1.16\linewidth]{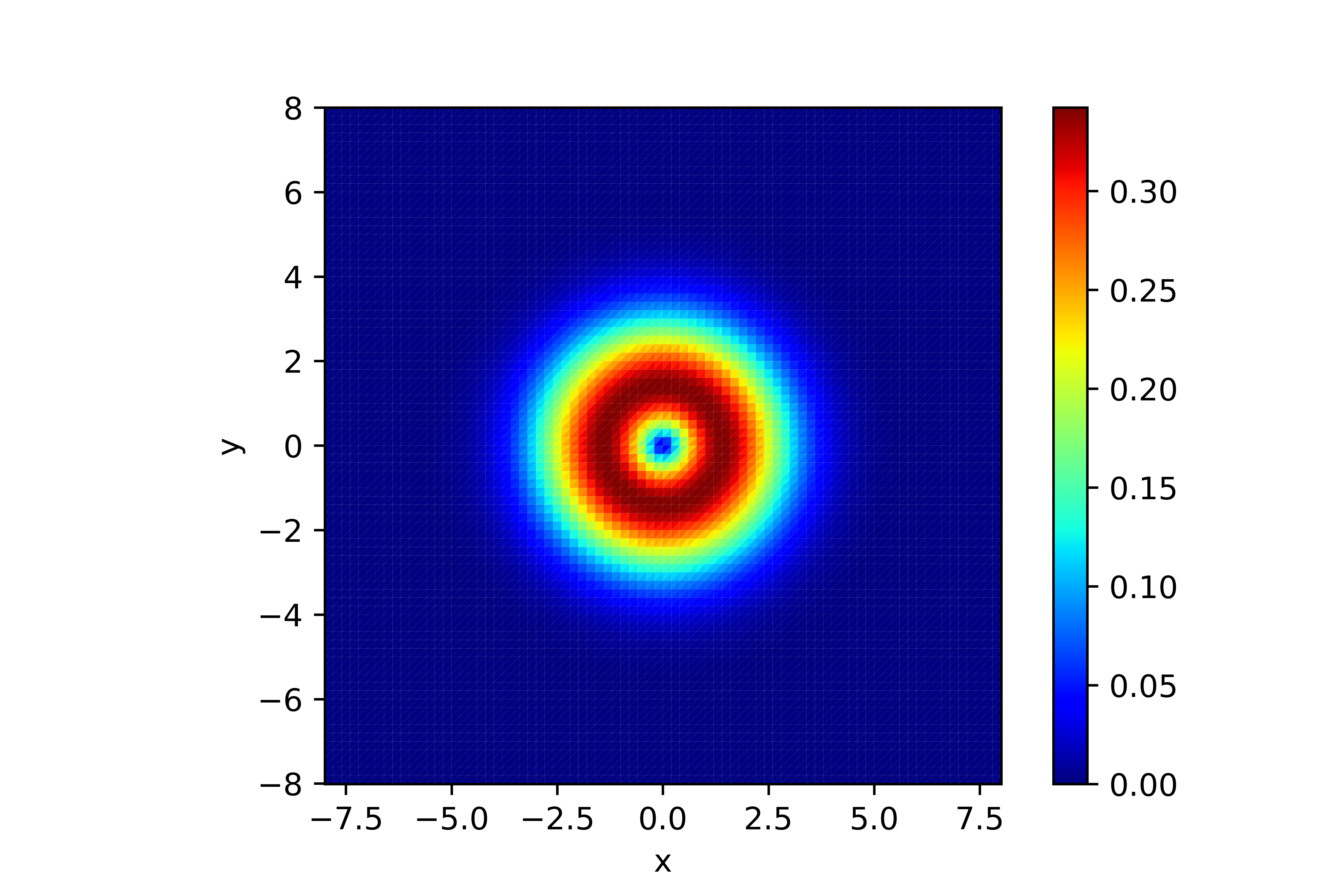}
		\caption{$t=0$}
	\end{subfigure}
	\centering
	\begin{subfigure}{0.325\linewidth}
		\centering
		\includegraphics[width=1.16\linewidth]{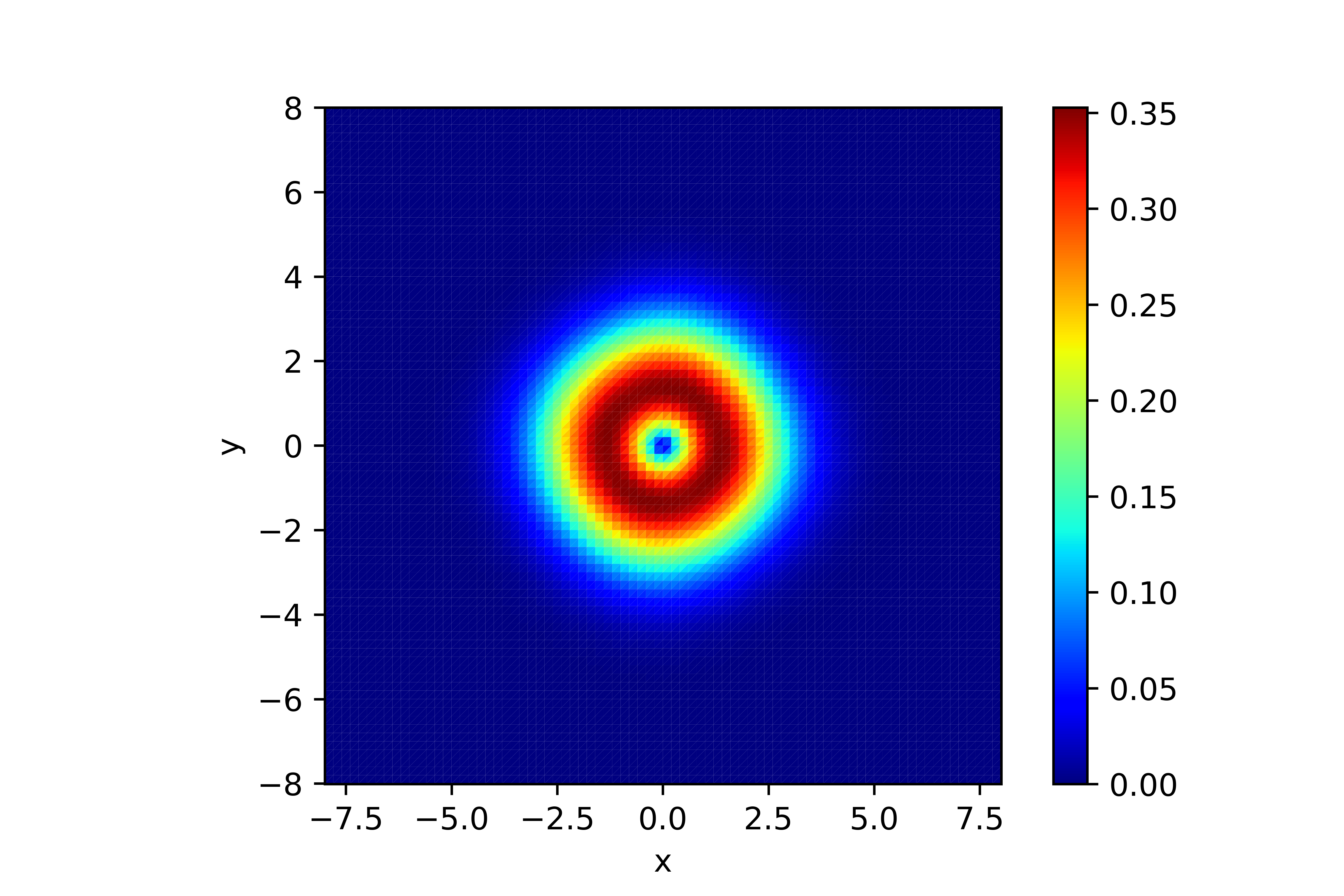}
		\caption{$t=5$}
	\end{subfigure}
	\centering
	\begin{subfigure}{0.325\linewidth}
		\centering
		\includegraphics[width=1.16\linewidth]{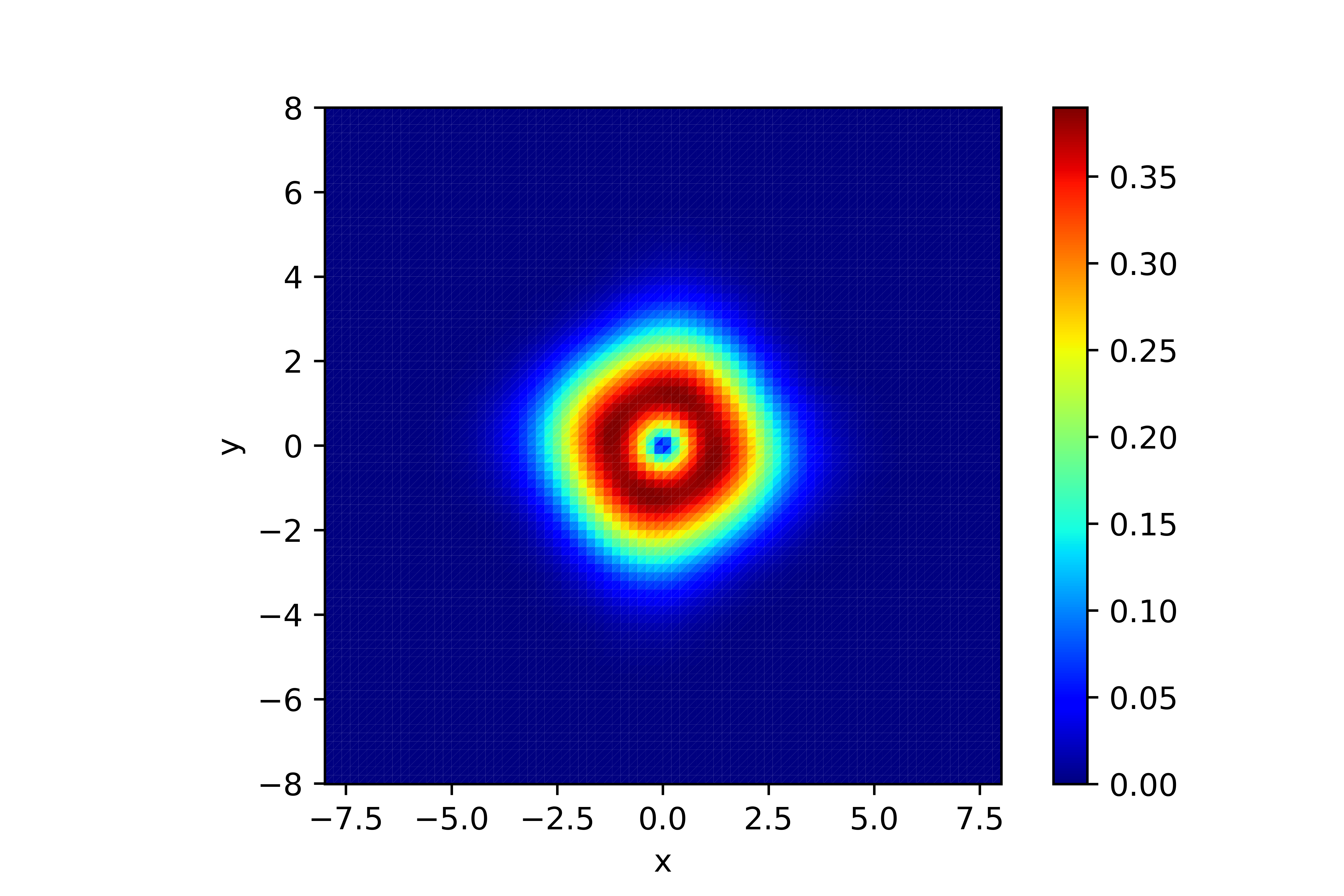}
		\caption{$t=10$}
	\end{subfigure}
	\caption{The patterns evolution of the wave function $|u(x_{1},x_{2},t)|$ with $V(x_{1},x_{2})=V_{1}(x_{1},x_{2})$.}
	\label{Case3Solution}
\end{figure}

\begin{figure}[htbp]
	\centering
	\begin{subfigure}{0.325\linewidth}
		\centering
		\includegraphics[width=1.16\linewidth]{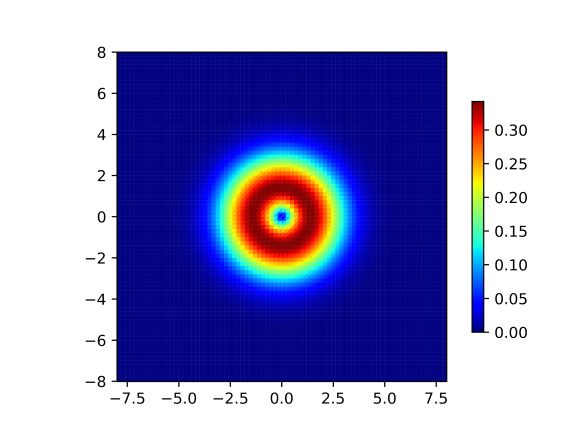}
		\caption{$t=0$}
	\end{subfigure}
	\centering
	\begin{subfigure}{0.325\linewidth}
		\centering
		\includegraphics[width=1.16\linewidth]{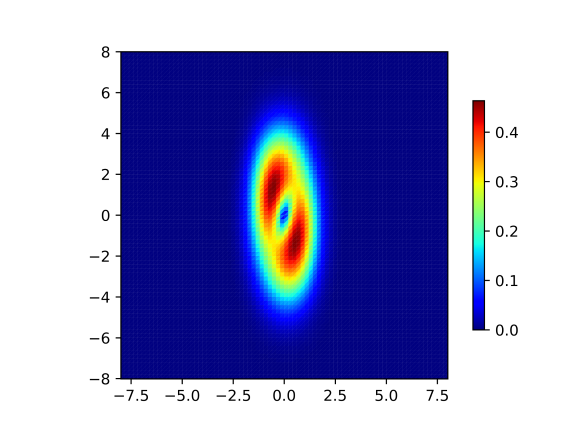}
		\caption{$t=1$}
	\end{subfigure}
	\centering
	\begin{subfigure}{0.325\linewidth}
		\centering
		\includegraphics[width=1.16\linewidth]{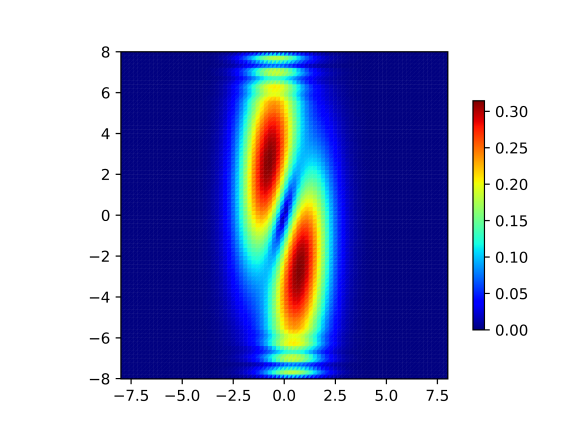}
		\caption{$t=2$}
	\end{subfigure}
	\begin{subfigure}{0.325\linewidth}
		\centering
		\includegraphics[width=1.16\linewidth]{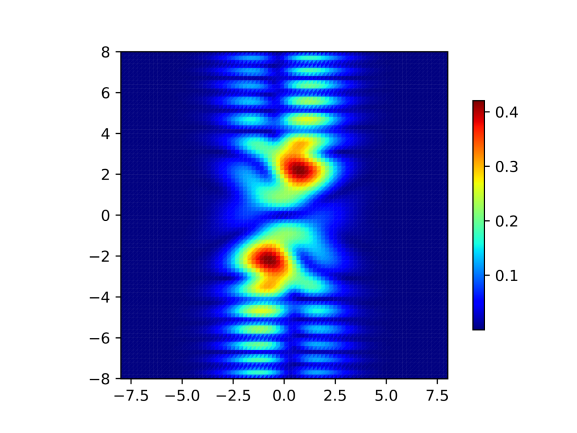}
		\caption{$t=5$}
	\end{subfigure}
	\centering
	\begin{subfigure}{0.325\linewidth}
		\centering
		\includegraphics[width=1.16\linewidth]{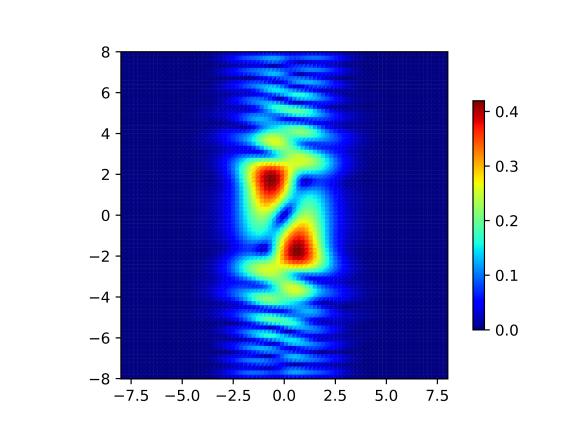}
		\caption{$t=7.5$}
	\end{subfigure}
	\centering
	\begin{subfigure}{0.325\linewidth}
		\centering
		\includegraphics[width=1.16\linewidth]{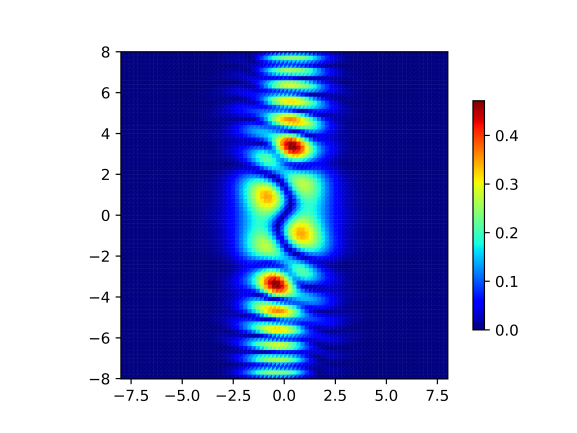}
		\caption{$t=10$}
	\end{subfigure}
	\caption{The patterns evolution of the wave function $|u(x_{1},x_{2},t)|$ with $V(x_{1},x_{2})=V_{2}(x_{1},x_{2})$.}
	\label{Case1Solutionb}
\end{figure}

{\bf  Test case 4.}
We present the evolution of the solution in \Cref{Case4Solution} - \Cref{Case1Solutionb} for the external potentials $V(x_{1},x_{2})=V_{i}(x_{1},x_{2})$, $i=0, 1,2$, respectively, using a mesh with $NC=80$ and a time step of $\tau=1\times 10^{-3}$, based on $\mathbb{Q}^2$ polynomials.

We first conduct numerical tests for the case with a zero potential, i.e., $V(x_{1},x_{2})=V_{0}(x_{1},x_{2})$. 
\Cref{Case4Solution} shows the patterns of the wave function $|u(x,y,t)|$ at time $t=0, 5, 10$, from which we can find that the pattern of the initial solution has evolved but not significantly, and the pattern evolves around the center of the pattern.

Next, we introduce different external potentials under the same conditions to observe the resulting changes in the solution.  This allows us to evaluate the performance of the proposed numerical method by comparing our results with similar findings in the literature.

We present the evolution of the solution in \Cref{Case3Solution} with potential $V(x_{1},x_{2}) = V_{1}(x_{1},x_{2})$ at times $t=0$, $t=5$, and $t=10$.
With the external potential $V_1$, the solution exhibits a pattern similar to that seen with zero potential. Notably, the pattern with $V_0$ at $t=10$ (see \Cref{Case4Solution}(c)) and the pattern with $V_1$ at $t=5$ (see \Cref{Case3Solution}(b)) are quite similar. This suggests that the external potential $V_1$ accelerates the evolution of patterns, particularly around the center of the pattern, compared to the zero potential case. Additionally, similar patterns of evolution to those in \Cref{Case3Solution} were also observed in \cite{Wang2018splitting}.

We also introduce a different external potential $V(x_{1},x_{2})=V_2(x_{1},x_{2})$ for problem \eqref{SPex}. The evolution of patterns is presented in \Cref{Case1Solutionb} at different times from $t=0$ to $t=10$. Under the influence of the external potential $V_2$, the patterns are driven away from the center, and similar patterns were also observed in \cite{Yi2022mass}.

{\bf Test case 5.} We compare the performance of the proposed relaxation Crank-Nicolson finite element algorithm \eqref{fullyVdiscrete}, or \Cref{rCNalg}, with the iterative method (IM) from \cite{Yi2022mass} by solving the Schr\"{o}dinger-Poisson problem \eqref{SPex} with $V(x_{1},x_{2})=V_{2}(x_{1},x_{2})$. 

\begin{figure}[htbp]
	\centering
	\begin{subfigure}{0.4\linewidth}
			\centering
			\includegraphics[width=1\linewidth]{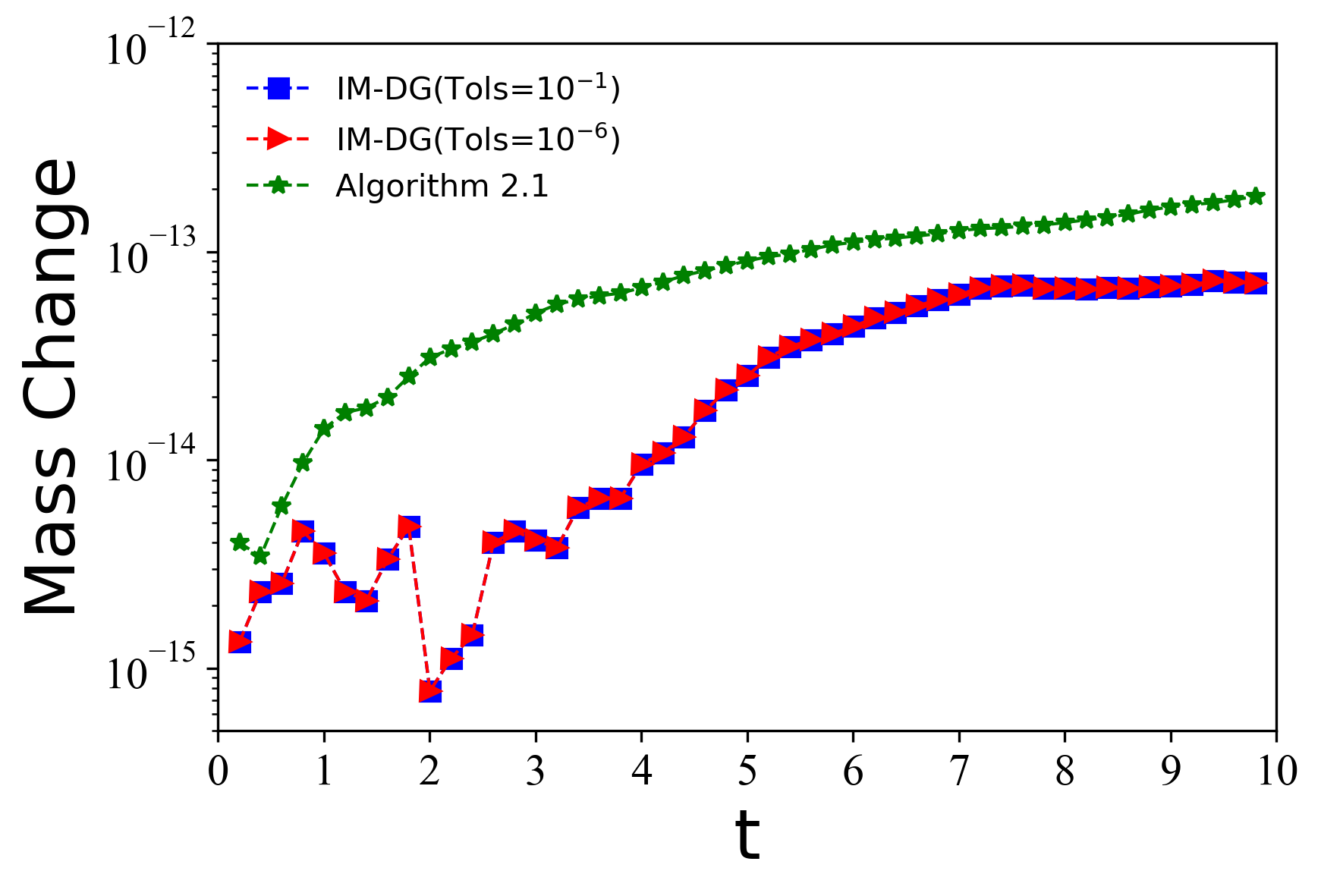}
		\end{subfigure}
	\centering 
	\begin{subfigure}{0.4\linewidth}
			\centering
			\includegraphics[width=1\linewidth]{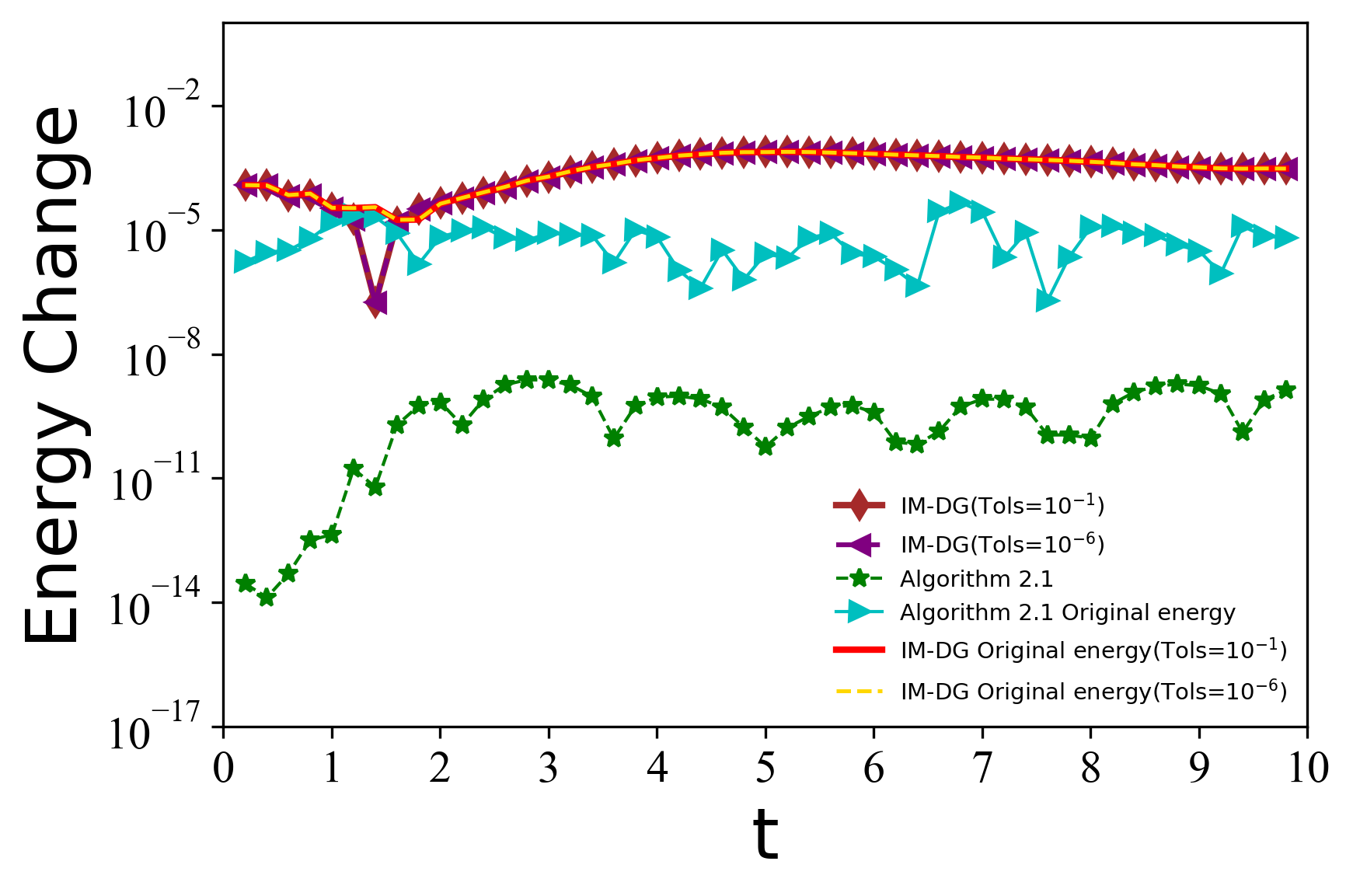}
		\end{subfigure}
	\caption{The patterns evolution of the mass and energy for \Cref{rCNalg} and IM-DG with $\tau=0.001$.}
	\label{Test5aSolution}
\end{figure}

\begin{figure}[htbp]
\centering
\begin{subfigure}{0.325\linewidth}
	\centering
	\includegraphics[width=1.16\linewidth]{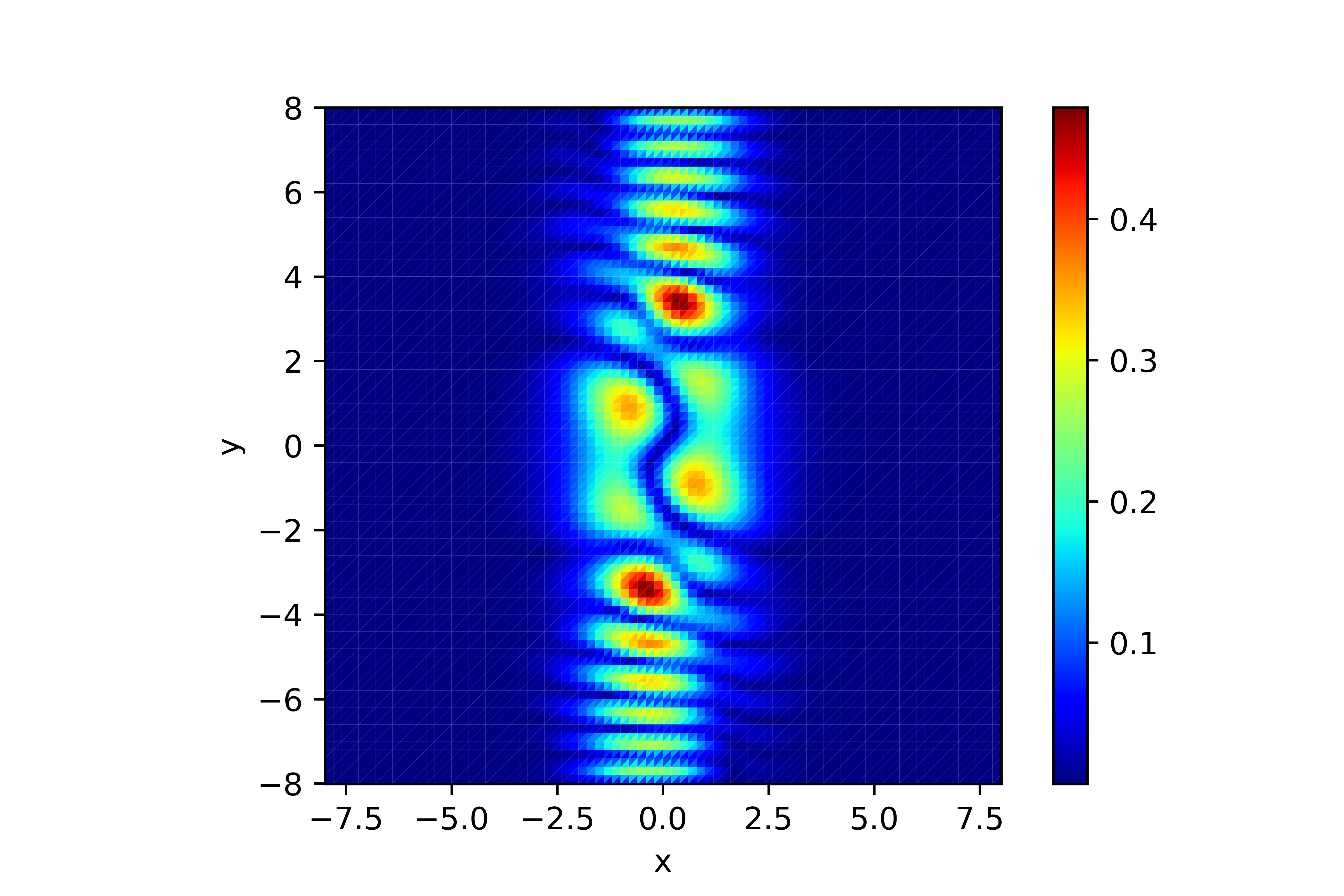}
	\caption{IM-DG (Tols $=10^{-1}$)}
\end{subfigure}
\centering
\begin{subfigure}{0.325\linewidth}
	\centering
	\includegraphics[width=1.16\linewidth]{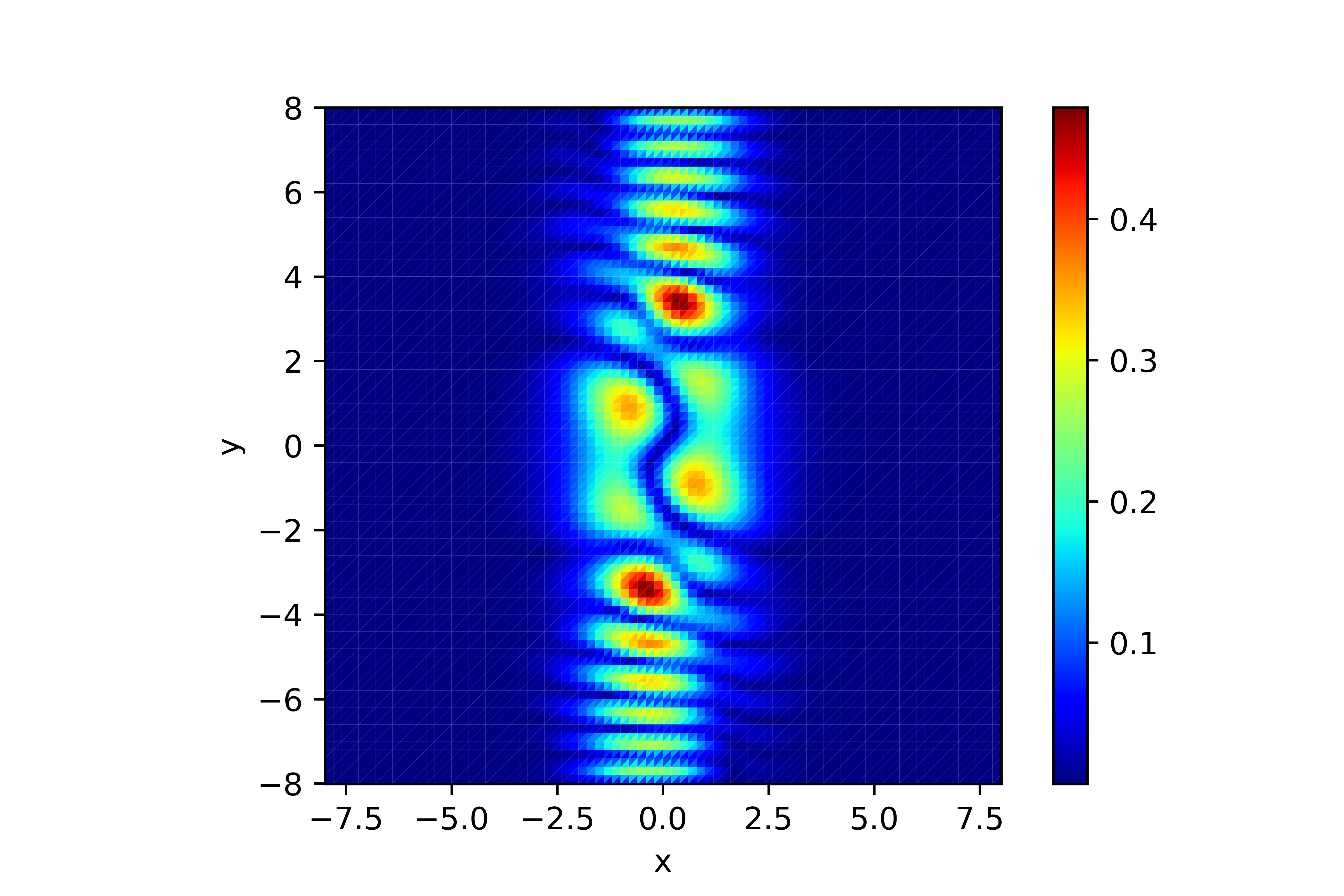}
	\caption{IM-DG (Tols $=10^{-6}$)}
\end{subfigure}
\centering
\begin{subfigure}{0.325\linewidth}
    \centering
   \includegraphics[width=1.16\linewidth]{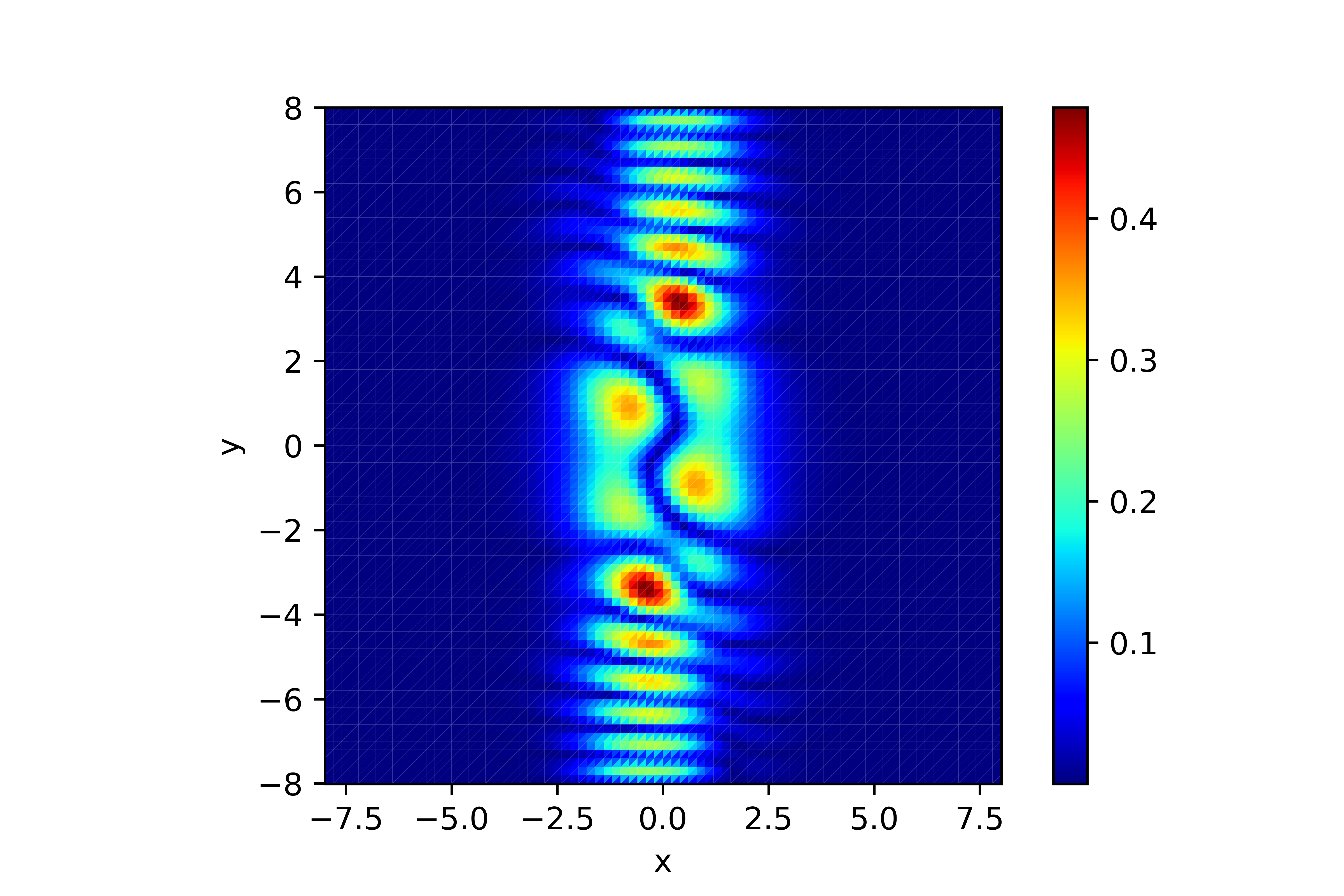}
    \caption{\Cref{rCNalg}}
\end{subfigure}
\vskip\baselineskip
\begin{subfigure}{0.325\linewidth}
    \centering
    \includegraphics[width=1.16\linewidth]{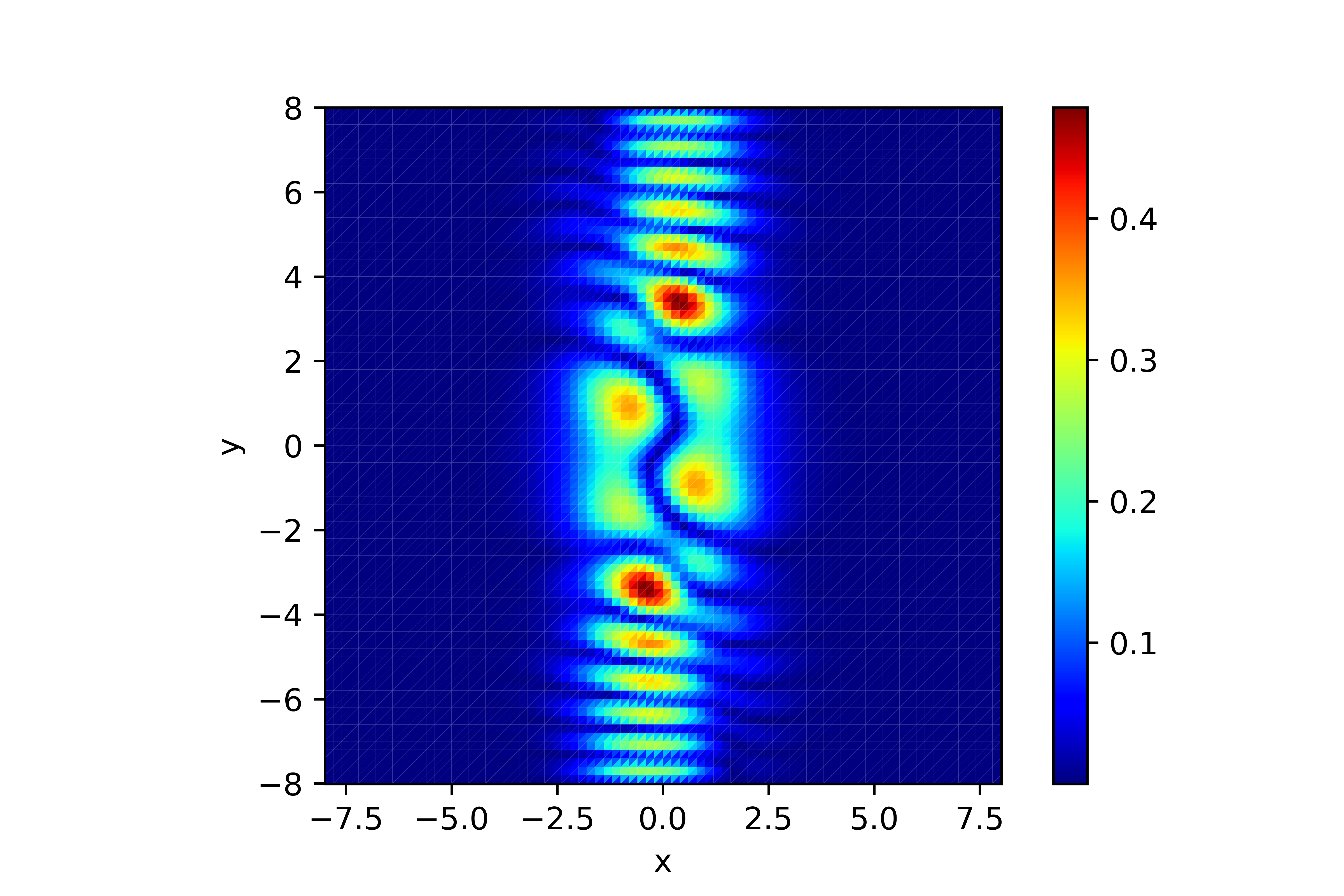}
    \caption{IM-FEM (Two-step iteration)}
\end{subfigure}
\centering
\begin{subfigure}{0.325\linewidth}
    \centering
    \includegraphics[width=1.16\linewidth]{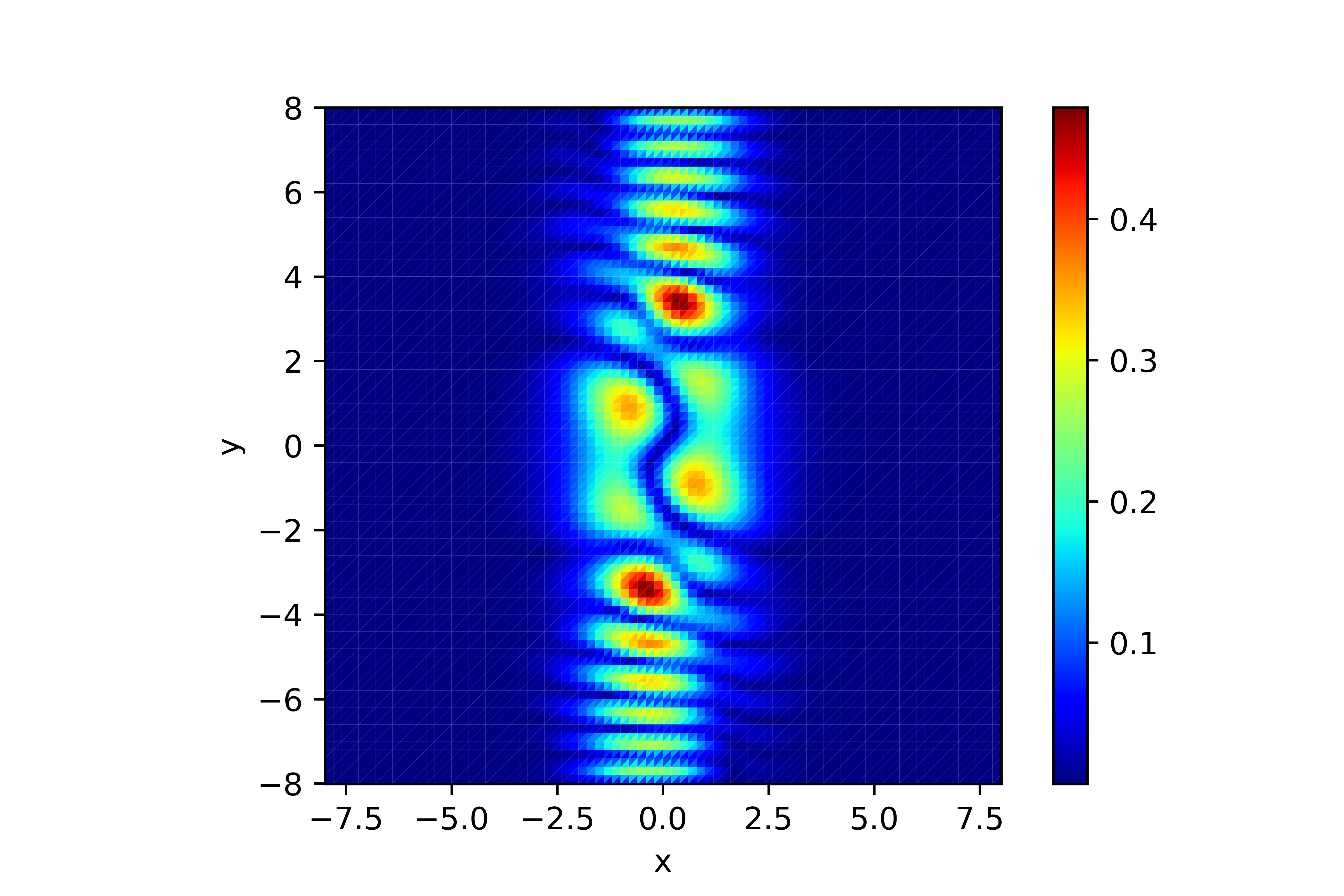}
    \caption{IM-FEM (Tols $=10^{-1}$)}
\end{subfigure}
    \centering
    \begin{subfigure}{0.325\linewidth}
   \centering
    \includegraphics[width=1.16\linewidth]{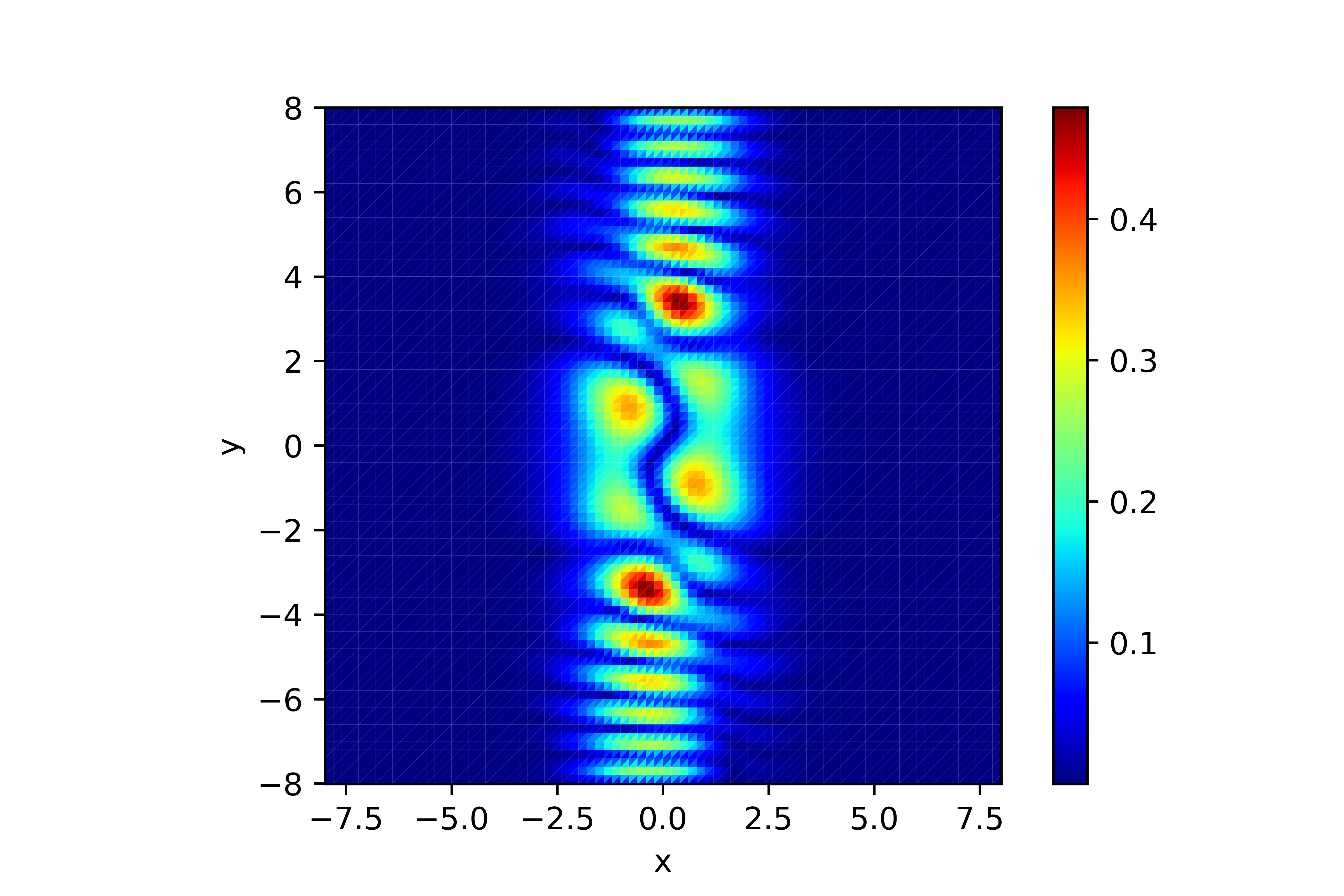}
    \caption{IM-FEM (Tols $=10^{-6}$)}
\end{subfigure}
\caption{The patterns of the wave function $|u(x_{1},x_{2},t)|$ at $t=10$ with $V_{2}(x_{1},x_{2})$.}
\label{Othermethod}
\end{figure}

First, we compare the performance of \Cref{rCNalg} with that of the IM using DG discretization (IM-DG) from \cite{Yi2022mass}. The parameters are set as follows: time step $\tau = 0.001$, mesh size $NC \times NC = 80 \times 80$, and $\mathbb{Q}^2$ polynomials. In the DG discretization, the penalty parameters are $\beta_{0} = 10$ and $\beta_{1} = {1}/{12}$. For the IM-DG method, the iteration is terminated when the prescribed tolerance (Tols $=10^{-1}$ or $10^{-6}$) is reached.
The corresponding solution patterns at $t = 10$ are shown in \Cref{Othermethod}(a)-(b), and they are comparable to that of \Cref{rCNalg} as shown in \Cref{Case1Solutionb}(f).
The corresponding CPU times of \Cref{rCNalg} and IM-DG is reported in \Cref{TableComTimeb}, showing that \Cref{rCNalg} is significantly more efficient, while the IM-DG method requires substantially more computational time. The evolution of mass, modified energy, and original energy is presented in \Cref{Test5aSolution}. Both methods conserve mass well. \Cref{rCNalg} preserves the modified energy with high accuracy, and the original energy is also conserved, though with a slightly larger error. In contrast, the IM-DG method exhibits noticeably larger relative errors in conserving both the modified and original energies compared to \Cref{rCNalg}.

\begin{table}[!ht]
	\caption{The computational time with $T=10$ and $V(x_{1},x_{2})=V_{2}(x_{1},x_{2})$.}
	\centering
	\begin{tabular}[c]{ccccccc}
		\toprule
		\Cref{rCNalg}   & IM-DG & IM-DG &  \\
		(linear, no iteration) & (Tols $=10^{-1}$)  & (Tols $=10^{-6}$)  \\
		\midrule
		 33411.52s &  195226.01s &  407897.47s \\
		\bottomrule
	\end{tabular}
	\label{TableComTimeb}
\end{table}

\begin{figure}[htbp]
	\centering
	\begin{subfigure}{0.4\linewidth}
		\centering
		\includegraphics[width=1\linewidth]{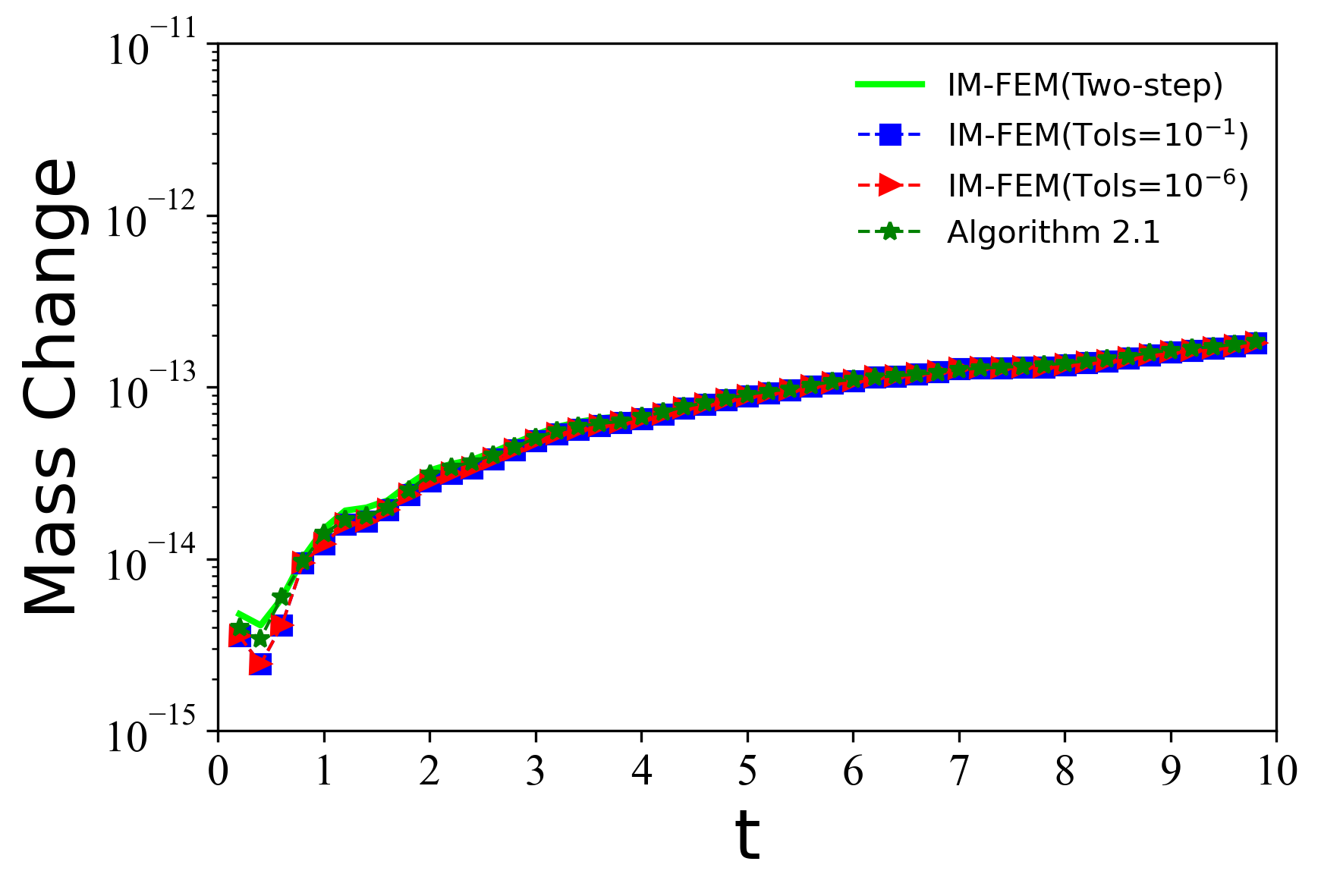}
	\end{subfigure}
	\centering
	\begin{subfigure}{0.4\linewidth}
		\centering
		\includegraphics[width=1\linewidth]{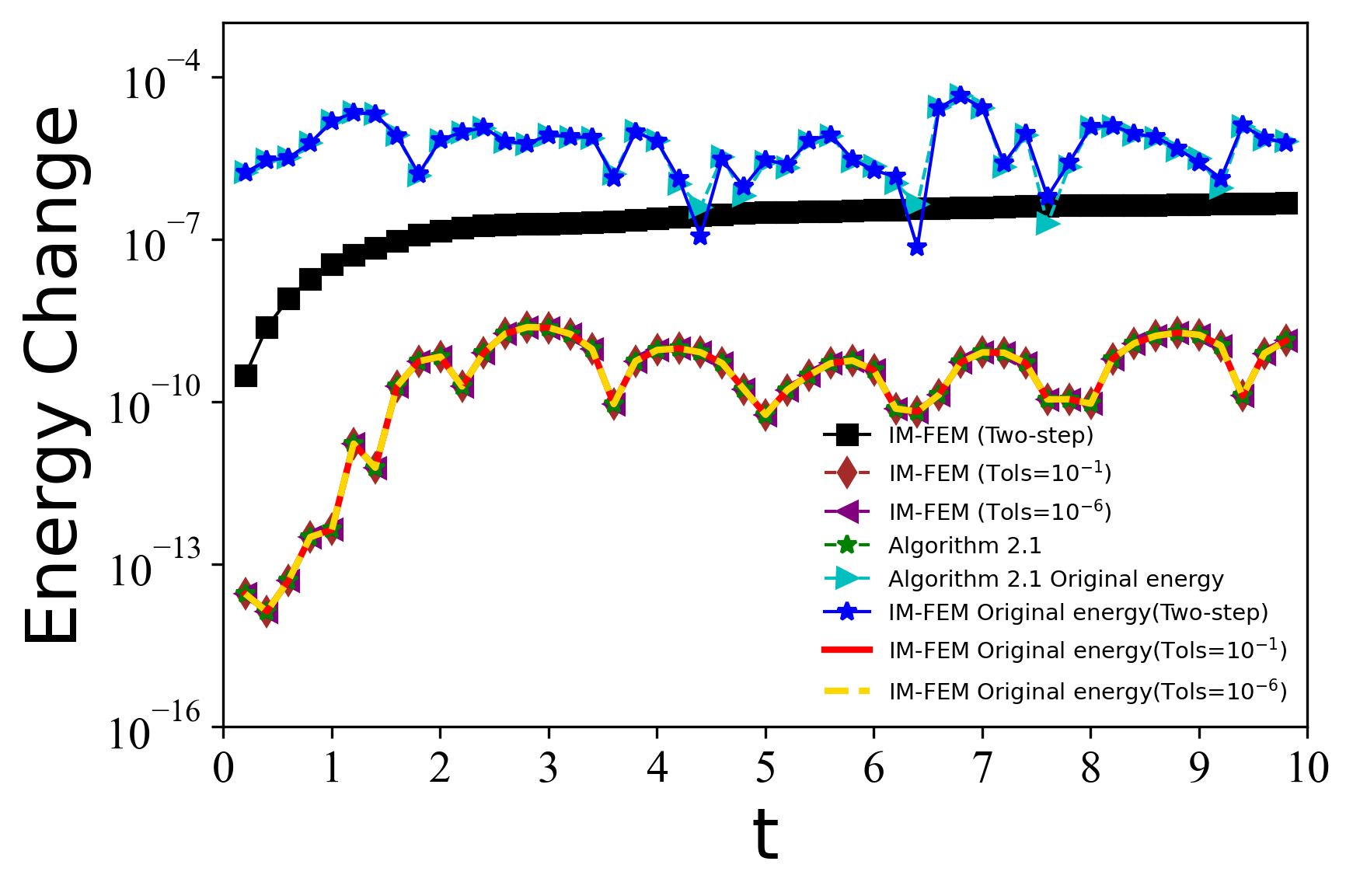}
	\end{subfigure}
	\caption{The patterns evolution of the mass and energy for \Cref{rCNalg} and IM-FEM with $\tau=0.001$.}
	\label{Test5Solutionb}
\end{figure}

\begin{figure}[htbp]
	\centering
	\begin{subfigure}{0.4\linewidth}
			\centering
			\includegraphics[width=1\linewidth]{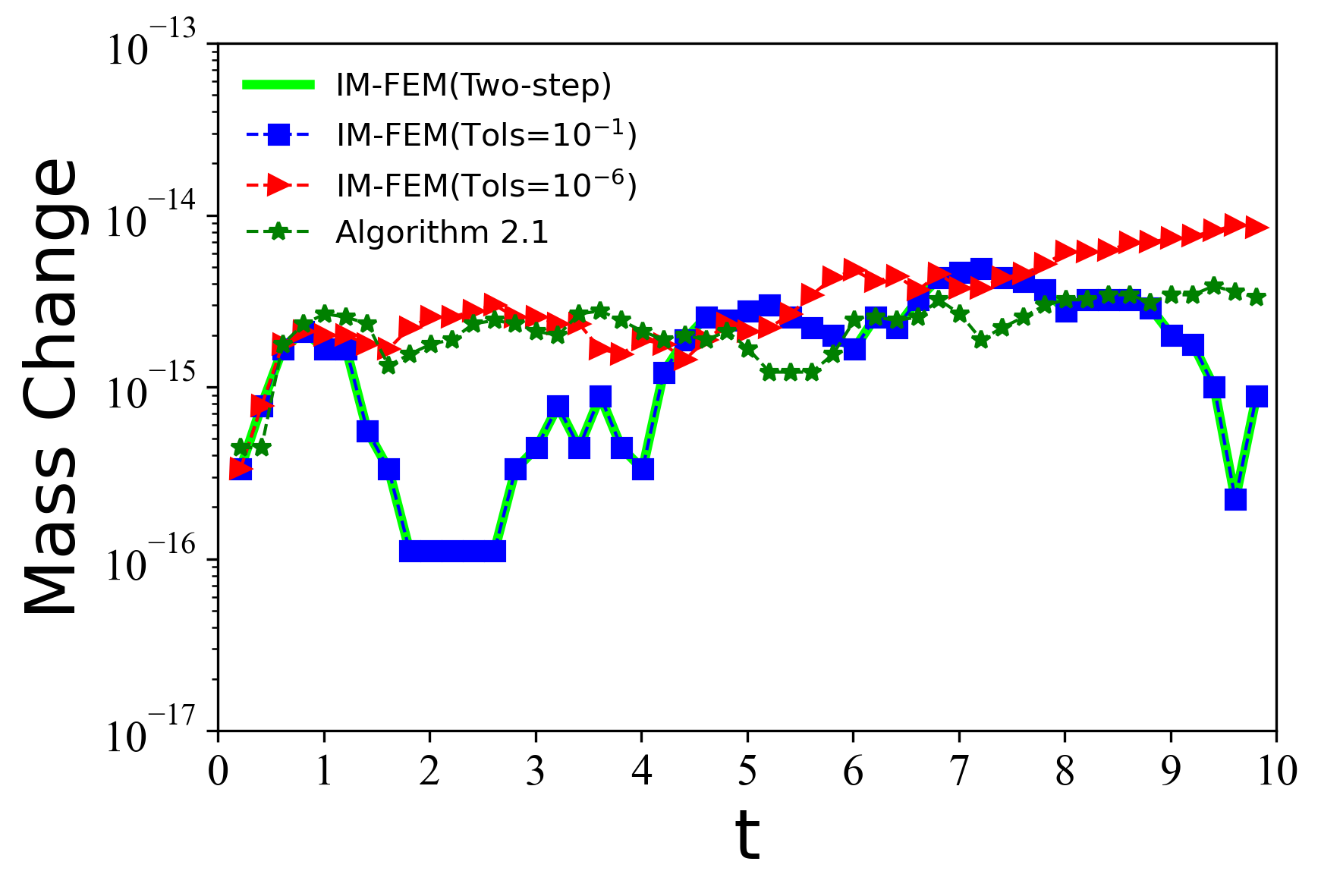}
		\end{subfigure}
	\centering
	\begin{subfigure}{0.4\linewidth}
			\centering
			\includegraphics[width=1\linewidth]{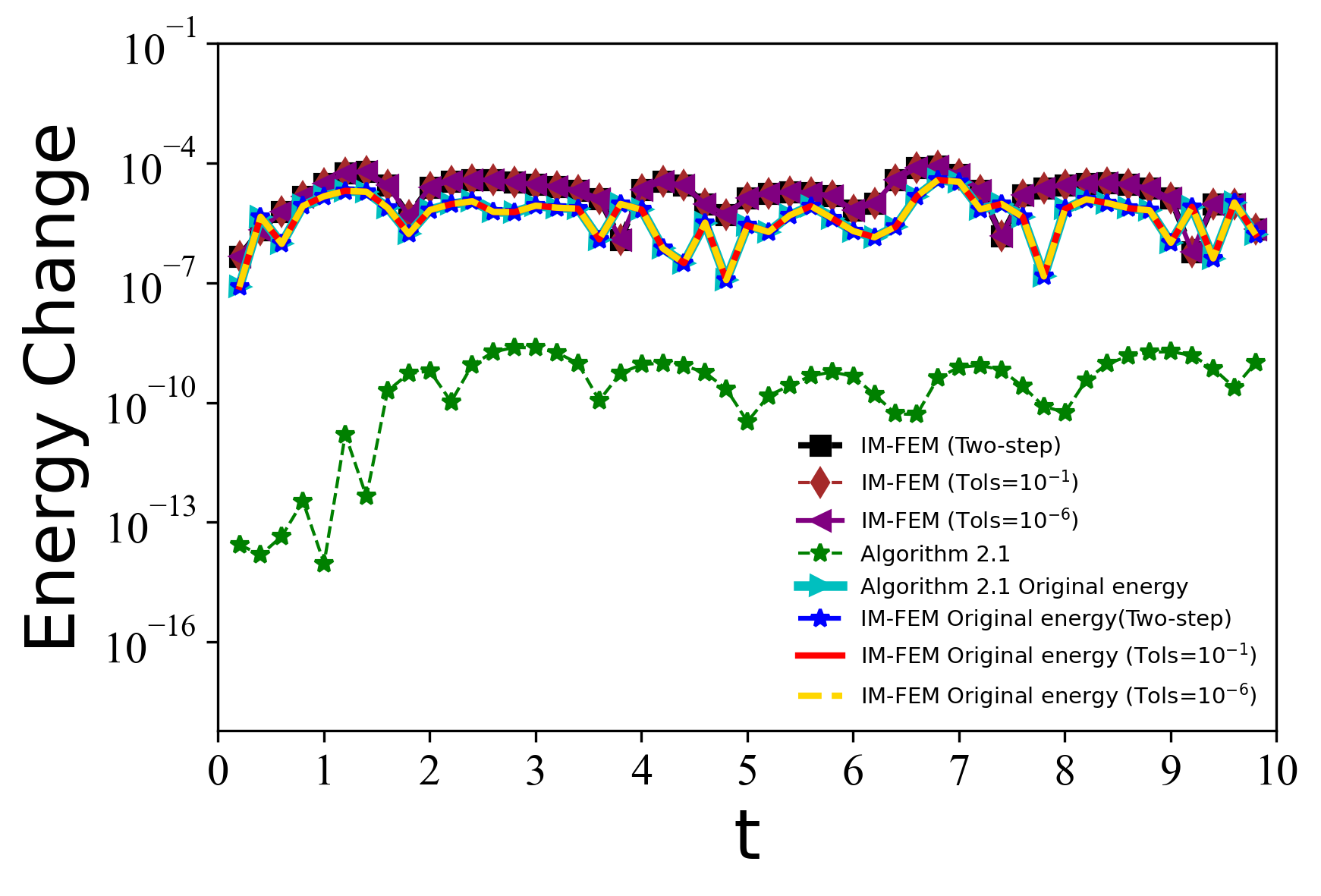}
		\end{subfigure}
		\caption{The patterns evolution of the mass and energy for \Cref{rCNalg} and IM-FEM with $\tau=0.01$.}
	\label{Test5Solutiona}
\end{figure}

Second, to eliminate the influence of the DG discretization and to provide a fairer comparison with \Cref{rCNalg}, we consider a IM-FEM variant, obtained by replacing the DG discretization in IM-DG from \cite{Yi2022mass} with the finite element method (FEM). This modification allows a larger time step for IM to produce a comparable final pattern. Specifically, we consider the time step both $\tau = 0.001$ and $\tau = 0.01$, mesh size $NC \times NC = 80 \times 80$, and employ $\mathbb{Q}^2$ polynomials for both \Cref{rCNalg} and IM-FEM.
In IM-FEM, the iteration is terminated either after two fixed steps or upon reaching the prescribed tolerance (Tols $=10^{-1}$ or $10^{-6}$). The solution patterns at $t = 10$ with time step $\tau=0.01$ are presented in \Cref{Othermethod}(c)–(f), and they are comparable to the pattern obtained by \Cref{rCNalg} with time step $\tau=0.001$, as shown in \Cref{Case1Solutionb}(f). The corresponding CPU times, reported in \Cref{TableComTimea} and \Cref{TableComTime}, indicate that \Cref{rCNalg} is the most efficient, while IM-FEM requires at least twice as much CPU time of \Cref{rCNalg}.
The evolution of mass, modified energy, and original energy is presented in \Cref{Test5Solutionb} and \Cref{Test5Solutiona}. Both methods conserve mass well. 
\Cref{rCNalg} preserves the modified energy with high accuracy, while the original energy is also conserved, albeit with slightly larger errors. In contrast, \Cref{Test5Solutionb} and \Cref{Test5Solutiona} demonstrate that IM-FEM requires a smaller time step and smaller iteration tolerance to preserve its modified energy, and it exhibits larger relative errors in conserving the original energy.



\begin{table}[!ht]
	\caption{{{The computational time at $T=10$ with $\tau=0.001$ and $V(x_{1},x_{2})=V_{2}(x_{1},x_{2})$.}}}
	\centering
	\begin{tabular}[c]{ccccccc}
		\toprule
		\Cref{rCNalg} & IM-FEM &   IM-FEM &  IM-FEM  &  \\
		(linear, no iteration)  & (two-step iteration)  & (Tols $=10^{-1}$)  & (Tols $=10^{-6}$)  \\
		\midrule
		 33411.52s & 64815.95s & 91914.44s  & 92866.31s \\
		\bottomrule
	\end{tabular}
	\label{TableComTimea}
\end{table}

\begin{table}[!ht]
\caption{The computational time at $T=10$ with $\tau = 0.01$ and $V(x_{1},x_{2})=V_{2}(x_{1},x_{2})$.}
\centering
\begin{tabular}[c]{ccccccc}
    \toprule
    \Cref{rCNalg} & IM-FEM &   IM-FEM &  IM-FEM  &  \\
  (linear, no iteration)  & (two-step iteration)  & (Tols $=10^{-1}$)  & (Tols $=10^{-6}$)  \\
    \midrule
     3488.17s & 6397.19s & 9334.41s  &  11717.49s\\
    \bottomrule
\end{tabular}
\label{TableComTime}
\end{table}


\section{Concluding remarks}\label{secConclusion}
A structure-preserving relaxation Crank-Nicolson finite element method has been proposed for the Schr\"{o}dinger-Poisson equation that contains the self-repulsion $|u|^2u$ in the Schr\"{o}dinger equation and the charge density $|u|^2$ in the Poisson equation, relying on a decoupled system that is equivalent to the original equation.
The fully discrete scheme is linear and is easy to implement without resorting to any iteration method. In addition, the finite element approximation is demonstrated to be both mass and modified energy conservative, irrespective of the mesh and time step. Optimal $L^{2}$ error estimates are established for the fully discrete scheme with second order accuracy in time and $(k+1)$th accuracy in space.
Numerical tests have been presented to verify the effectiveness and robustness of the proposed method.
The proposed relaxation Crank-Nicolson finite element method is a very competitive algorithm for solving the Schr\"{o}dinger-Poisson equation.

The spatial discretization utilized in this paper is the finite element method, it is noteworthy that the DG method \cite{Yi2022mass} can also be a viable alternative, in which the Poisson equation can be solved by the direct DG (DDG) method \cite{YHL14, YHL18}.
The proposed scheme preserves mass and a modified energy. Developing efficient numerical methods that preserve the original energy remains an important and challenging problem, which we leave for future work.
In the case of the three-dimensional Schr\"{o}dinger-Poisson equation, the self-repulsion term is substituted by $|u|^{4/3}u$. Extending the current findings to encompass this scenario could be an intriguing direction for future research, which we intend to pursue.



\section*{Acknowledgments}
Yi's research was partially supported by NSFC Project (12431014). Yin’s research was supported by the University of Texas at El Paso Startup Award.

\end{document}